\documentclass[a4paper]{article}

\usepackage{a4}
\usepackage{amsmath, amssymb}   % the inevitable ams packages for beautiful math
\usepackage{mathtools}
\usepackage{graphicx,color}                   % graphicx package for graphics includes
\usepackage{float}  
\usepackage{xspace}         
\usepackage[scientific-notation=true]{siunitx}         

\usepackage{tikz}
\tikzstyle{startstop} = [rectangle, rounded corners,minimum width=3cm, minimum height=1cm, text centered, text width=5cm, draw=black]
\tikzstyle{process} = [rectangle, minimum width=3cm, minimum height=1cm, text width=6 cm,text centered, draw=black]
\tikzstyle{decision} = [diamond, minimum width=3cm, minimum height=1cm, text width=3cm, text centered, draw=black]
\tikzstyle{arrow} = [thick,->,>=stealth]
\usetikzlibrary{shapes.geometric, arrows}
\usetikzlibrary{positioning} 
\usetikzlibrary{calc}
\usepackage{pgfplots}
\pgfplotsset{width=10cm,compat=1.9}
\usepackage[a4paper, total={6.5in, 8.5in}]{geometry}

\usepackage[justification=centering]{caption}

\usepackage{hyperref}
\usepackage[all]{hypcap}

\usetikzlibrary{spy}
\usetikzlibrary{shapes.geometric, arrows}
\newcommand{\bb}{{\boldsymbol b}}

\newcommand{\bn}{\boldsymbol n} 
\newcommand{\bt}{\boldsymbol t}

\newcommand{\mL}{{\mathrm L}}
\newcommand{\mH}{{\mathrm H}}
\newcommand{\mW}{{\mathrm W}}
\newcommand{\mD}{{\mathrm D}}
\newcommand{\mN}{{\mathrm N}}
\newcommand{\AS}{{\mathrm {AS}}}

\definecolor{light_gray}{gray}{0.75}
\colorlet{light_blue}{blue!20}

\newcommand{\blist}{\begin{list}{}{\itemsep0.0ex\parsep0.1ex\topsep0.2ex\leftmargin1.6em\labelwidth1.3em}}

\newtheorem{lemma}{Lemma}
\newtheorem{theorem}[lemma]{Theorem}
\newtheorem{remark}[lemma]{Remark}
\definecolor{light_gray}{gray}{0.75}
\definecolor{lighter_gray}{gray}{0.5}
\colorlet{light_blue}{blue!20}
\definecolor{dark_green}{rgb}{0.0, 0.6, 0.0}
\definecolor{royal_blue}{rgb}{0.0, 0.22, 0.66}
\definecolor{salmon}{rgb}{1.0, 0.55, 0.41}
\definecolor{gold}{rgb}{0.8, 0.63, 0.21}
\definecolor{navy_blue}{rgb}{0.0, 0.0, 0.5}
\definecolor{crimson}{rgb}{0.79, 0.0, 0.09}
\definecolor{amethyst}{rgb}{0.6, 0.4, 0.8}
\definecolor{alizarin}{rgb}{0.82, 0.1, 0.26}
\definecolor{amaranth}{rgb}{0.9, 0.17, 0.31}
\definecolor{azure}{rgb}{0.0, 0.5, 1.0}
\definecolor{canaryyellow}{rgb}{0.82, 0.41, 0.12}
\definecolor{carrotorange}{rgb}{0.8, 0.33, 0.0}
\definecolor{cadmiumgreen}{rgb}{0.0, 0.42, 0.24}
\definecolor{copper}{rgb}{0.72, 0.45, 0.2}
\definecolor{aqua}{rgb}{0.5, 1.0, 0.83}
\definecolor{awesome}{rgb}{1.0, 0.13, 0.32}
\definecolor{candyapplered}{rgb}{1.0, 0.03, 0.0}
\definecolor{caribbeangreen}{rgb}{0.0, 0.8, 0.6}

\title{Adaptive Finite Elements with Algebraic Stabilization for Convection-Dominated Transport}
\author{Naveed Ahmed, \footnote{Center for Applied Mathematics and Bioinformatics
Department of Mathematics and Natural Sciences, Gulf University for Science and Technology, Mubarak Al-Abdullah Area/West Mishref, Kuwait, \texttt{ahmed.n@gust.edu.kw}},  Abhinav Jha\footnote{Department of Mathematics, Indian Institute of Technology Gandhinagar, Palaj, Gandhinagar, 382055, Gujarat, India \texttt{abhinav.jha@iitgn.ac.in}}}

\date{}

\begin{document}
\maketitle

\begin{abstract}
We present a numerical investigation of residual-based a posteriori error
estimation for finite element discretizations of convection--diffusion
equations stabilized by algebraic flux correction and related algebraic
stabilization techniques. In particular, we consider AFC schemes employing
the BJK and Monolithic Convex (MC) limiters and algebraically stabilized
methods including MUAS, SMUAS, and the BBK approach. The performance of the
estimators and limiters are studied on adaptively refined meshes for several
two-dimensional test problems, including boundary layers, interior layers,
and a nonlinear convection problem with solution-dependent transport.

The experiments assess accuracy, preservation of the discrete maximum
principle, adaptive mesh behaviour, and computational efficiency. The
results show that the interaction between stabilization and a posteriori
error estimation depends strongly on mesh alignment and on the character of
the convection field. In particular, for problems with moving or curved layers, the behaviour of the limiters differs significantly: strongly
upwind-biased limiters provide the most accurate solutions, while smoother
algebraic stabilizations lead to more efficient nonlinear iterations. The
study also indicates that residual-based estimators remain reliable for both
linear and nonlinear problems but may react to changes in limiter activation
during adaptive refinement.

Overall, the numerical results clarify the practical behaviour of several
widely used stabilization techniques within an adaptive framework and
highlight aspects that are not yet fully explained by the current theory,
particularly for nonlinear transport problems.
\end{abstract}

\textbf{Keywords}: steady-state convection-diffusion-reaction equations; algebraically stabilized finite element methods; a posteriori estimator; adaptive grid refinement
\section{Introduction}\label{sec:intro}
This paper is devoted to the study of convection--diffusion--reaction (CDR) equations, which arise in a wide range of applications including, transport phenomena and species concentration models. The governing partial differential equation is given by
\begin{alignat}{3}\label{eq:cdr_eqn}
-\varepsilon \Delta u + \bb \cdot \nabla u + c u &= f \qquad &&\text{in } \Omega, \nonumber \\
u &= u_{\mD} \qquad &&\text{on } \Gamma_{\mD}, \nonumber \\
\varepsilon \nabla u \cdot \bn &= u_{\mN} \qquad &&\text{on } \Gamma_{\mN}.
\end{alignat}

Here, $\varepsilon > 0$ is the diffusion coefficient, $\bb$ is the convection (or advective) field, $c$ is the reaction coefficient, and $f$ is a source or sink term. The functions $u_{\mD}$ and $u_{\mN}$ specify the Dirichlet and Neumann boundary conditions, respectively. The equation is posed on a polygonal domain $\Omega \subset \mathbb{R}^d$ with $d \geq 2$, whose boundary $\Gamma$ is partitioned into disjoint Dirichlet and Neumann components: $\Gamma = \Gamma_{\mD} \cup \Gamma_{\mN}$, with $\Gamma_{\mD} \cap \Gamma_{\mN} = \emptyset$.

In the convection-dominated regime, where $L \|\bb\|_{\infty, \Omega} \gg \varepsilon$ and $L$ is the characteristic length scale (e.g., the domain diameter), the solution of Eq.~\eqref{eq:cdr_eqn} typically develops sharp gradients in narrow regions of the domain—known as interior and boundary layers. Moreover, solutions to Eq.~\eqref{eq:cdr_eqn} are known to satisfy maximum principles. This motivates the development of numerical methods that not only capture these layers accurately but also preserve a discrete analog of the maximum principle, known as the Discrete Maximum Principle (DMP). Standard discretization techniques, such as the finite difference or finite element methods, often fail in both respects: they may introduce spurious under- and overshoots due to violation of the DMP, and they may inadequately resolve sharp layers. As a result, the numerical solutions may be physically unreliable and unsuitable for practical use, see \cite{RST08}.

Keeping these two properties in mind, stabilization techniques were developed. One of the most prominent stabilization schemes is the Streamline Upwind Petrov Galerkin (SUPG) scheme developed by Brooks and Hughes \cite{BH81}. They approximate these layers well but fail to satisfy the DMP; hence, under and overshoots can be observed in the vicinity of layers. Other linear schemes, such as Upwind\cite{HB79}, Local Projection Stabilization (LPS) \cite{Kn09},  Spurious Oscillations at Layers Diminishing (SOLD) \cite{JK07_1} have the same issues. We refer to \cite{RST08} for an overview of these methods. Nonlinear stabilization schemes are a small class of techniques that not only approximate these layers but also satisfy the DMP. The idea of non-linear stabilization can be traced back to the work of Zalesak \cite{Zal79}. However, it has gained popularity through the work of Kuzmin and co-authors \cite{Ku06, Ku07, KH15} and is referred to as \emph{Algebraic Flux Correction} (AFC) schemes. The numerical analysis of these schemes has been developed in \cite{BJK16}, and further improvements have been proposed in \cite{BJK17, Ku20, KKJ24}. Recently, a more generalized approach has been presented in \cite{JK21, Kn23} called \emph{Algebraic Stabilizations}.  Here the symmetric conditions of the limiters are dropped and an upwind type method is proposed. We refer to \cite{BJK22} for a detailed review of these methods.

Another way to resolve the layers is to use non-equidistant meshes. Generally, there are two approaches to using non-equidistant meshes: predetermined meshes, such as the Shishkin mesh \cite{MOS12}, or adaptive meshes generated by a posteriori error estimators \cite{BR78}. In this paper, we focus on a posteriori error estimators. Much work has been done to develop a posteriori error estimators for Eq.~\eqref{eq:cdr_eqn} and linear stabilization methods. One of the first numerical studies was done in \cite{Joh00}, where it was shown that none of the estimators was robust to $\varepsilon$. By robustness, we mean that the constants appearing in the estimators should be independent of convection dominance. A robust estimator was developed in \cite{Ver05}, where a dual-norm approach was used, and its generalization to linear stabilization schemes was developed in \cite{TV15}. In \cite{JN13}, a robust estimator was also developed for the SUPG scheme, but in the natural norm of the system, that is, the SUPG norm, but here the analysis relied on certain assumptions. Lastly, a posteriori estimators in different norms, such as $\mL^1, \mL^2$, and $\mL^{\infty}$ can be found in \cite{HGMF06, HFD08,HDF11,DFK23}, respectively.

Very limited work exists at the intersection of a posteriori error estimation and AFC schemes. The first notable contribution in this direction was made in~\cite{ABR17}, where a fully computable, albeit non-robust, estimator was developed for linearity-preserving algebraic stabilization techniques. Subsequently, a residual-based non-robust estimator, independent of the choice of flux limiters, was proposed in~\cite{Jha20} specifically for AFC schemes. This work also integrated the SUPG-based estimator from~\cite{JN13} within the AFC framework. Both studies focused on the energy norm of the underlying system.

Further developments were presented in~\cite{JJK23}, which extended the residual-based estimator to non-conforming meshes. Three algebraic stabilization schemes were examined: the AFC method with the Kuzmin limiter~\cite{Ku06}, the AFC method with the BJK limiter~\cite{BJK17}, and the Monotone Upwind-type Algebraically Stabilized (MUAS) method~\cite{JK21}. Since the residual estimator predates the MUAS formulation, the stabilization contribution in MUAS was omitted in the analysis. More recently,~\cite{Jha24} extended the analysis from~\cite{Jha20} to a more general algebraic stabilization framework that includes the MUAS method and its improved variant, the Symmetric Monotone Upwind-type Algebraically Stabilized (SMUAS) method, introduced in~\cite{Kn23}. While~\cite{Jha24} primarily focused on theoretical analysis and included only basic numerical experiments, it demonstrated promising performance for SMUAS in preliminary studies. However, that study did not include another algebraic stabilized method which is based on the local smoothness indicator known as the BBK method~\cite{BBK17}, which are incorporated in the present work.

The present work builds on these foundations by offering a more comprehensive numerical evaluation of a posteriori error estimators and various algebraic flux limiters. The main contributions of this study are as follows:
\begin{enumerate}
    \item[$\bullet$] A detailed numerical assessment in two dimensions of AFC and algebraic stabilization schemes for adaptively refined grids.
    \item[$\bullet$] Inclusion of numerical studies for the Monolithic Convex Limiting (MC) scheme~\cite{Ku20}, as well as the BBK method~\cite{BBK17, BJK25}, which were absent in prior analyses.
    \item[$\bullet$] Numerical experiments involving nonlinear problems. To the best of our knowledge, prior numerical studies on AFC schemes have been limited to linear examples; in this work, we present results for nonlinear cases as well.
\end{enumerate}

It is worth noting that previous numerical assessments of these methods have largely been conducted on uniformly refined meshes; see~\cite{JJ18, JJ19, Jha22, JKP22}. In contrast, the current work focuses on adaptively refined meshes, offering new insights into their performance in more practical computational settings.

The structure of the paper is as follows. In Section~\ref{sec:algebraic_stab}, we provide a unified overview of the algebraic stabilization schemes considered in this study. Section~\ref{sec:apost} presents an overview of the important theorems for the residual-based a posteriori error estimators. In Section~\ref{sec:numres}, we carry out detailed numerical experiments, focusing on adaptive mesh refinement and a range of relevant parameters. Finally, Section~\ref{sec:summary} summarizes the main findings.

\section{Algebraic Stabilization Schemes}\label{sec:algebraic_stab}
Throughout this paper standard notions for Sobolev spaces and their norms are used, \cite{Ada75}.  Let $\Omega\subset \mathbb{R}^d$, $d\geq 2$ be a measurable set, then the $\mL^2$-inner product is defined by $\left(\cdot, \cdot\right)$ and its duality pairing by $\langle \cdot,\cdot \rangle$. The norm (semi-norm) on $\mW^{m,p}(\Omega)$ is denoted by $\|\cdot\|_{m,p,\Omega}(|\cdot|_{m,p,\Omega})$ with the convention $\mW^{m,2}(\Omega)=\mH^m(\Omega)$ and $\|\cdot\|_{m,\Omega}=\|\cdot\|_{m,2,\Omega}$ (similarly $|\cdot|_{m,\Omega}=|\cdot|_{m,2,\Omega}$).

It is well-known that under the assumption
\begin{equation}\label{eq:cdr_cond}
c(x)-\frac{1}{2}\nabla \cdot \bb(x)\geq \sigma>0,
\end{equation}
Eq.~\eqref{eq:cdr_eqn} possess an unique weak solution $u\in \mathcal{C}\left(\overline{\Omega}\right)\cap \mH_{\mD}^1(\Omega)$ satisfying
$$
a(u,v)=\left\langle f,v\right\rangle+\left\langle g,v\right\rangle_{\Gamma_{\mN}}\qquad \forall v\in\mH_{0,\mD}^1(\Omega),
$$
with 
\begin{alignat*}{2}
a(u,v)&=\varepsilon \left( \nabla u,\nabla v\right)+\left(\bb\cdot \nabla u,v\right)+(cu,v),\\
\mH_{\mD}^1(\Omega)&=\left\lbrace v\in \mH^1(\Omega): v|_{\Gamma_{\mD}}=u_{\mD}\right\rbrace,\\
\mH_{0,\mD}^1(\Omega)&=\left\lbrace v\in \mH^1(\Omega): v|_{\Gamma_{\mD}}=0\right\rbrace,
\end{alignat*}
and $\langle\cdot,\cdot\rangle_{\Gamma_{\mN}}$ is the duality pairing restricted to the Neumann boundary, e.g., see \cite[Sec.~III.1.1]{RST08}.

The algebraic stabilizations schemes for Eq.~\eqref{eq:cdr_eqn} reads as (see \cite{JK21}): Find $u_h\in W_h\left( \subset \mathcal{C}(\overline{\Omega})\cap \mH_{\mD}^1(\Omega)\right)$ such that
\begin{equation}\label{eq:as_scheme}
a_{\AS}(u_h;u_h,v_h)=\left\langle f,v_h\right\rangle+\left\langle g,v_h\right\rangle_{\Gamma_{\mN}}\qquad \forall v_h\in V_h,
\end{equation}
with $a_{\AS}(\cdot;\cdot,\cdot):\mH_{\mD}^1(\Omega)\times\mH_{\mD}^1(\Omega)\times \mH_{0,\mD}^1(\Omega)\rightarrow \mathbb{R}$ such that
$$
a_{\AS}(u;v,w)=a(u,v)+b_h(u;v,w),
$$
where $W_h,V_h$ are linear finite-dimensional subspaces of $\mathcal{C}(\overline{\Omega})\cap \mH_{\mD}^1(\Omega)$ and $\mathcal{C}\left(\overline{\Omega}\right)\cap \mH_{0,\mD}^1(\Omega)$, respectively and $b_h(\cdot;\cdot,\cdot)$ is the nonlinear stabilization term.

For AFC schemes such as those presented in~\cite{BJK16, BJK17, BJKR18}, the stabilization term is defined as
\begin{equation}\label{eq:afc_stab}
b_h(u;v,w) = \sum_{i,j=1}^N \left(1 - \alpha_{ij}(u)\right) d_{ij} \left(v(x_j) - v(x_i)\right) w(x_i), \qquad \forall\, u,v,w \in \mathcal{C}(\overline{\Omega}),
\end{equation}
where $\mathbb{D} = \{ d_{ij} \}$ is the artificial diffusion matrix, given by
\begin{equation}\label{eq:diff_matrix}
d_{ij} = -\max \left\{ a_{ij},\, 0,\, a_{ji} \right\}, \quad i \neq j, \qquad 
d_{ii} = -\sum_{i \neq j} d_{ij}.
\end{equation}
Here, $a_{ij}$ are the entries of the stiffness matrix associated with Eq.~\eqref{eq:cdr_eqn}, $N$ denotes the total number of nodes, and $\alpha_{ij}(u) \in [0,1]$ are the \emph{symmetric} solution-dependent flux limiters. In this work, we consider two types of limiters: the \emph{BJK} limiter~\cite{BJK17}, and the \emph{Monolithic Convex} (MC) limiter~\cite{Ku20}.

For the algebraic stabilization schemes proposed in~\cite{JK21, Kn23}, the symmetry condition on the limiters is relaxed. The stabilization term is defined as
\begin{equation}\label{eq:as_stab}
b(u;v,w) = \sum_{i,j=1}^N b_{ij}(u)\, v(x_j)\, w(x_i) \qquad \forall\, u,v,w \in \mathcal{C}(\overline{\Omega}),
\end{equation}
where
\begin{equation*}
b_{ij}(u) = -\max \left\{ \left(1 - \alpha_{ij}(u)\right)a_{ij},\, 0,\, \left(1 - \alpha_{ji}(u)\right)a_{ji} \right\}, \quad i \neq j, \qquad
b_{ii}(u) = -\sum_{i \neq j} b_{ij}(u),
\end{equation*}
and $\alpha_{ij}(u) \in [0,1]$ are non-symmetric flux limiters. 

The BBK method presented in \cite{BBK17} can also be considered as an algebraic stabilized method with $b_{ij}(u)$ defined by
\begin{equation}\label{eq:as_bbk}
b_{ij}(u) = d_{ij}\gamma_{ij}(u), \qquad  i \neq j, \qquad
b_{ii}(u) = -\sum_{i \neq j} b_{ij}(u).
\end{equation}
The concrete values of $\gamma_{ij}(u)$ would be given later.

We observe that in both Eq.~\eqref{eq:afc_stab} and Eq.~\eqref{eq:as_stab}, the stabilization term is linear in its second and third arguments.

In all the methods considered in this paper, the computation of limiters follows a two-step approach. First, the Galerkin finite element discretization is applied, incorporating Neumann boundary conditions. Once the limiters are computed, the Dirichlet boundary conditions are then imposed in the standard manner.

In the following, we describe the flux limiters and algebraic stabilization schemes that are studied in this work.

\subsection{BJK Limiter}

This method was proposed in \cite{BJK17}, and it is shown that, for $\mathbb{P}_1$ elements, the method can be made linearity preserving.  

The computation of the limiter starts with a pre-processing step by setting
\[
a_{ji} := 0 \quad \mbox{if  } a_{ij} < 0,\quad i = 1, \ldots, M,\  j = M + 1, \ldots, N,
\]
where $M$ denote the number of non-Dirichlet degrees of freedom, compare \cite[Eq.~(2.4)]{BJK17}. For a real number $a$, denote $a^+=\max\{a,0\}$ and $a^-=\min\{a,0\}$ then, the computation proceeds as follows:
\begin{enumerate}
    \item Compute 
    \begin{equation*}
       P_i^+  = \sum_{j=1}^N \left( d_{ij}(u_j-u_i)\right)^+,\qquad
       P_i^-  = \sum_{j=1}^N \left( d_{ij}(u_j-u_i)\right)^-.
    \end{equation*}
    \item\label{it:BJK2} Compute 
       \begin{equation*}
       Q_i^+ = q_i\left(u_i-u_i^{\max}\right), \qquad Q_i^- = q_i\left(u_i-u_i^{\min}\right),
    \end{equation*}
    with
    \[
    u_i^{\max}  =  \max_{j\in N_i\cup\{i\}}u_j,\qquad
    u_i^{\min}   = \min_{j\in N_i\cup\{i\}}u_j,\qquad
    q_i  = \sum_{j\in N_i}\gamma_i d_{ij},
    \]
    where 
\begin{equation}\label{eq:Ni_def}
   N_i=\{j\in\{1,\dots,N\}\setminus\{i\}\,:\,\,
             a_{ij}\neq0\,\,\,\mbox{or}\,\,\,a_{ji}>0\}
\end{equation}
    and $\gamma_i$ is a positive constant.
For $\mathbb{P}_1$ elements, a value of $\gamma_i$ is used 
which guarantees the linearity preservation, see \cite[Eq.~6.5]{BJK17} for details.

    \item Compute
    \begin{equation*}
       R_i^+ = \min\left\{ 1,\frac{Q_i^+}{P_i^+} \right\}, \quad R_i^- =
\min\left\{ 1,\frac{Q_i^-}{P_i^-} \right\},\qquad i=1,\dots, M.
    \end{equation*}
If \(P_i^+\) or \(P_i^-\) is zero, one sets \(R_i^+=1\) or \(R_i^-=1\), respectively. The values for \(R_i^+\) and \(R_i^-\) are set to $1$ also for Dirichlet nodes. 
    \item Compute
    \[ \overline{\alpha}_{ij} = \begin{cases} 
        R_i^+ & \mbox{ if } d_{ij}(u_j-u_i)>0,\\
        1  & \mbox{ if } d_{ij}(u_j-u_i)=0, \\
        R_i^- & \mbox{ if } d_{ij}(u_j-u_i)<0,
      \end{cases}\quad i,j=1,\dots, N.
   \]
\item Finally, one sets
\begin{equation*}
\alpha_{ij} = \min \left\lbrace\overline{\alpha}_{ij}, \overline{\alpha}_{ji}\right\rbrace,\quad i,j=1,\dots,N.
\end{equation*}
\end{enumerate}
In \cite{BJK17}, it was shown that the method is DMP satisfying on arbitrary grids.  In \cite{JJ19}, it was observed that solving the nonlinear system of equations arising from the BJK limiter was harder than for other AFC methods.

\subsection{Monolithic Convex Limiter}
In \cite{Ku20}, an algebraic stabilization technique was introduced for hyperbolic conservation laws. This method is demonstrated to yield a well-defined steady-state limit for the nonlinear discrete problem, irrespective of the chosen time step. Consequently, it presents a viable approach for discretizing steady-state equations. Although \cite{Ku20} does not incorporate diffusion effects, \cite{Jha22} employed the monolithic convex (MC) limiter to solve a steady-state convection-diffusion equation, yielding encouraging outcomes. The computation of these limiters is described by the following formula:
\begin{equation}\label{MClimiter}
\alpha_{ij} =
	\begin{cases} 
	\min\left\lbrace 1, \min \left\lbrace \frac{2d_{ij} \left(\bar{u}_{ij}-u_i^{\max}\right)}{d_{ij}(u_j-u_i)}, \frac{2d_{ij} \left(u_j^{\min}-\bar{u}_{ji}\right)}{d_{ij}(u_j-u_i)} \right\rbrace \right\rbrace &\text{if } d_{ij}(u_j-u_i)>0, \\
	0 &\text{if } d_{ij}(u_j-u_i)=0, \\
	\min\left\lbrace 1, \min \left\lbrace \frac{2d_{ij} \left(\bar{u}_{ij}-u_i^{\min}\right)}{d_{ij}(u_j-u_i)}, \frac{2d_{ij} \left(u_j^{\max}-\bar{u}_{ji}\right)}{d_{ij}(u_j-u_i)}\right\rbrace \right\rbrace &\text{otherwise},
	\end{cases}
\end{equation}
where $\bar u_{ij}$ are intermediate states defined by
$$2d_{ij} \bar{u}_{ij}  =  d_{ij} ( u_i+u_j ) + a_{ij} (u_j-u_i)$$
and
\begin{align}
\label{eq:u_max_min}
u_i^{\max} = \max_{j\in S_i}u_j,\ \ u_i^{\min} = \min_{j\in S_i}u_j.
\end{align}
In the last formula, 
$S_i=\{j\in\{1,\ldots,N\}\,:\,a_{ij}\ne 0\}$ that is the integer set containing the indices of node $i$ and its nearest neighbors. Note that the definition of $u_i^{\max}$ and $u_i^{\min}$ can be
changed to ensure linearity preservation \cite[Section 6.1]{Ku20}. A local DMP for the algebraic system corresponding to the discretization of the steady-state (transport) problem was also proved in this paper. Numerical studies on uniformly refined grids were performed in \cite{JKP22} and it was observed that the efficiency of the MC limiter is better than the BJK limiter.

\subsection{MUAS Method}
Till now, we have considered symmetric limiters. For the algebraic stabilizations schemes proposed in \cite{JK21}, the symmetric condition is dropped on $\alpha_{ij}$. This is an upwind type method referred to as Monotone Upwind-Type Algebraically Stabilized (MUAS) and the limiters are computed as follows:
\begin{enumerate}
    \item Compute 
    \begin{equation*}
       P_i^+  =  \sum_{j=1, a_{ij}>0}^N a_{ij}(u_i-u_j)^+,  \quad
       P_i^-  =  \sum_{j=1, a_{ij}>0}^N  a_{ij}(u_i-u_j)^-.
    \end{equation*}
    \item Compute 
       \begin{equation*}
       Q_i^+ =\sum_{j=1}^N \max\left\{|a_{ij}|,a_{ji}\right\}(u_j-u_i)^+, \quad 
       Q_i^-=\sum_{j=1}^N \max\left\{|a_{ij}|,a_{ji}\right\}(u_j-u_i)^-.
    \end{equation*}
    \item Compute
    \begin{equation*}
       R_i^+ = \min\left\{ 1,\frac{Q_i^+}{P_i^+} \right\}, \quad R_i^- =
\min\left\{ 1,\frac{Q_i^-}{P_i^-} \right\},\qquad i=1,\dots, M.
    \end{equation*}
     If \(P_i^+\) or \(P_i^-\) is zero, one sets \(R_i^+=1\) or \(R_i^-=1\), respectively.  The values of \(R_i^+\) and \(R_i^-\) are set to $1$ for Dirichlet nodes as well. 
    \item Define
    \[ \alpha_{ij} = \begin{cases} 
        R_i^+ & \mbox{ if } u_i>u_j ,\\
        1 & \mbox{ if } u_i=u_j,\\
        R_i^- & \mbox{ if } u_i<u_j,
      \end{cases} \quad  i,j=1,\dots,N.
   \]
\end{enumerate}
Local and global DMP have been shown in \cite{JK21} to be valid on arbitrary simplicial meshes, but this method is not linearity preserving.

\subsection{SMUAS Method}
Recently, a further extension of the MUAS method has been proposed in \cite{Kn23} called the Symmetrized Monotone Upwind Type Algebraically Stabilized (SMUAS) method, which is an upwind type linearity-preserving limiter that satisfies DMP on arbitrary meshes. The limiters are computed as follows:
\begin{enumerate}
    \item Compute 
    \begin{alignat*}{2}
       P_i^+  &= \sum_{j\in S_i}|d_{ij}|\left\lbrace (u_i-u_j)^+ + (u_i-u_{ij})^+\right\rbrace,\\
       P_i^-  &=  \sum_{j\in S_i}|d_{ij}|\left\lbrace (u_i-u_j)^- + (u_i-u_{ij})^-\right\rbrace.
     \end{alignat*}
    \item Compute 
       \begin{alignat*}{2}
       Q_i^+ &= \sum_{j\in S_i}\max\lbrace |a_{ij}|, a_{ji}\rbrace\left\lbrace (u_j-u_i)^+ + (u_{ij}-u_i)^+\right\rbrace, \\
       Q_i^- &= \sum_{j\in S_i}\max\lbrace |a_{ij}|, a_{ji}\rbrace\left\lbrace (u_j-u_i)^- + (u_{ij}-u_i)^-\right\rbrace.
    \end{alignat*}
    \item Compute
    \begin{equation*}
       R_i^+ = \min\left\{ 1,\frac{Q_i^+}{P_i^+} \right\}, \quad R_i^- =
\min\left\{ 1,\frac{Q_i^-}{P_i^-} \right\},\qquad i=1,\dots, M.
    \end{equation*}
     If \(P_i^+\) or \(P_i^-\) is zero, one sets \(R_i^+=1\) or \(R_i^-=1\), respectively.  The values of \(R_i^+\) and \(R_i^-\) are set to $1$ for Dirichlet nodes as well. 
    \item Define
    \[ \alpha_{ij} = \begin{cases} 
        R_i^+ & \mbox{ if } u_i>u_j ,\\
        1 & \mbox{ if } u_i=u_j,\\
        R_i^- & \mbox{ if } u_i<u_j,
      \end{cases} \quad  i,j=1,\dots,N.
   \]
\end{enumerate}
In the above definition $u_{ij}=u_i + \nabla u_h|_{K_i^j}\cdot (x_i-x_j)$ for all $j\in S_i$. Here $K_i^j$ is a mesh cell containing $x_i$ that is intersected by the half line $\lbrace x_i+\theta (x_i-x_j):\theta >0\rbrace$. In \cite{Kn23}, optimal convergence rates were observed numerically, whereas the Kuzmin, BJK, and the MUAS methods did not have the optimal rates.

\begin{remark}
For the computation of the SMUAS method, it is necessary to evaluate $\nabla u_h|_{K_i^j}$ at each nonlinear iteration. This step is computationally expensive. To mitigate this issue, at the beginning of the nonlinear loop, we store the cell number and the coordinates of the reflected point $x_j$, and subsequently use these stored values to compute the gradient and, in turn, $u_{ij}$.
\end{remark}

\begin{remark}
In the SMUAS method, if $i$ is a Neumann node, then the coefficient $\alpha_{ij}$ is computed using the MUAS method.
\end{remark}

\subsection{BBK Method}\label{subsec:bbk_method}
The last algebraic stabilization method we consider is the BBK method, originally presented in \cite{BBK17}. In Eq.~\eqref{eq:as_bbk}, the factors $\gamma_{ij}(u)$ are computed as
$$
\gamma_{ij}(u)
= \left[\max\left\{ \alpha_i(u),\, \alpha_j(u) \right\} \right]^p,
\qquad i,j=1,2,\dots,N,
$$
where
$$
\alpha_i(u)=
\begin{cases}
\dfrac{\left| \sum_{j\in S_i}(u_i - u_j) \right|}{\sum_{j\in S_i} |u_i - u_j|}, & 
\text{if } \sum_{j\in S_i} |u_i - u_j| \neq 0,\\[2ex]
0, & \text{otherwise},
\end{cases}
$$
and for Dirichlet nodes we set $\alpha_i(u)=0$ for $i = M+1,\dots, N$.  
The exponent $p$ controls the rate at which the artificial diffusion decays; following \cite{BJK25}, we choose $p=10$. A limited numerical study of the BBK method exists, and a comparison can be found in \cite{JKP22}.

\begin{remark}
The BBK method used here is a modification of the original BBK method proposed in \cite{BBK17} and later reformulated in \cite{BJK25}. This reformulation is necessary because the original method cannot be written in the standard AFC framework: it may be the case that $b_{ij}(u)=0$ even when $d_{ij}\neq 0$.
\end{remark}

\section{A Posteriori Error Estimators}\label{sec:apost}
This section reviews the a posteriori error estimators developed in \cite{Jha24}.

Before presenting the results, we provide some details about the mesh triangulation. Let $\{\mathcal{T}_h\}_{h > 0}$ denote a family of conforming triangulations of the domain $\Omega$, where each triangulation $\mathcal{T}_h$ consists of simplices and has characteristic mesh size 
\[
h := \max_{K \in \mathcal{T}_h} h_K.
\]
We assume that each triangulation is admissible, meaning that any two elements are either disjoint or share a complete $m$-dimensional face for some $0 \leq m \leq d - 1$. Furthermore, the mesh is assumed to be shape-regular.

The set of all edges is denoted by $\mathcal{E}_h$, and the set of edges associated with an individual element $K \in \mathcal{T}_h$ is denoted by $\mathcal{E}_h(K)$. The set of all faces is denoted by
\[
\mathcal{F}_h = \mathcal{F}_{h,\Omega} \cup \mathcal{F}_{h,\mD} \cup \mathcal{F}_{h,\mN},
\]
where $\mathcal{F}_{h,\Omega}$, $\mathcal{F}_{h,\mD}$, and $\mathcal{F}_{h,\mN}$ refer to the sets of interior, Dirichlet boundary, and Neumann boundary faces, respectively. In two dimensions, faces coincide with edges, and we have $\mathcal{F}_h = \mathcal{E}_h$.

Our analysis is carried out in the energy norm, which is given by
$$
\|u\|_a^2=\varepsilon\|\nabla u\|_{0,\Omega}^2+\sigma \|u\|_{0,\Omega}^2.
$$

\begin{theorem}\label{thm:global_upper_bound}
Let $u_h \in W_h$ be a solution of Eq.~\eqref{eq:as_scheme}. Then the global \emph{a posteriori} error estimate in the energy norm is given by
\begin{equation}\label{eq:post_upper_bound}
\|u - u_h\|_a^2 \leq \eta^2 := \eta_1^2 + \eta_2^2 + \eta_3^2,
\end{equation}
where
\begin{align*}
\eta_1^2 &:= \sum_{K \in \mathcal{T}_h} \min\left\{ \frac{4 C_\mathrm{I}^2}{\sigma},\ \frac{4 C_\mathrm{I}^2 h_K^2}{\varepsilon} \right\} \|R_K(u_h)\|_{0,K}^2, \\ 
\eta_2^2 &:= \sum_{F \in \mathcal{F}_h} \min \left\{ \frac{4 C_\mathrm{F}^2 h_F}{\varepsilon},\ \frac{4 C_\mathrm{F}^2}{\sigma^{1/2} \varepsilon^{1/2}} \right\} \|R_F(u_h)\|_{0,F}^2,\\
\eta_3^2 &:= \sum_{E \in \mathcal{E}_h} \min \left\{ \frac{4 \kappa_1 h_E^2}{\varepsilon},\ \frac{4 \kappa_2}{\sigma} \right\} |b_E|^2 h_E^{1-d} \|\nabla u_h \cdot \bt_E\|_{0,E}^2,
\end{align*}
and the element and face residuals are defined as follows:
\begin{align*}
R_K(u_h) &:= f + \varepsilon \Delta u_h - \bb \cdot \nabla u_h - c u_h \quad \text{in } K, \\
R_F(u_h) &:= 
\begin{cases}
- \varepsilon \left[\![ \nabla u_h \cdot \bn_F ]\!\right]_F & \text{if } F \in \mathcal{F}_{h,\Omega}, \\
u_{\mathrm{N}} - \varepsilon (\nabla u_h \cdot \bn_F) & \text{if } F \in \mathcal{F}_{h,\mN}, \\
0 & \text{if } F \in \mathcal{F}_{h,\mD},
\end{cases}
\end{align*}
where $\left[\![ \cdot ]\!\right]$ denotes the jump across the face $F$, and $\bn_F$ is the outward-pointing unit normal to $F$.
\end{theorem}

\begin{remark}
Several constants appear in Theorem~\ref{thm:global_upper_bound}, which we now describe. The constant $C_\mathrm{I}$ denotes the interpolation constant (see~\cite[Corollary~4.8.15]{BS08}), and $C_\mathrm{F}$ arises from the face residual estimates (see~\cite[Lemma~3.1]{Ver13}). The constants $\kappa_1$ and $\kappa_2$ are defined as
\[
\kappa_1 = C\,C_{\mathrm{edge,max}} \left(1 + \left(1 + C_{\mathrm{I}}\right)^2\right), \quad 
\kappa_2 = C\,C_{\mathrm{inv}}^2\,C_{\mathrm{edge,max}} \left(1 + \left(1 + C_{\mathrm{I}}\right)^2\right),
\]
where $C$ is a generic constant independent of the mesh size $h$, $C_{\mathrm{inv}}$ is the inverse inequality constant (see~\cite[Lemma~4.5.3]{BS08}), and $C_{\mathrm{edge,max}}$ is a computable constant introduced in~\cite[Remark~9]{Jha20}.

\end{remark}

\begin{remark}
We would like to highlight a few remarks regarding the a posteriori error bound. First, we note that the estimator is not robust to the parameter~$\varepsilon$, a well-known limitation of energy-norm estimates. An alternative is to employ duality-based techniques, as proposed in~\cite{TV15}; however, extending such approaches to nonlinear models presents significant challenges.

Next, we address the estimator's efficiency, specifically the local lower bound. A formal local lower bound for the estimator exists without requiring any assumptions on the stabilization term. It is termed \emph{formal} because the local error is not directly bounded by the residual as is standard in the literature, but rather by a term that decays at the optimal rate.
\end{remark}
\section{Numerical Studies}\label{sec:numres}
This section presents the numerical assessment of the a posteriori error estimator introduced in Sec.~\ref{sec:apost}. We compare different algebraic stabilization schemes from Sec.~\ref{sec:algebraic_stab} on adaptively refined grids using a set of performance metrics.

The general strategy for solving an a posteriori error problem follows the loop:
\begin{equation}
\textbf{SOLVE} \rightarrow \textbf{ESTIMATE} \rightarrow \textbf{MARK} \rightarrow \textbf{REFINE}.
\end{equation}
We now go over each step and provide relevant details.

\begin{enumerate}
    \item[$\bullet$] \textbf{SOLVE}: The solution procedure for the AFC scheme on uniform grids has been discussed in detail in \cite{JJ18, JJ19}. We provide a brief overview here. Let the system of equations corresponding to Eq.~\eqref{eq:as_scheme} be written as
    \begin{equation}\label{eq:mod_sys_0}
    \left(\mathbb{A} + \mathbb{B}(U)\right)U = F,
    \end{equation}
    where $\mathbb{A} = \lbrace a_{ij} \rbrace_{i,j=1}^N$ is the stiffness matrix obtained using $\mathbb{P}_1$ finite elements, $\mathbb{B}(U) = \lbrace b_{ij}(U) \rbrace_{i,j=1}^N$ is the stabilization matrix, and $F = \lbrace f_i \rbrace_{i=1}^N$ is the discretized right-hand side. 
    
    Instead of solving this nonlinear system directly, we solve the modified system
    \begin{alignat}{2}
    \left(\mathbb{A} + \mathbb{D}\right)\tilde{U}^{\nu} &= F - \left(\mathbb{D} + \mathbb{B}(U^{\nu-1})\right)U^{\nu-1}, \qquad \nu \geq 1, \nonumber \\
    U^{\nu} &= \omega \tilde{U}^{\nu} + (1 - \omega)U^{\nu-1}, \label{eq:mod_sys}
    \end{alignat}
    where $\mathbb{D} = \lbrace d_{ij} \rbrace_{i,j=1}^N$ is the artificial diffusion matrix (see Eq.~\eqref{eq:diff_matrix}), $\omega$ is a damping parameter, and $\nu$ denotes the $\nu^{\text{th}}$ iteration step. 
    
    The rationale for solving the modified system is that the matrix $\mathbb{A} + \mathbb{D}$ is a constant M-matrix of non-negative type. Hence, it can be factorized once, stored, and reused within the iterative loop, resulting in computational efficiency. 
    
    If Eq.~\eqref{eq:cdr_eqn} is nonlinear, then the matrix $\mathbb{A}$ also becomes nonlinear, and we directly solve Eq.~\eqref{eq:mod_sys_0}. This approach is referred to as the \emph{fixed-point matrix} strategy in \cite{JJ19}. In this case too, we employ a dynamic damping parameter $\omega$. The stopping criterion is either a maximum of $10^4$ nonlinear iterations or when the $\ell_2$ norm of the residual falls below $10^{-8}\sqrt{\# \mathrm{dofs}}$, where $\#\mathrm{dofs}$ denotes the total number of degrees of freedom.

    \item[$\bullet$] \textbf{ESTIMATE}: To estimate the error in each mesh cell $K$, we use the global upper bound $\eta$ presented in Theorem~\ref{thm:global_upper_bound}. We observe that the estimator includes certain constants. For simplicity, we set $C_{\mathrm{I}}$, $C_{\mathrm{F}}, C$, and $C_{\mathrm{inv}}$ to unity. The constant $C_{\mathrm{edge,max}}$ is computable for 2D triangular elements and is given by
    \begin{equation*}
    C_{\mathrm{edge,max}} = \frac{4\sqrt{2}\left(1+\sqrt{2}\right)|K|}{1 - C_{\cos}\rho_K^3},
    \end{equation*}
    where $\rho_K$ is the diameter of the largest ball inscribed in $K$, $|K|$ is the area of the element $K$, and $C_{\cos} \geq \cos(\theta_i)$ for $i = 1, 2, 3$, with $\theta_i$ denoting the angle opposite the edge $E_i$ of the triangle $K$ (see \cite{Jha21}).

    \item[$\bullet$] \textbf{MARK}: For marking mesh cells, we use the maximum strategy, also known as Dörfler marking \cite{Dor96}. A common issue that arises when marking mesh cells in convection-dominated problems is that only a few cells with high error are marked, which deteriorates the effectiveness of the adaptive algorithm. 

    To overcome this, we adopt a dual marking strategy: in addition to the maximum strategy, we enforce a minimum percentage of mesh cells that must be refined, as suggested in \cite{Joh00}. This strategy involves two parameters: $\mathtt{ref\_tol}$ and $\mathtt{min\_ref}$. 

    Let $\eta_{\max} = \max_{K \in \mathcal{T}_h} \eta_K$. A mesh cell $K$ is marked for refinement if
    \begin{equation*}
    \eta_K \geq \mathtt{ref\_tol} \cdot \eta_{\max}.
    \end{equation*}
    If at least $\mathtt{min\_ref}$ fraction of the cells are marked, the process is complete. Otherwise, we relax $\mathtt{ref\_tol}$ by multiplying it by a factor of $0.8$ and repeat the marking process.

    In our numerical studies, we set $\mathtt{ref\_tol} = 0.5$ and $\mathtt{min\_ref} = 0.05$, corresponding to marking at least $5\%$ of the cells.

    \item[$\bullet$]  \textbf{REFINE}: For mesh refinement, we follow the red-green refinement process \cite{Ver13}. A grid with red-green refinement consists of regularly refined cells arising from red refinement and closure cells originating from green refinement. Both cell types may be marked for refinement by the error indicator.

    In the first step of the refinement process, the parents of closure cells are marked for refinement if any of their children are marked. Note that the parents of closure cells are regularly refined cells from a coarser grid. All closure cells are then removed so that only regularly refined cells remain, and all marked cells are refined using regular refinement. Finally, the resulting grid is closed.

    An alternative approach is to use grids with hanging nodes. Theoretical and numerical results regarding such grids were presented in \cite{JJK23}. It was shown that the AFC scheme with the Kuzmin limiter fails to satisfy the discrete maximum principle (DMP) on grids refined with red-green refinement, whereas all other methods satisfy the DMP on both red-green refined grids and grids with hanging nodes. However, the error estimates and performance of the nonlinear loop remained comparable for both refinement strategies. 

    Therefore, in this work, unless otherwise mentioned, we present results only using grids refined via red-green refinement.
\end{enumerate}

\begin{figure}
    \centering
    \includegraphics[width=0.15\linewidth]{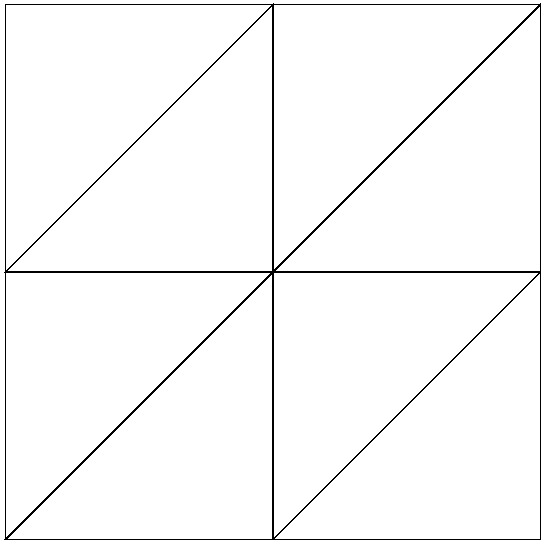}\hspace*{2em}
    \includegraphics[width=0.15\linewidth]{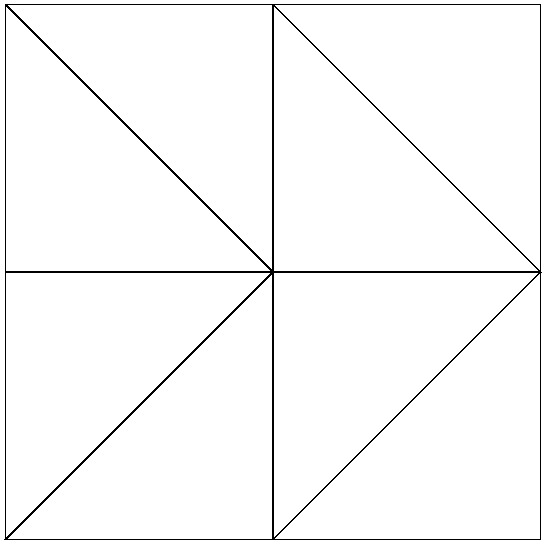}\hspace*{2em}
    \includegraphics[width=0.15\linewidth]{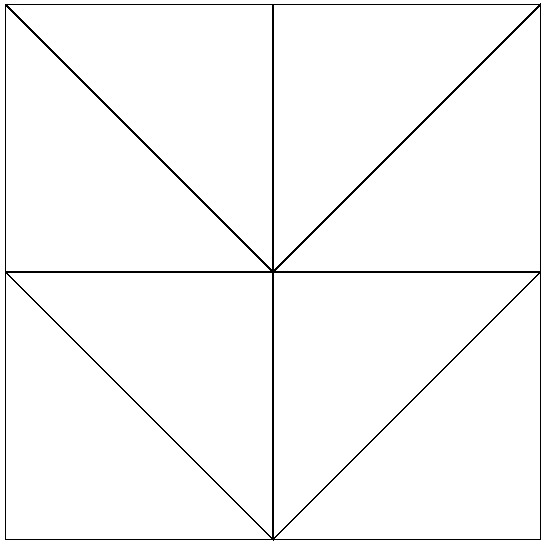}\hspace*{2em}
    \includegraphics[width=0.15\linewidth]{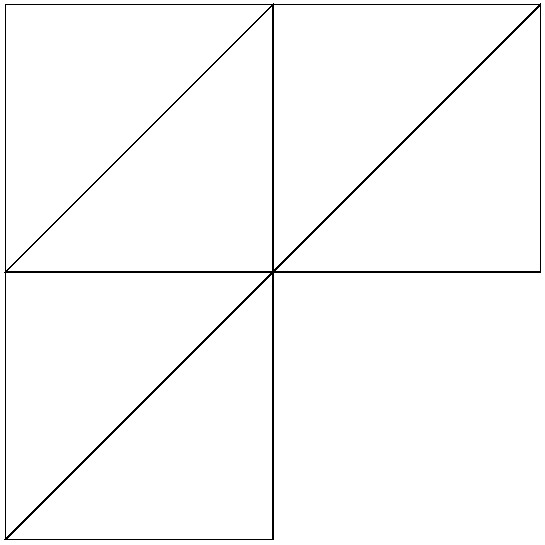}
    \caption{Coarsest grids for the numerical studies, Grid~1,2,3,4 (left to right)}
    \label{fig:grids}
\end{figure}

To ensure a fair comparison of the results, we terminate the adaptive loop once the total number of degrees of freedom reaches $2.5\times 10^5$. For problems defined on $(0,1)^2$ we consider three different types of grids, whose level-0 configurations are shown in Fig.~\ref{fig:grids} (see Grid~1,2,3).  Grid~4 corresponds to the L-shaped domain $[0,1]^2 \setminus [0.5,1] \times [0,0.5]$. 

For our adaptive algorithm, we begin with two uniform refinement steps, then proceed to adaptive refinement. The resulting systems of equations are solved using the sparse direct solver \textsc{UMFPACK}~\cite{Dav04}. All simulations were carried out using the finite element library \textsc{ParMooN}~\cite{WB16}.

The numerical assessment is performed based on the following criteria:
\begin{enumerate}
    \item \textbf{Effectivity index:} The effectivity index provides a measure of the performance of the error estimator. It is defined as
    \begin{equation}
    \gamma_{\mathrm{eff}} = \frac{\eta}{\|u - u_h\|_a}.
    \end{equation}
    Ideally, we seek $\gamma_{\mathrm{eff}} \approx 1$, independently of the choice of $\varepsilon$ (robustness).

    \item \textbf{Accuracy of the solution:} Plots of the error in the solution and the gradient of the error in the $\mathrm{L}^2$ norm are presented.

    \item \textbf{Sharpness of the layer:} The thickness of the interior layer is plotted against the total number of degrees of freedom ($\#\mathrm{dofs}$). For this, we compute 
\begin{equation}\label{eq:smear}
\mathrm{smear}_{\mathrm{int}}=y_2-y_1,
\end{equation}
along $x=a$, where $y_1$ is the $y-$ coordinate of the first point on the cut line $(a,y_1)$ with $u_h(a,y_1)\geq 0.1$ and $y_2$ is the $y-$ coordinate of the first point with $u_h(a,y_2)\geq 0.9$. 

    \item \textbf{Quality of the adaptively refined grids:} Visual representations of the solution and the corresponding adaptively refined grids are presented.

    \item \textbf{Efficiency of the nonlinear solver:} Since the algebraic stabilization schemes are nonlinear in nature, we plot the number of iterations and rejections with respect to the total $\#\mathrm{dofs}$. As the computational cost of a rejection in the nonlinear loop is equivalent to that of an iteration, both are plotted for comparison.

    \item \textbf{Computational cost:} The total time taken to solve the problem is reported.
\end{enumerate}

%\subsection{Known Boundary Layer}\label{ex:example_1_ABR17}
%\input{example_1_ABR17.tex}

\subsection{Solution with Regular Boundary Layer}\label{ex:boundary_layer}
\begin{figure}[tbp]
    \centering
    \includegraphics[width=0.3\linewidth]{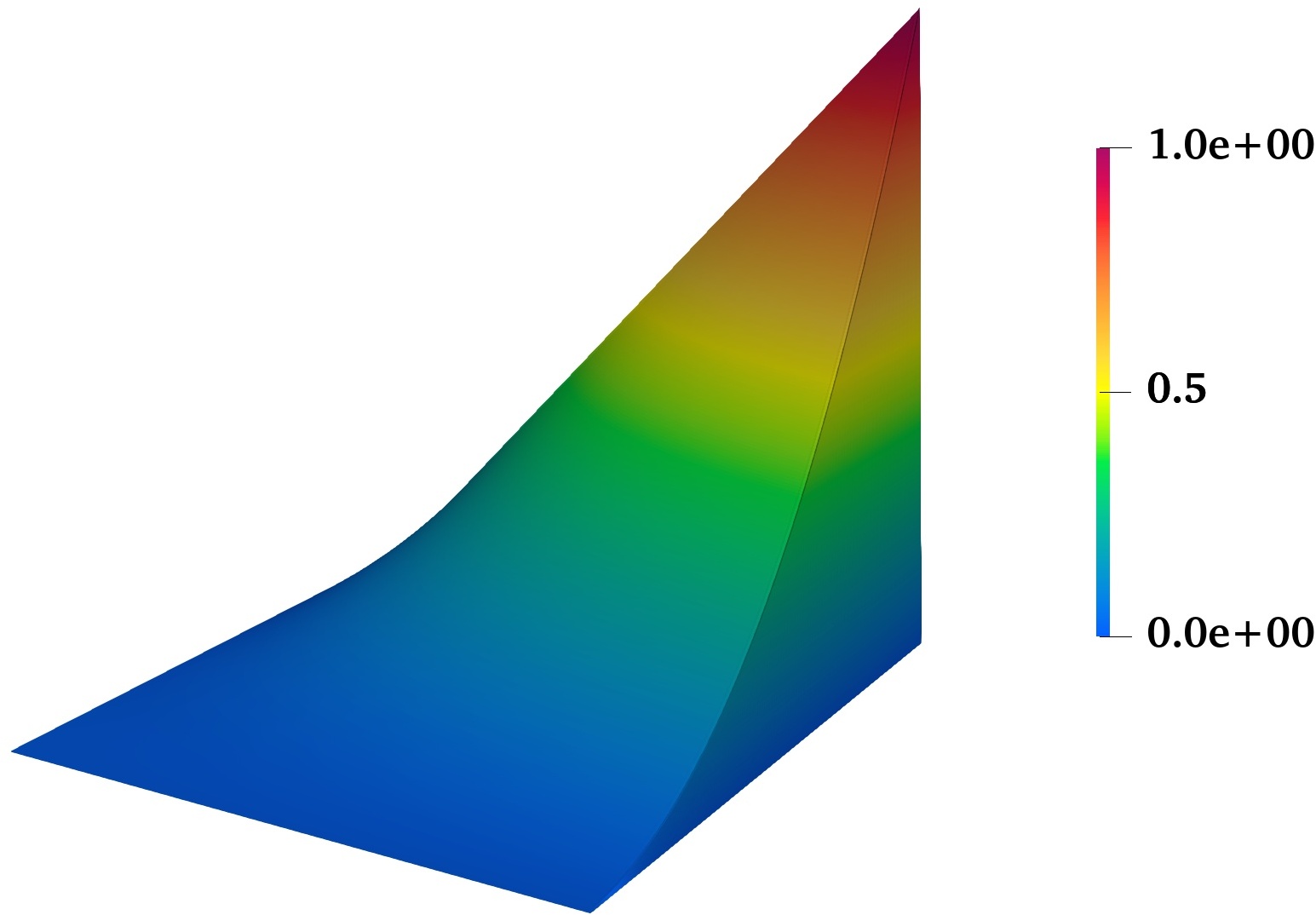}
    \caption{Example~\ref{ex:boundary_layer}: Numerical solution with the MUAS method on Grid~1 with $16641$ degrees of freedom.}
    \label{fig:boundary_layer}
\end{figure}

Consider Eq.~\eqref{eq:cdr_eqn} with $\varepsilon = 10^{-2}$, $\bb = (2,3)^{\top}$, $c = 1$, 
$\Omega = (0,1)^2$, and $\Gamma_{\mD} = \Gamma$. 
The right-hand side and the Dirichlet boundary condition are chosen such that
$$
u(x,y) = xy^2 
- y^2 \exp\left( \frac{2(x - 1)}{\varepsilon} \right) 
- x \exp\left( \frac{3(y - 1)}{\varepsilon} \right)
+ \exp\left( \frac{2(x - 1) + 3(y - 1)}{\varepsilon} \right)
$$
is the exact solution; see Fig.~\ref{fig:boundary_layer}. 
The solution exhibits regular boundary layers along the outflow boundaries $x = 1$ and $y = 1$.

\begin{figure}[tbp]
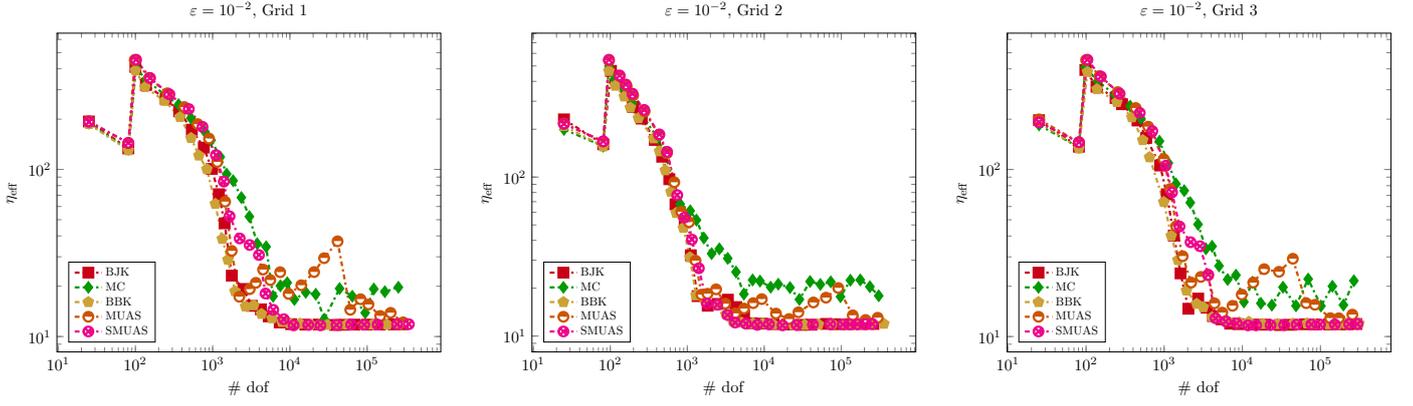

\centerline{
}
\caption{Example~\ref{ex:boundary_layer}: Effectivity index for Grids~1–3 (left to right). }
\label{fig:effec_index_boundary_layer_grid_1}
\end{figure}

We use this example to evaluate both the effectivity index and the accuracy of the numerical solution.  
Fig.~\ref{fig:effec_index_boundary_layer_grid_1} plots the effectivity index for Grids~1, 2, and~3. 
The effectivity index for the MC limiter is the largest, approaching~20, whereas for all other methods, the values are comparable and close to~12 on sufficiently refined grids.  
Furthermore, we observe robustness with respect to the choice of grid.

\begin{figure}[tbp]
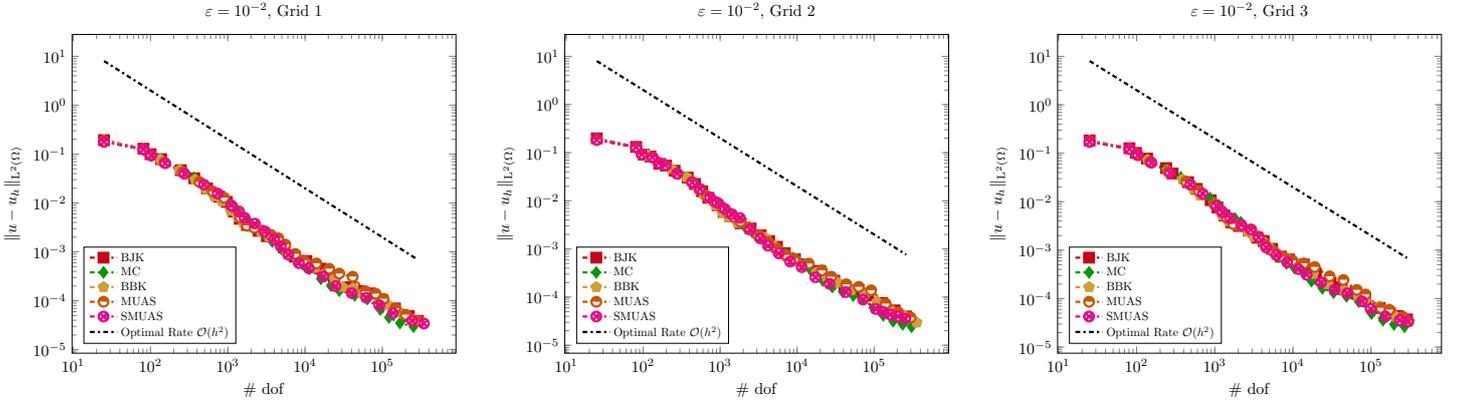

\centerline{
%
}
\caption{Example~\ref{ex:boundary_layer}: Error in the $\mL^2(\Omega)$ norm for Grids~1–3 (left to right).}
\label{fig:l2_error_boundary_layer_grid_1}
\end{figure}

\begin{figure}[tbp]
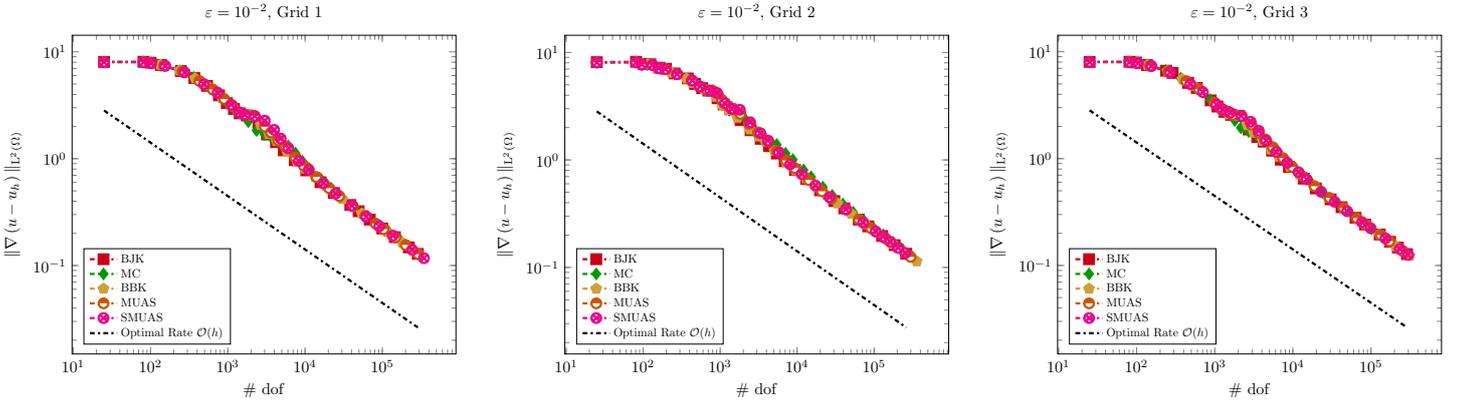

\centerline{
%
}
\caption{Example~\ref{ex:boundary_layer}: Error of the gradient in the $\mL^2(\Omega)$ norm for Grids~1–3 (left to right).}
\label{fig:h1_error_boundary_layer_grid_1}
\end{figure}

Figs.~\ref{fig:l2_error_boundary_layer_grid_1} and \ref{fig:h1_error_boundary_layer_grid_1} 
show the $\mL^2$ error of the solution and the $\mL^2$ error of its gradient, respectively.  
We observe that all limiters exhibit optimal convergence rates for all three grids.

\begin{figure}[tbp]
\centerline{
\begin{tikzpicture}[scale=0.6]
\begin{loglogaxis}[
    legend pos=south west, xlabel = $\#\ \mathrm{dof}$, ylabel=$\eta_3$,
    legend cell align ={left}, title = {$\varepsilon=10^{-2}$, Grid~1},
    legend style={nodes={scale=0.75, transform shape}}]
\addplot[color=crimson,  mark=square*, line width = 0.5mm, dashdotted,,mark options = {scale= 1.5, solid}]
coordinates{( 25.0 , 154.952 )( 81.0 , 105.737 )( 100.0 , 326.815 )( 137.0 , 241.654 )( 247.0 , 174.937 )( 369.0 , 126.644 )( 538.0 , 82.9213 )( 769.0 , 52.6861 )( 973.0 , 33.1074 )( 1203.0 , 20.4815 )( 1428.0 , 12.2666 )( 1771.0 , 5.10039 )( 2499.0 , 3.16538 )( 3215.0 , 1.70626 )( 4213.0 , 1.30014 )( 5320.0 , 0.799271 )( 7382.0 , 0.385636 )( 10381.0 , 0.192091 )( 15868.0 , 0.12911 )( 23601.0 , 0.100698 )( 38598.0 , 0.0385208 )( 49440.0 , 0.0259502 )( 68782.0 , 0.00682706 )( 100140.0 , 0.0145553 )( 145678.0 , 0.00752662 )( 220575.0 , 0.00453123 )( 288016.0 , 0.00264885 )};
\addlegendentry{BJK} 
\addplot[color=dark_green,  mark=diamond*, line width = 0.5mm, dashdotted,,mark options = {scale= 1.5, solid}]
coordinates{( 25.0 , 150.771 )( 81.0 , 104.427 )( 102.0 , 322.059 )( 137.0 , 258.246 )( 251.0 , 181.329 )( 362.0 , 140.994 )( 522.0 , 99.1102 )( 798.0 , 63.7392 )( 953.0 , 50.0285 )( 1254.0 , 32.7563 )( 1517.0 , 23.9315 )( 1840.0 , 18.9062 )( 2352.0 , 12.326 )( 2986.0 , 8.64017 )( 3773.0 , 5.25863 )( 4924.0 , 4.50308 )( 6057.0 , 1.74163 )( 7532.0 , 1.86658 )( 9235.0 , 1.66202 )( 11419.0 , 0.967508 )( 16029.0 , 0.925183 )( 21610.0 , 0.726879 )( 27814.0 , 0.238008 )( 43994.0 , 0.551083 )( 64249.0 , 0.370969 )( 94383.0 , 0.172022 )( 121913.0 , 0.314412 )( 169425.0 , 0.248497 )( 253695.0 , 0.220242 )};
\addlegendentry{MC} 
%\addplot[color=amethyst,  mark=triangle*, line width = 0.5mm, dashdotted,,mark options = {scale= 1.5, solid}]
%coordinates{( 25.0 , 101.369 )( 81.0 , 64.7105 )( 100.0 , 198.87 )( 163.0 , 138.485 )( 238.0 , 102.123 )( 341.0 , 79.1911 )( 469.0 , 57.1077 )( 588.0 , 42.9445 )( 724.0 , 31.74 )( 968.0 , 18.2944 )( 1158.0 , 9.92929 )( 1381.0 , 4.82123 )( 1649.0 , 2.84253 )( 2132.0 , 2.19931 )( 2640.0 , 1.43409 )( 3287.0 , 1.008 )( 4633.0 , 0.626009 )( 6267.0 , 0.411921 )( 8286.0 , 0.279284 )( 14557.0 , 0.15625 )( 20572.0 , 0.0885824 )( 27452.0 , 0.0504996 )( 38593.0 , 0.0272592 )( 67998.0 , 0.0155509 )( 89130.0 , 0.0085391 )( 129235.0 , 0.00549378 )( 179717.0 , 0.00353685 )( 235558.0 , 0.0024526 )( 317502.0 , 0.00115638 )};
%\addlegendentry{BBK} 
\addplot[color=gold,  mark=pentagon*, line width = 0.5mm, dashdotted,,mark options = {scale= 1.5, solid}]
coordinates{( 25.0 , 149.175 )( 81.0 , 102.764 )( 100.0 , 310.737 )( 133.0 , 235.772 )( 239.0 , 172.252 )( 378.0 , 116.124 )( 525.0 , 74.9561 )( 664.0 , 53.3447 )( 843.0 , 36.856 )( 1083.0 , 19.1888 )( 1333.0 , 10.2297 )( 1576.0 , 6.83348 )( 1947.0 , 3.7448 )( 2665.0 , 2.05022 )( 3404.0 , 1.75575 )( 4288.0 , 1.11304 )( 5855.0 , 0.692037 )( 7877.0 , 0.389307 )( 9967.0 , 0.253457 )( 13517.0 , 0.191042 )( 18359.0 , 0.128971 )( 23885.0 , 0.0903112 )( 30415.0 , 0.062916 )( 41853.0 , 0.045341 )( 54295.0 , 0.0271977 )( 89019.0 , 0.0143064 )( 117020.0 , 0.00935776 )( 180938.0 , 0.00488997 )( 252501.0 , 0.00467344 )};
\addlegendentry{BBK} 
\addplot[color=carrotorange,  mark=halfcircle*, line width = 0.5mm, dashdotted,,mark options = {scale= 1.5, solid}]
coordinates{( 25.0 , 155.704 )( 81.0 , 114.151 )( 98.0 , 358.375 )( 151.0 , 263.273 )( 260.0 , 186.688 )( 430.0 , 125.978 )( 640.0 , 82.7924 )( 897.0 , 55.3008 )( 1147.0 , 34.4075 )( 1449.0 , 16.9648 )( 1764.0 , 7.7824 )( 2216.0 , 3.31414 )( 3021.0 , 3.10614 )( 3621.0 , 3.10421 )( 4494.0 , 3.32957 )( 5571.0 , 2.41742 )( 7480.0 , 2.20022 )( 9622.0 , 1.1981 )( 14044.0 , 1.12548 )( 20562.0 , 1.13081 )( 27729.0 , 1.20617 )( 41422.0 , 1.2726 )( 60507.0 , 0.246777 )( 82182.0 , 0.293359 )( 103932.0 , 0.222503 )( 147413.0 , 0.112466 )( 201252.0 , 0.11297 )( 292779.0 , 0.0352837 )};
\addlegendentry{MUAS} 
\addplot[color=magenta,  mark=otimes, line width = 0.5mm, dashdotted,,mark options = {scale= 1.5, solid}]
coordinates{( 25.0 , 154.06 )( 81.0 , 114.699 )( 102.0 , 357.168 )( 155.0 , 260.989 )( 276.0 , 180.748 )( 491.0 , 114.397 )( 739.0 , 74.4163 )( 1103.0 , 38.8812 )( 1395.0 , 22.8738 )( 1653.0 , 13.2041 )( 2237.0 , 9.24864 )( 3008.0 , 7.62221 )( 3977.0 , 5.34128 )( 4836.0 , 2.19036 )( 6087.0 , 1.13252 )( 8331.0 , 0.541417 )( 11175.0 , 0.187344 )( 17272.0 , 0.133777 )( 24818.0 , 0.0622233 )( 40268.0 , 0.0354421 )( 62122.0 , 0.0179832 )( 91539.0 , 0.00987084 )( 136982.0 , 0.00434573 )( 243674.0 , 0.00544994 )( 345984.0 , 0.00164995 )};
\addlegendentry{SMUAS} 
\addplot[color=black,  line width = 0.5mm, dashdotted,,mark options = {scale= 1.0, solid}]
coordinates{( 25.0 , 7.999999999999998 )( 81.0 , 2.469135802469136 )( 100.0 , 1.9999999999999996 )( 137.0 , 1.4598540145985401 )( 247.0 , 0.8097165991902833 )( 369.0 , 0.5420054200542005 )( 538.0 , 0.37174721189591087 )( 769.0 , 0.260078023407022 )( 973.0 , 0.2055498458376156 )( 1203.0 , 0.16625103906899416 )( 1428.0 , 0.14005602240896362 )( 1771.0 , 0.11293054771315642 )( 2499.0 , 0.08003201280512202 )( 3215.0 , 0.062208398133748066 )( 4213.0 , 0.047472110135295516 )( 5320.0 , 0.037593984962406006 )( 7382.0 , 0.027092928745597395 )( 10381.0 , 0.01926596666987766 )( 15868.0 , 0.01260398285858331 )( 23601.0 , 0.008474217194186688 )( 38598.0 , 0.005181615627752734 )( 49440.0 , 0.004045307443365695 )( 68782.0 , 0.0029077374894594514 )( 100140.0 , 0.0019972039145196726 )( 145678.0 , 0.0013728908963604664 )( 220575.0 , 0.0009067210699308625 )( 288016.0 , 0.0006944058663407588 )};
\addlegendentry{Optimal Rate $\mathcal{O}(h^2)$} 
\end{loglogaxis}
\end{tikzpicture}\hspace*{1em}
\begin{tikzpicture}[scale=0.6]
\begin{loglogaxis}[
    legend pos=south west, xlabel = $\#\ \mathrm{dof}$, ylabel=$\eta_3$,
    legend cell align ={left}, title = {$\varepsilon=10^{-2}$, Grid~2},
    legend style={nodes={scale=0.75, transform shape}}]
\addplot[color=crimson,  mark=square*, line width = 0.5mm, dashdotted,,mark options = {scale= 1.5, solid}]
coordinates{( 25.0 , 188.367 )( 81.0 , 130.103 )( 102.0 , 364.151 )( 124.0 , 307.851 )( 162.0 , 238.396 )( 198.0 , 194.224 )( 256.0 , 148.633 )( 379.0 , 97.8021 )( 472.0 , 68.0126 )( 584.0 , 44.6733 )( 703.0 , 29.2611 )( 927.0 , 19.8504 )( 1112.0 , 9.85658 )( 1363.0 , 4.03195 )( 1846.0 , 2.43266 )( 2490.0 , 2.03873 )( 3372.0 , 1.96885 )( 4331.0 , 1.35238 )( 5468.0 , 0.803702 )( 6903.0 , 0.320186 )( 9383.0 , 0.218071 )( 13195.0 , 0.149951 )( 19911.0 , 0.100438 )( 30406.0 , 0.0398742 )( 41019.0 , 0.0224291 )( 63328.0 , 0.0307829 )( 83949.0 , 0.0196883 )( 127336.0 , 0.0180556 )( 184946.0 , 0.0115908 )( 261546.0 , 0.00249089 )};
\addlegendentry{BJK} 
\addplot[color=dark_green,  mark=diamond*, line width = 0.5mm, dashdotted,,mark options = {scale= 1.5, solid}]
coordinates{( 25.0 , 161.236 )( 81.0 , 126.102 )( 96.0 , 387.241 )( 113.0 , 322.051 )( 151.0 , 260.624 )( 194.0 , 203.716 )( 259.0 , 155.046 )( 414.0 , 95.8104 )( 499.0 , 70.4541 )( 680.0 , 39.8764 )( 808.0 , 29.0691 )( 1092.0 , 21.0216 )( 1316.0 , 15.2248 )( 1643.0 , 10.2754 )( 2073.0 , 7.5322 )( 2611.0 , 6.37116 )( 3359.0 , 4.59656 )( 4301.0 , 3.33209 )( 5444.0 , 2.00509 )( 7070.0 , 2.23923 )( 8815.0 , 1.90362 )( 12129.0 , 1.34249 )( 15214.0 , 1.22953 )( 21603.0 , 0.919869 )( 28638.0 , 0.576436 )( 39371.0 , 0.72403 )( 53562.0 , 0.570492 )( 73161.0 , 0.505124 )( 97832.0 , 0.307242 )( 131542.0 , 0.381901 )( 176748.0 , 0.327212 )( 235485.0 , 0.243219 )( 306961.0 , 0.169599 )};
\addlegendentry{MC} 
%\addplot[color=amethyst,  mark=triangle*, line width = 0.5mm, dashdotted,,mark options = {scale= 1.5, solid}]
%coordinates{( 25.0 , 92.5466 )( 81.0 , 65.45 )( 101.0 , 196.216 )( 133.0 , 141.968 )( 179.0 , 111.452 )( 219.0 , 91.0608 )( 266.0 , 79.6403 )( 410.0 , 57.6383 )( 563.0 , 37.7965 )( 688.0 , 28.9394 )( 879.0 , 18.6275 )( 1094.0 , 10.0561 )( 1318.0 , 4.31877 )( 1609.0 , 2.24385 )( 2229.0 , 1.6743 )( 2937.0 , 0.987011 )( 4426.0 , 0.528219 )( 5447.0 , 0.355668 )( 7637.0 , 0.222929 )( 10467.0 , 0.148172 )( 17277.0 , 0.0968016 )( 26547.0 , 0.0472719 )( 34842.0 , 0.0215271 )( 59079.0 , 0.00902408 )( 105469.0 , 0.0083088 )( 160874.0 , 0.00177324 )( 280932.0 , 0.00207846 )};
%\addlegendentry{BBK} 
\addplot[color=gold,  mark=pentagon*, line width = 0.5mm, dashdotted,,mark options = {scale= 1.5, solid}]
coordinates{( 25.0 , 170.089 )( 81.0 , 125.752 )( 96.0 , 363.058 )( 116.0 , 292.412 )( 152.0 , 230.518 )( 182.0 , 193.695 )( 232.0 , 151.988 )( 365.0 , 99.8492 )( 436.0 , 75.4823 )( 518.0 , 53.5505 )( 614.0 , 36.3734 )( 718.0 , 25.4978 )( 889.0 , 17.9418 )( 1082.0 , 9.19216 )( 1299.0 , 4.05194 )( 1677.0 , 3.04811 )( 2441.0 , 1.97819 )( 3173.0 , 1.25457 )( 3897.0 , 0.745904 )( 5077.0 , 0.490905 )( 7155.0 , 0.306503 )( 9525.0 , 0.211979 )( 13038.0 , 0.140289 )( 17939.0 , 0.101867 )( 24015.0 , 0.0451851 )( 33559.0 , 0.0282635 )( 50552.0 , 0.0182586 )( 76341.0 , 0.027631 )( 109487.0 , 0.010689 )( 138471.0 , 0.00702822 )( 176137.0 , 0.00636528 )( 248512.0 , 0.00878026 )( 358141.0 , 0.00593367 )};
\addlegendentry{BBK} 
\addplot[color=carrotorange,  mark=halfcircle*, line width = 0.5mm, dashdotted,,mark options = {scale= 1.5, solid}]
coordinates{( 25.0 , 176.312 )( 81.0 , 135.147 )( 96.0 , 419.577 )( 132.0 , 333.481 )( 163.0 , 267.051 )( 198.0 , 228.414 )( 279.0 , 161.991 )( 428.0 , 100.458 )( 545.0 , 67.8903 )( 678.0 , 40.8783 )( 837.0 , 25.6847 )( 1054.0 , 18.4081 )( 1282.0 , 8.73117 )( 1529.0 , 4.20235 )( 1832.0 , 3.91652 )( 2370.0 , 3.48592 )( 3113.0 , 2.03868 )( 4199.0 , 1.18959 )( 5465.0 , 1.15667 )( 7183.0 , 1.11526 )( 9295.0 , 1.05605 )( 14210.0 , 0.494019 )( 20985.0 , 0.270674 )( 27792.0 , 0.344237 )( 43496.0 , 0.387181 )( 63382.0 , 0.357857 )( 91557.0 , 0.373877 )( 141259.0 , 0.12003 )( 207268.0 , 0.0643741 )( 297813.0 , 0.0661403 )};
\addlegendentry{MUAS} 
\addplot[color=magenta,  mark=otimes, line width = 0.5mm, dashdotted,,mark options = {scale= 1.5, solid}]
coordinates{( 25.0 , 175.382 )( 81.0 , 135.809 )( 96.0 , 424.277 )( 132.0 , 334.677 )( 159.0 , 273.445 )( 194.0 , 234.841 )( 275.0 , 166.348 )( 434.0 , 101.664 )( 547.0 , 69.1312 )( 741.0 , 33.4353 )( 915.0 , 22.7675 )( 1149.0 , 13.0588 )( 1414.0 , 7.22718 )( 1805.0 , 3.30411 )( 2456.0 , 2.46642 )( 3336.0 , 1.33102 )( 4194.0 , 0.671469 )( 5794.0 , 0.407623 )( 8124.0 , 0.202006 )( 11502.0 , 0.171398 )( 17293.0 , 0.0811082 )( 26418.0 , 0.0484185 )( 42937.0 , 0.023639 )( 72033.0 , 0.0147916 )( 104800.0 , 0.00965916 )( 130524.0 , 0.00803963 )( 162200.0 , 0.00895126 )( 205277.0 , 0.00320215 )( 254471.0 , 0.00158608 )};
\addlegendentry{SMUAS} 
\addplot[color=black,  line width = 0.5mm, dashdotted,,mark options = {scale= 1.0, solid}]
coordinates{( 25.0 , 7.999999999999998 )( 81.0 , 2.469135802469136 )( 102.0 , 1.9607843137254903 )( 124.0 , 1.6129032258064515 )( 162.0 , 1.234567901234568 )( 198.0 , 1.0101010101010102 )( 256.0 , 0.7812500000000001 )( 379.0 , 0.5277044854881267 )( 472.0 , 0.42372881355932196 )( 584.0 , 0.34246575342465757 )( 703.0 , 0.2844950213371267 )( 927.0 , 0.2157497303128371 )( 1112.0 , 0.17985611510791363 )( 1363.0 , 0.14673514306676452 )( 1846.0 , 0.10834236186348863 )( 2490.0 , 0.08032128514056226 )( 3372.0 , 0.05931198102016608 )( 4331.0 , 0.046178711613945975 )( 5468.0 , 0.0365764447695684 )( 6903.0 , 0.028972910328842527 )( 9383.0 , 0.021315144410103383 )( 13195.0 , 0.015157256536566882 )( 19911.0 , 0.010044698910150167 )( 30406.0 , 0.0065776491481944355 )( 41019.0 , 0.00487578926838782 )( 63328.0 , 0.0031581606872157658 )( 83949.0 , 0.0023823988373893675 )( 127336.0 , 0.0015706477351259663 )( 184946.0 , 0.0010813967320190759 )( 261546.0 , 0.000764683841465746 )};
\addlegendentry{Optimal Rate $\mathcal{O}(h^2)$} 
\end{loglogaxis}
\end{tikzpicture}\hspace*{1em}
\begin{tikzpicture}[scale=0.6]
\begin{loglogaxis}[
    legend pos=south west, xlabel = $\#\ \mathrm{dof}$, ylabel=$\eta_3$,
    legend cell align ={left}, title = {$\varepsilon=10^{-2}$, Grid~3},
    legend style={nodes={scale=0.75, transform shape}}]
\addplot[color=crimson,  mark=square*, line width = 0.5mm, dashdotted,,mark options = {scale= 1.5, solid}]
coordinates{( 25.0 , 157.243 )( 81.0 , 108.216 )( 98.0 , 314.818 )( 136.0 , 235.557 )( 239.0 , 177.816 )( 289.0 , 155.667 )( 456.0 , 100.261 )( 589.0 , 70.5876 )( 887.0 , 36.7546 )( 1081.0 , 21.4996 )( 1330.0 , 10.5248 )( 1622.0 , 5.35523 )( 2041.0 , 2.3106 )( 2726.0 , 2.29432 )( 3497.0 , 1.51725 )( 4220.0 , 0.97567 )( 5411.0 , 0.644877 )( 6994.0 , 0.304518 )( 8915.0 , 0.202343 )( 13787.0 , 0.133805 )( 19670.0 , 0.101363 )( 30395.0 , 0.0340011 )( 41994.0 , 0.0173064 )( 63451.0 , 0.0172746 )( 84087.0 , 0.0146628 )( 129247.0 , 0.00929763 )( 175063.0 , 0.00682453 )( 220472.0 , 0.0027722 )( 288698.0 , 0.00153192 )};
\addlegendentry{BJK} 
\addplot[color=dark_green,  mark=diamond*, line width = 0.5mm, dashdotted,,mark options = {scale= 1.5, solid}]
coordinates{( 25.0 , 146.41 )( 81.0 , 107.146 )( 98.0 , 324.391 )( 141.0 , 253.495 )( 235.0 , 185.599 )( 368.0 , 134.907 )( 503.0 , 98.2261 )( 648.0 , 78.2862 )( 871.0 , 52.4101 )( 1130.0 , 32.0215 )( 1416.0 , 22.2995 )( 1806.0 , 16.8285 )( 2158.0 , 12.1058 )( 2711.0 , 8.18727 )( 3422.0 , 5.23727 )( 4101.0 , 4.95741 )( 5030.0 , 3.28396 )( 6487.0 , 2.17612 )( 8288.0 , 1.93607 )( 10476.0 , 0.909145 )( 13509.0 , 1.22117 )( 17259.0 , 0.635564 )( 23858.0 , 0.489162 )( 33165.0 , 0.642997 )( 48640.0 , 0.31954 )( 69699.0 , 0.502532 )( 101558.0 , 0.213296 )( 142358.0 , 0.301519 )( 196969.0 , 0.155689 )( 268083.0 , 0.24034 )};
\addlegendentry{MC} 
%\addplot[color=amethyst,  mark=triangle*, line width = 0.5mm, dashdotted,,mark options = {scale= 1.5, solid}]
%coordinates{( 25.0 , 95.3586 )( 81.0 , 65.5142 )( 111.0 , 175.183 )( 156.0 , 131.775 )( 236.0 , 96.6513 )( 330.0 , 77.0275 )( 446.0 , 58.9491 )( 556.0 , 44.2778 )( 707.0 , 31.9599 )( 914.0 , 20.17 )( 1109.0 , 11.0442 )( 1354.0 , 4.81757 )( 1613.0 , 2.80241 )( 2112.0 , 2.19189 )( 2629.0 , 1.38476 )( 3253.0 , 0.952525 )( 4569.0 , 0.562672 )( 6137.0 , 0.386106 )( 8220.0 , 0.269647 )( 12230.0 , 0.172541 )( 18789.0 , 0.0868037 )( 27996.0 , 0.0439766 )( 42340.0 , 0.023557 )( 74939.0 , 0.017086 )( 93266.0 , 0.0137066 )( 135823.0 , 0.00531842 )( 187737.0 , 0.00268632 )( 296798.0 , 0.00224089 )};
%\addlegendentry{BBK} 
\addplot[color=gold,  mark=pentagon*, line width = 0.5mm, dashdotted,,mark options = {scale= 1.5, solid}]
coordinates{( 25.0 , 152.916 )( 81.0 , 105.4 )( 104.0 , 303.331 )( 138.0 , 227.231 )( 248.0 , 166.761 )( 386.0 , 113.815 )( 532.0 , 71.683 )( 650.0 , 52.2756 )( 988.0 , 20.4837 )( 1232.0 , 10.8375 )( 1500.0 , 6.838 )( 1929.0 , 3.79552 )( 2606.0 , 2.2283 )( 3162.0 , 1.70393 )( 4098.0 , 0.943454 )( 5256.0 , 0.718293 )( 6600.0 , 0.462525 )( 9168.0 , 0.296827 )( 12114.0 , 0.272857 )( 18269.0 , 0.142877 )( 28092.0 , 0.0925852 )( 44580.0 , 0.0724618 )( 67549.0 , 0.05178 )( 84272.0 , 0.0296441 )( 128315.0 , 0.0192565 )( 163371.0 , 0.0164136 )( 211044.0 , 0.0203995 )( 310772.0 , 0.0175567 )};
\addlegendentry{BBK} 
\addplot[color=carrotorange,  mark=halfcircle*, line width = 0.5mm, dashdotted,,mark options = {scale= 1.5, solid}]
coordinates{( 25.0 , 159.204 )( 81.0 , 115.47 )( 100.0 , 356.687 )( 149.0 , 268.207 )( 255.0 , 189.294 )( 432.0 , 123.764 )( 625.0 , 80.4928 )( 987.0 , 38.5935 )( 1195.0 , 22.3787 )( 1425.0 , 12.2629 )( 1724.0 , 7.19516 )( 2217.0 , 4.44565 )( 2947.0 , 4.09938 )( 3648.0 , 3.53964 )( 4405.0 , 1.67582 )( 5542.0 , 1.03047 )( 7488.0 , 1.02352 )( 9829.0 , 1.13758 )( 13804.0 , 1.18495 )( 19032.0 , 1.23849 )( 29497.0 , 0.918828 )( 44353.0 , 0.92396 )( 65689.0 , 0.293217 )( 84952.0 , 0.243029 )( 133642.0 , 0.0986772 )( 171817.0 , 0.0853526 )( 256278.0 , 0.088158 )};
\addlegendentry{MUAS} 
\addplot[color=magenta,  mark=otimes, line width = 0.5mm, dashdotted,,mark options = {scale= 1.5, solid}]
coordinates{( 25.0 , 152.639 )( 81.0 , 115.591 )( 104.0 , 355.603 )( 155.0 , 264.727 )( 269.0 , 182.285 )( 487.0 , 108.523 )( 705.0 , 70.828 )( 1050.0 , 33.523 )( 1260.0 , 20.7546 )( 1572.0 , 11.8799 )( 2177.0 , 8.76885 )( 2871.0 , 7.27958 )( 3654.0 , 3.7537 )( 4545.0 , 0.910139 )( 6135.0 , 0.546433 )( 8467.0 , 0.320398 )( 11802.0 , 0.135917 )( 15016.0 , 0.106352 )( 22992.0 , 0.071732 )( 34279.0 , 0.0361788 )( 49613.0 , 0.0194005 )( 76632.0 , 0.0126469 )( 99006.0 , 0.00984315 )( 152824.0 , 0.00400975 )( 233337.0 , 0.0037837 )( 297601.0 , 0.00174591 )};
\addlegendentry{SMUAS} 
\addplot[color=black,  line width = 0.5mm, dashdotted,,mark options = {scale= 1.0, solid}]
coordinates{( 25.0 , 7.999999999999998 )( 81.0 , 2.469135802469136 )( 98.0 , 2.0408163265306123 )( 136.0 , 1.4705882352941173 )( 239.0 , 0.8368200836820082 )( 289.0 , 0.6920415224913494 )( 456.0 , 0.43859649122807015 )( 589.0 , 0.3395585738539899 )( 887.0 , 0.2254791431792559 )( 1081.0 , 0.18501387604070305 )( 1330.0 , 0.15037593984962402 )( 1622.0 , 0.12330456226880396 )( 2041.0 , 0.09799118079372857 )( 2726.0 , 0.07336757153338225 )( 3497.0 , 0.057191878753217046 )( 4220.0 , 0.04739336492890995 )( 5411.0 , 0.036961744594344856 )( 6994.0 , 0.02859593937660853 )( 8915.0 , 0.022434099831744245 )( 13787.0 , 0.014506419090447522 )( 19670.0 , 0.010167768174885612 )( 30395.0 , 0.006580029610133245 )( 41994.0 , 0.00476258513120922 )( 63451.0 , 0.00315203858095223 )( 84087.0 , 0.0023784889459726242 )( 129247.0 , 0.001547424698445612 )( 175063.0 , 0.0011424458623466981 )( 220472.0 , 0.0009071446714321999 )( 288698.0 , 0.0006927654504014576 )};
\addlegendentry{Optimal Rate $\mathcal{O}(h^2)$} 
\end{loglogaxis}
\end{tikzpicture}}
\caption{Example~\ref{ex:boundary_layer}: Values of $\eta_3$ for Grids~1–3 (left to right).}
\label{fig:eta3_error_boundary_layer_grid_1}
\end{figure}
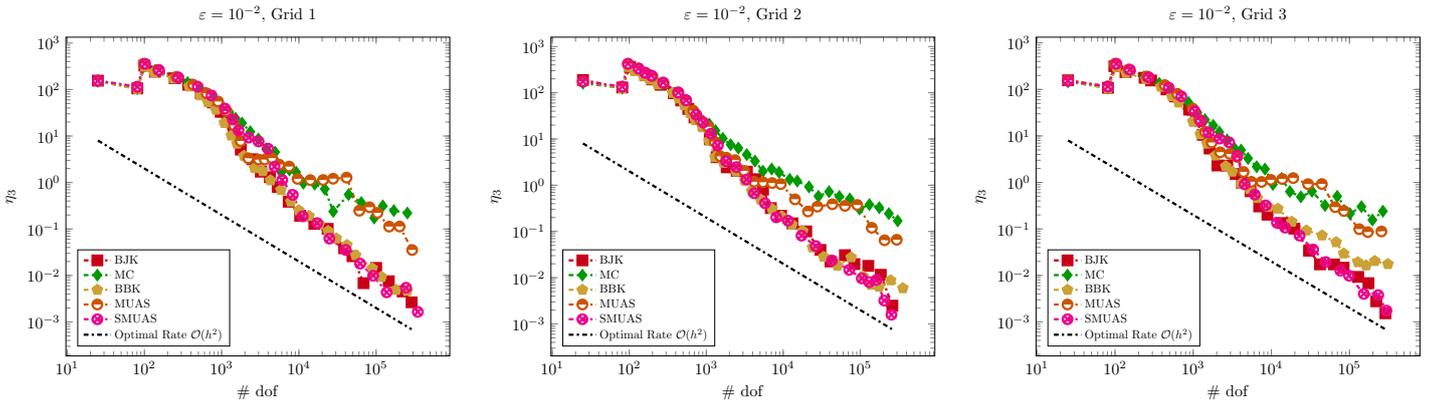

Next, we investigate the decay of $\eta_3$, the term arising from the stabilization part.  
In Fig.~\ref{fig:eta3_error_boundary_layer_grid_1}, we observe that $\eta_3$ decays at a rate of 
$\mathcal{O}(h^2)$ for all methods except the monolithic and MUAS limiters.  
This deterioration is independent of the choice of grid.  
Although the convergence rate of $\eta_3$ deteriorates on sufficiently refined meshes, this behaviour does not adversely affect the accuracy of the solution itself.  
One of the motivations for the SMUAS method was to achieve optimal convergence on arbitrary grids \cite{Kn23}, and here we clearly see that SMUAS outperforms the MUAS method.

From this example, we conclude that the a posteriori error estimator delivers comparable results across all limiters and is robust with respect to the choice of grid. Based on this example, the BJK, BBK, and the MUAS method seem to be the most promising.

\subsection{Solution with Corner Singularity on L-shaped Domain}\label{ex:corner_singularity}
\begin{figure}[tbp]
    \centering
    \includegraphics[width=0.4\linewidth]{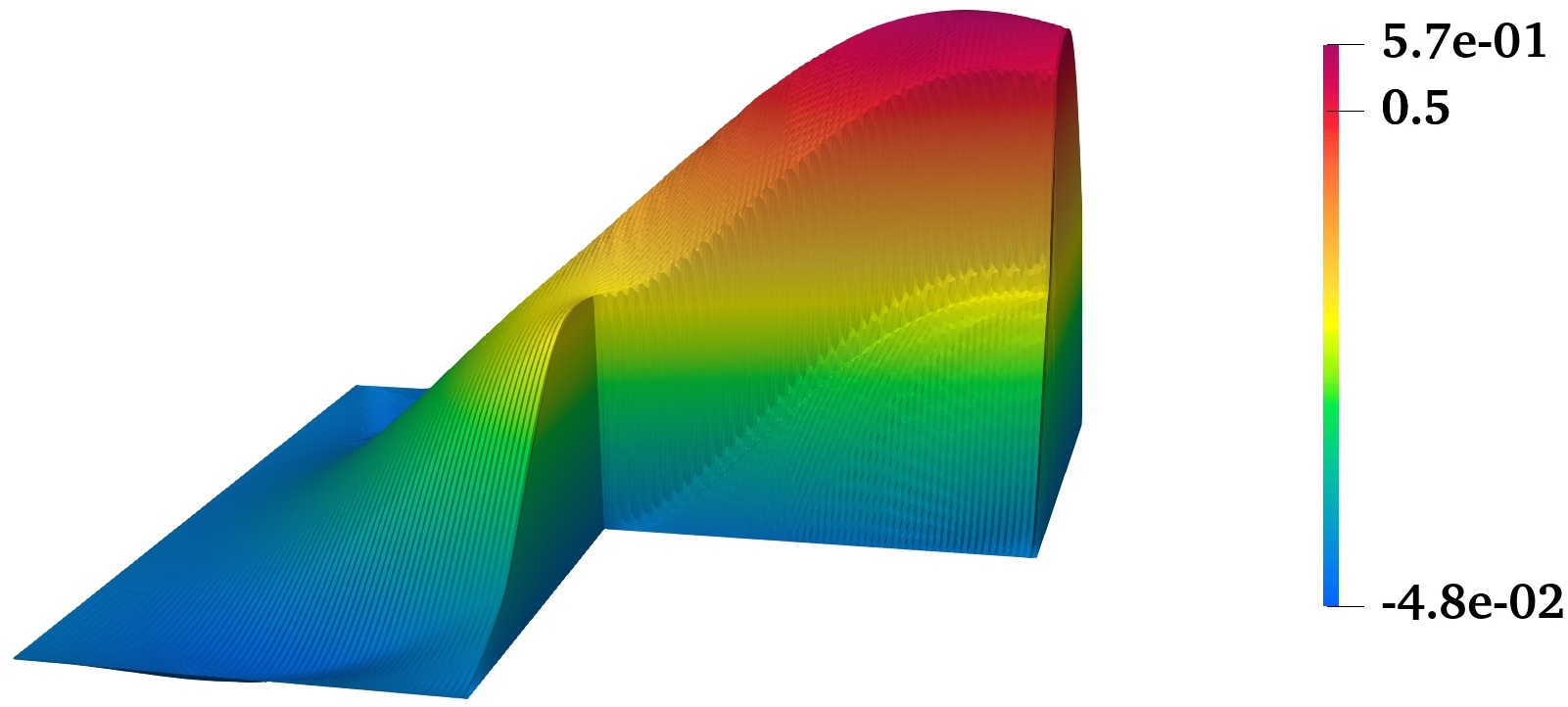}
    \caption{Example~\ref{ex:corner_singularity}: Numerical solution obtained with algebraic stabilization and the MUAS limiter on Grid~4 with $49{,}655$ degrees of freedom.}
    \label{fig:corner_singularity}
\end{figure}

Consider Eq.~\eqref{eq:cdr_eqn} on the L-shaped domain
\[
\Omega = (0,1)^2 \setminus [0.5,1] \times [0,0.5].
\]
The coefficients are given by
\[
\varepsilon = 10^{-6}, \quad \bb = (3,1)^{\top}, \quad c = 1, \quad 
f = 100\, r(r - 0.5)\bigl(r - \sqrt{2}/2\bigr),
\]
where
\[
r = \sqrt{(x - 0.5)^2 + (y - 0.5)^2}.
\]
The entire boundary is subject to Dirichlet conditions, i.e., $\Gamma_{\mD} = \Gamma$, with boundary data $u_{\mD} = 0$ (see Fig.~\ref{fig:corner_singularity}).

An analytical solution to this problem is not known. Nevertheless, the solution exhibits several characteristic features:
\begin{itemize}
    \item a corner singularity at the re-entrant point $(0.5,0.5)$,
    \item an interior layer aligned with the convection direction,
    \item boundary layers along the top boundary $y = 1$ and along the vertical segments
    \[
    \left\lbrace (0.5,y) : 0 < y \leq 0.5 \right\rbrace, \quad
    \left\lbrace (1,y) : 0.5 < y \leq 1 \right\rbrace.
    \]
\end{itemize}

\begin{figure}[tbp]
    \centering
    \includegraphics[width=0.3\linewidth]{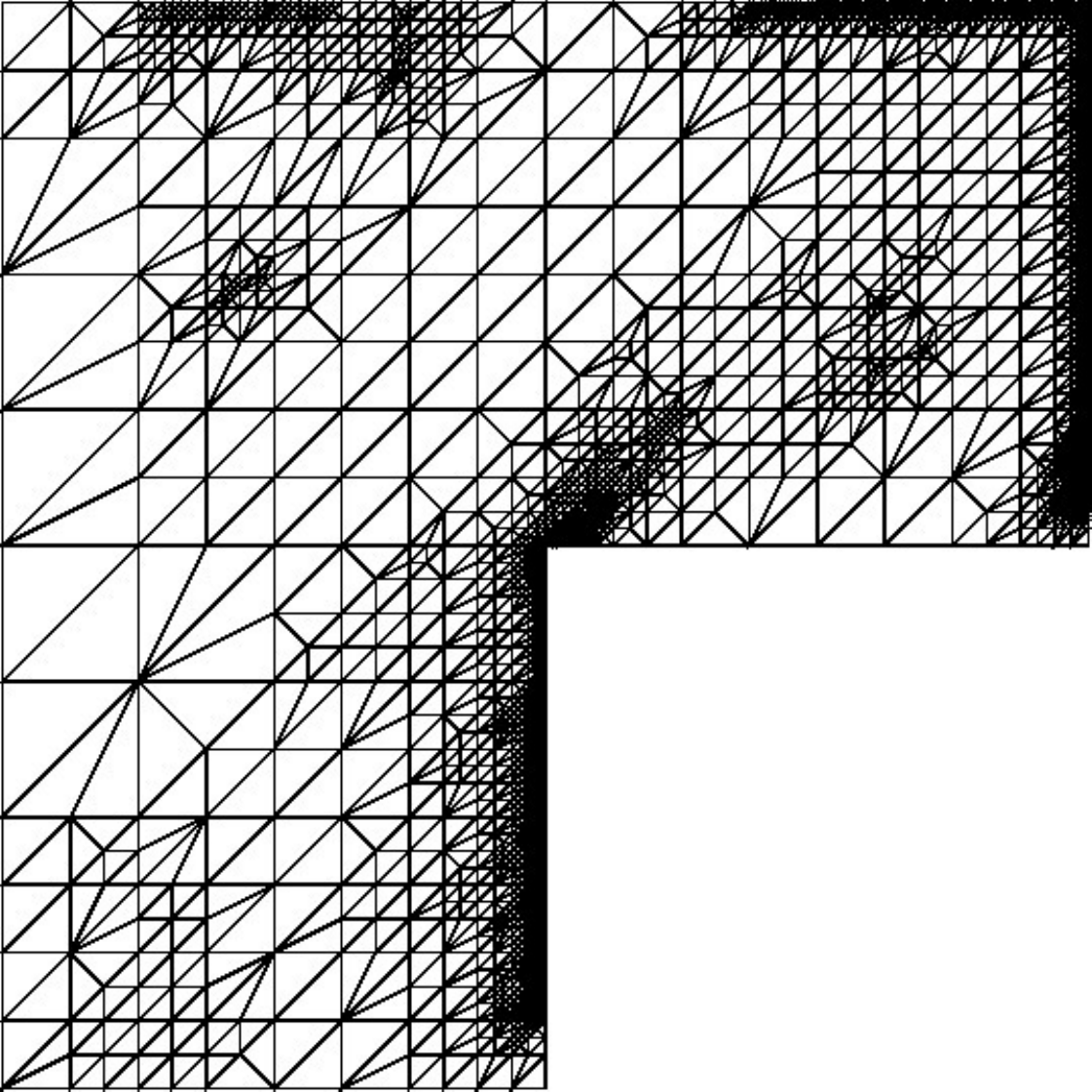}\hspace*{1em}
    \includegraphics[width=0.3\linewidth]{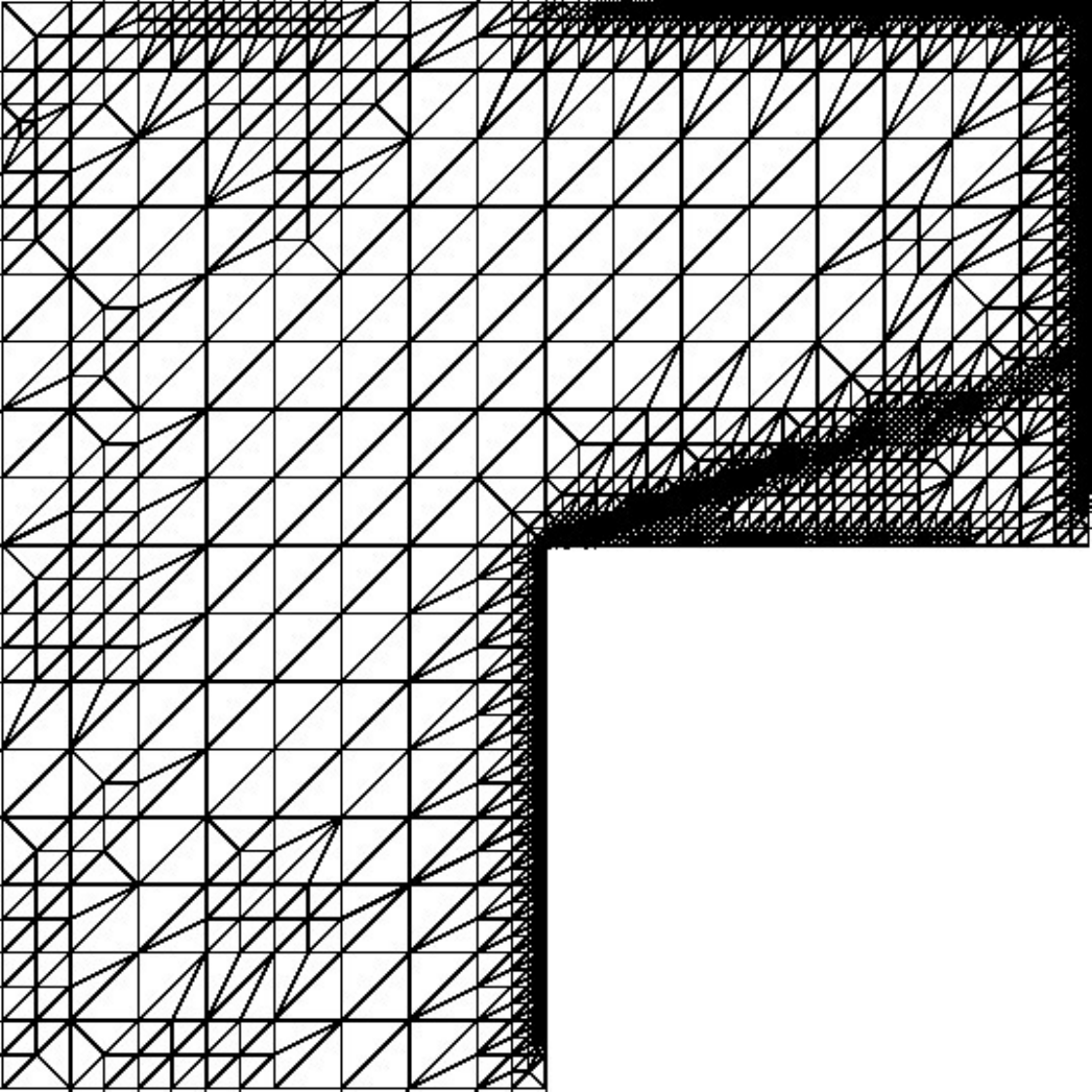}
    \caption{Example~\ref{ex:corner_singularity}: Adaptive grids for the AFC methods with $\#\mathrm{dofs} \approx 2.5 \times 10^5$. 
    BJK limiter (left) and MC limiter (right).}
    \label{fig:adaptive_afc}
\end{figure}

\begin{figure}[tbp]
    \centering
    \includegraphics[width=0.3\linewidth]{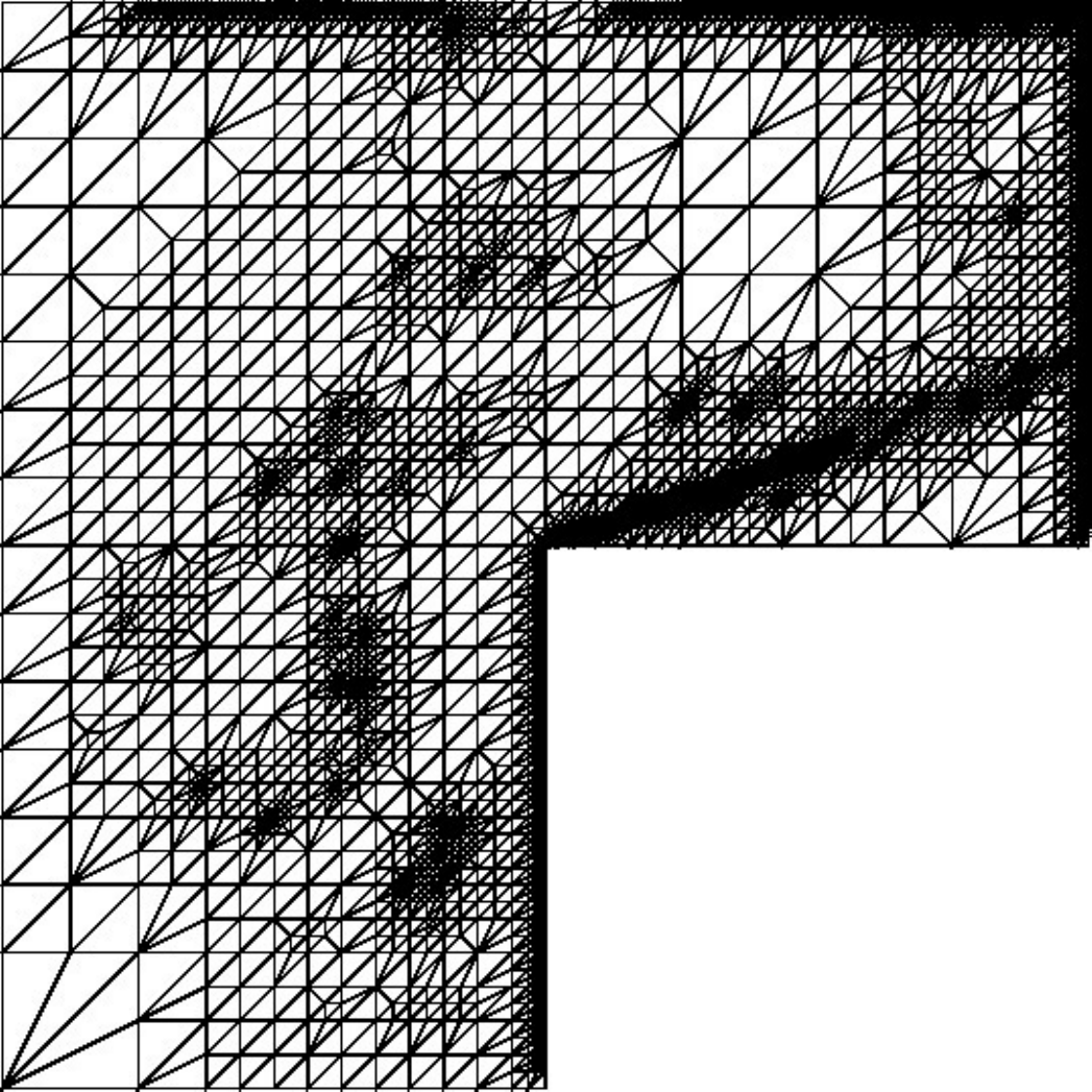}\hspace*{1em}
    \includegraphics[width=0.3\linewidth]{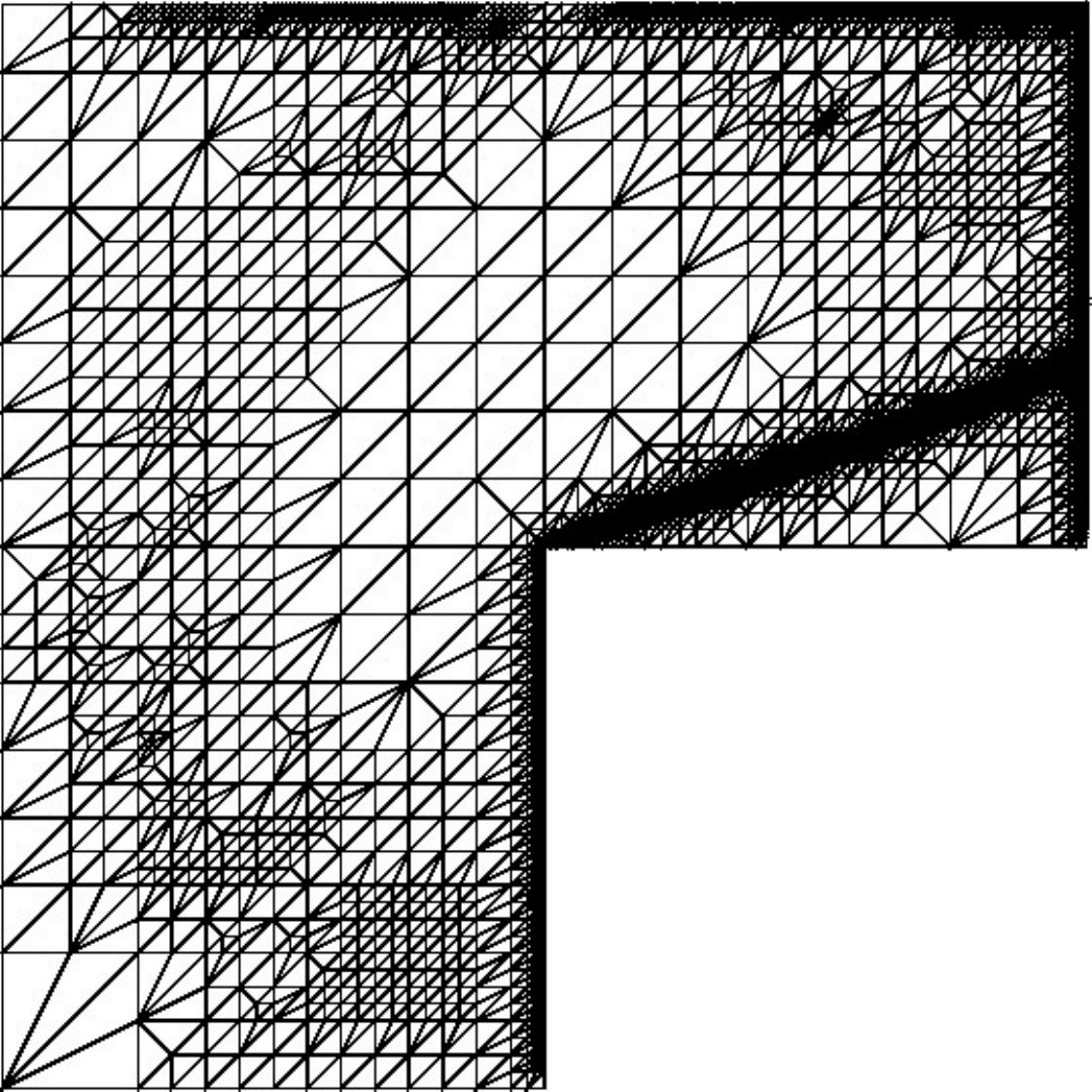}\hspace*{1em}
    \includegraphics[width=0.3\linewidth]{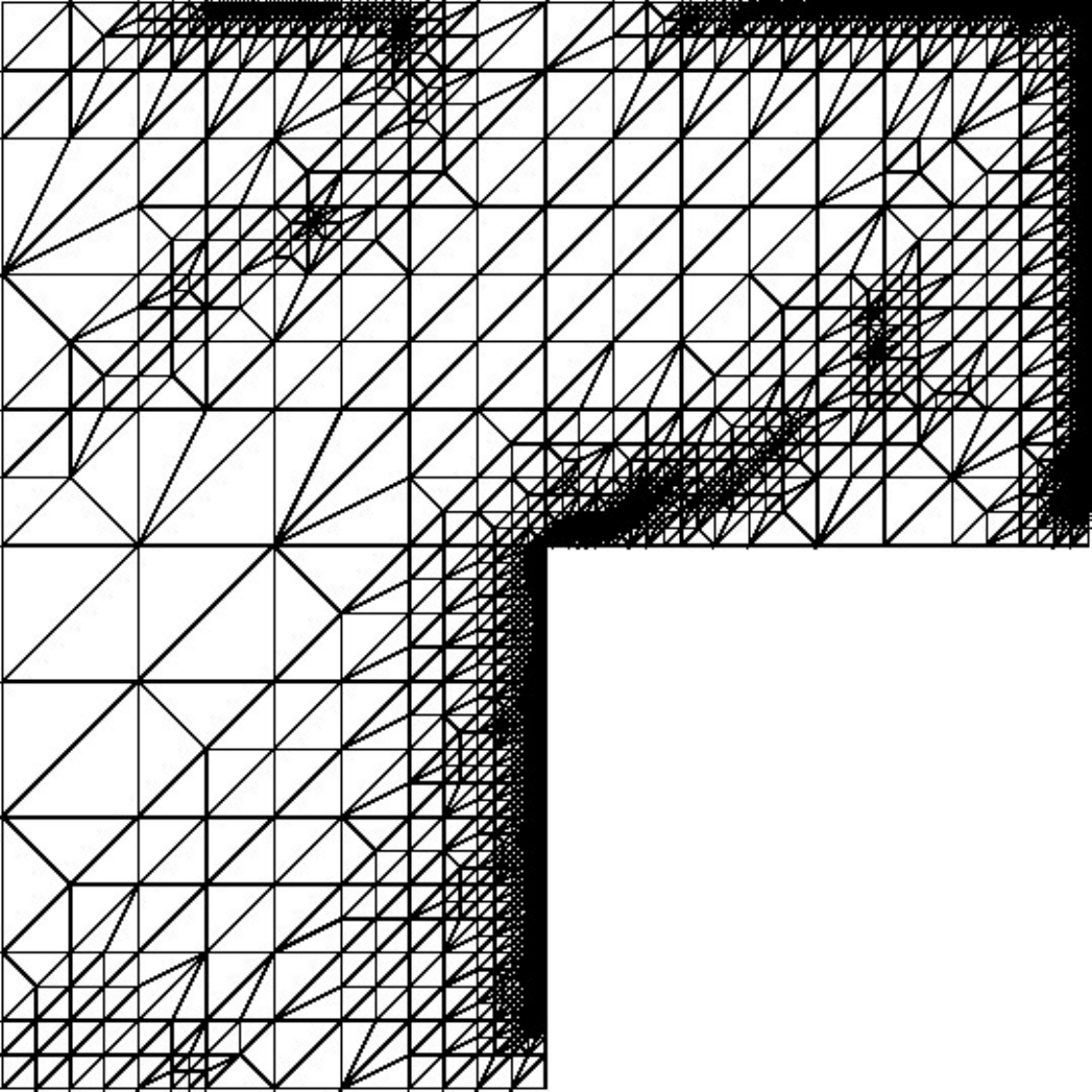}
    \caption{Example~\ref{ex:corner_singularity}: Adaptive grids for the algebraically stabilized methods with $\#\mathrm{dofs} \approx 2 \times 10^5$. 
    MUAS (left), SMUAS (middle), and BBK (right).}
    \label{fig:adaptive_as}
\end{figure}

Due to the presence of the re-entrant corner, the solution does not belong to $\mH^2(\Omega)$. 
The adaptively refined grids obtained with all limiters are shown in Fig.~\ref{fig:adaptive_afc} and Fig.~\ref{fig:adaptive_as} for $\#\mathrm{dofs} \approx 2.5 \times 10^5$. 
Among the AFC schemes, the BJK limiter fails to adequately refine the interior layer; instead, most refinement is concentrated near the boundary layers and the re-entrant corner. 
In contrast, the MC limiter performs better, resolving both the boundary layers and the interior layer more effectively.

For the algebraically stabilized schemes, the MUAS and SMUAS methods provide the best grid adaptation, with the SMUAS method performing slightly better overall. 
The BBK method struggles to adequately refine the interior region. 
In comparison with the AFC schemes, the SMUAS method appears to be the most promising in terms of adaptive grid quality.

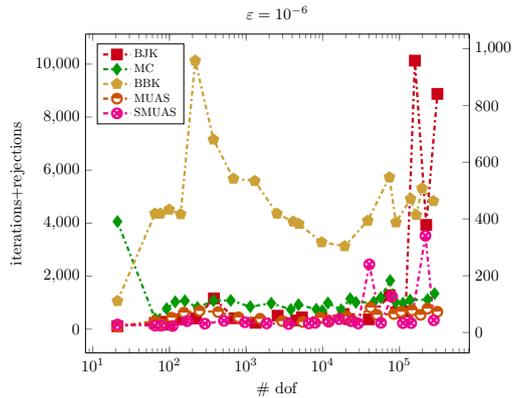
\begin{figure}[tbp]
\centerline{\begin{tikzpicture}[scale=0.6]
\begin{semilogxaxis}[
    scaled y ticks=false, 
    axis y line*=left, 
    axis x line=none, 
    xlabel = $\#\ \mathrm{dof}$, 
    ylabel = iterations+rejections,
    legend pos=north west, 
    legend cell align={left}, 
    title={$\varepsilon=10^{-6}$}, 
    legend style={nodes={scale=0.75, transform shape}}
]
\addplot[color=crimson,  mark=square*, line width = 0.5mm, dashdotted,,mark options = {scale= 1.5, solid}]
coordinates{( 21.0 , 110.0 )( 65.0 , 230.0 )( 78.0 , 278.0 )( 101.0 , 169.0 )( 146.0 , 370.0 )( 225.0 , 406.0 )( 388.0 , 1161.0 )( 721.0 , 412.0 )( 1368.0 , 241.0 )( 2667.0 , 508.0 )( 4672.0 , 330.0 )( 5555.0 , 446.0 )( 10011.0 , 432.0 )( 19595.0 , 572.0 )( 23058.0 , 380.0 )( 41987.0 , 371.0 )( 78404.0 , 1281.0 )( 95654.0 , 626.0 )( 127656.0 , 670.0 )( 168861.0 , 10130.0 )( 237826.0 , 3933.0 )( 333187.0 , 8875.0 )};
\label{BL_BJK_ITE}
\end{semilogxaxis}
\begin{semilogxaxis}[
    scaled y ticks=false, 
    axis y line*=right, 
    xlabel={$\#\ \mathrm{dof}$}, 
    legend pos=north west, 
    legend cell align={left}, 
    legend style={nodes={scale=0.75, transform shape}}
]
\addlegendimage{/pgfplots/refstyle=BL_BJK_ITE}
\addlegendentry{BJK} 
\addplot[color=dark_green,  mark=diamond*, line width = 0.5mm, dashdotted,,mark options = {scale= 1.5, solid}]
coordinates{( 21.0 , 391.0 )( 65.0 , 48.0 )( 92.0 , 84.0 )( 119.0 , 108.0 )( 159.0 , 112.0 )( 232.0 , 88.0 )( 376.0 , 112.0 )( 631.0 , 113.0 )( 1134.0 , 91.0 )( 2130.0 , 103.0 )( 3837.0 , 81.0 )( 4826.0 , 98.0 )( 8239.0 , 82.0 )( 9877.0 , 77.0 )( 11679.0 , 104.0 )( 17387.0 , 80.0 )( 22770.0 , 119.0 )( 27062.0 , 106.0 )( 46170.0 , 107.0 )( 57288.0 , 123.0 )( 75805.0 , 183.0 )( 89365.0 , 107.0 )( 111346.0 , 103.0 )( 137623.0 , 114.0 )( 232856.0 , 116.0 )( 287371.0 , 137.0 )};
\addlegendentry{MC} 
%\addplot[color=amethyst,  mark=triangle*, line width = 0.5mm, dashdotted,,mark options = {scale= 1.5, solid}]
%coordinates{( 21.0 , 61.0 )( 65.0 , 101.0 )( 100.0 , 112.0 )( 132.0 , 91.0 )( 212.0 , 92.0 )( 326.0 , 135.0 )( 567.0 , 119.0 )( 1014.0 , 123.0 )( 1869.0 , 111.0 )( 3601.0 , 114.0 )( 5759.0 , 109.0 )( 6967.0 , 136.0 )( 11247.0 , 125.0 )( 14783.0 , 128.0 )( 23674.0 , 136.0 )( 27989.0 , 130.0 )( 47264.0 , 153.0 )( 58675.0 , 124.0 )( 71260.0 , 125.0 )( 111469.0 , 112.0 )( 136109.0 , 150.0 )( 217705.0 , 140.0 )( 277705.0 , 153.0 )};
%\addlegendentry{BBK} 
\addplot[color=gold,  mark=pentagon*, line width = 0.5mm, dashdotted,,mark options = {scale= 1.5, solid}]
coordinates{( 21.0 , 111.0 )( 65.0 , 419.0 )( 76.0 , 419.0 )( 99.0 , 433.0 )( 140.0 , 417.0 )( 217.0 , 958.0 )( 375.0 , 680.0 )( 684.0 , 542.0 )( 1297.0 , 534.0 )( 2527.0 , 419.0 )( 4111.0 , 391.0 )( 4871.0 , 383.0 )( 9679.0 , 318.0 )( 19301.0 , 304.0 )( 38606.0 , 394.0 )( 74027.0 , 547.0 )( 90198.0 , 388.0 )( 139752.0 , 471.0 )( 165769.0 , 414.0 )( 197043.0 , 508.0 )( 277667.0 , 463.0 )};
\addlegendentry{BBK} 
\addplot[color=carrotorange,  mark=halfcircle*, line width = 0.5mm, dashdotted,,mark options = {scale= 1.5, solid}]
coordinates{( 21.0 , 22.0 )( 65.0 , 40.0 )( 78.0 , 30.0 )( 107.0 , 52.0 )( 154.0 , 70.0 )( 247.0 , 76.0 )( 432.0 , 72.0 )( 797.0 , 53.0 )( 1523.0 , 48.0 )( 2960.0 , 40.0 )( 5461.0 , 37.0 )( 9274.0 , 52.0 )( 11675.0 , 43.0 )( 20521.0 , 56.0 )( 24638.0 , 41.0 )( 42128.0 , 88.0 )( 50873.0 , 62.0 )( 84844.0 , 65.0 )( 105247.0 , 73.0 )( 146340.0 , 78.0 )( 188672.0 , 62.0 )( 232559.0 , 84.0 )( 312534.0 , 73.0 )};
\addlegendentry{MUAS} 
\addplot[color=magenta,  mark=otimes, line width = 0.5mm, dashdotted,,mark options = {scale= 1.5, solid}]
coordinates{( 21.0 , 28.0 )( 65.0 , 24.0 )( 78.0 , 24.0 )( 110.0 , 23.0 )( 170.0 , 40.0 )( 290.0 , 31.0 )( 514.0 , 41.0 )( 966.0 , 36.0 )( 1850.0 , 32.0 )( 3594.0 , 30.0 )( 6580.0 , 32.0 )( 7921.0 , 34.0 )( 11968.0 , 38.0 )( 16407.0 , 47.0 )( 23400.0 , 40.0 )( 28991.0 , 32.0 )( 40480.0 , 240.0 )( 57397.0 , 34.0 )( 80062.0 , 127.0 )( 110893.0 , 33.0 )( 142063.0 , 33.0 )( 218299.0 , 341.0 )( 280439.0 , 43.0 )};
\addlegendentry{SMUAS} 
\end{semilogxaxis}
\end{tikzpicture}}
\caption{Example~\ref{ex:corner_singularity}: Number of iterations and rejections. 
The left $y$-axis corresponds to the BJK limiter, whereas the right $y$-axis corresponds to the other limiters.}
\label{fig:iteration_corner_singularity}
\end{figure}

Figure~\ref{fig:iteration_corner_singularity} displays the number of iterations and rejections in the nonlinear iterative loop. 
The left $y$-axis corresponds to the BJK limiter, while the right $y$-axis corresponds to the remaining limiters. 
Clearly, the BJK limiter performs the worst among the tested methods, whereas the MUAS and SMUAS limiters achieve the best results.  
Similar observations for the BJK limiter under uniform refinement were reported in~\cite{JJ19}.

Table~\ref{tab:corner_time} reports the computing time and the number of adaptive loops required by the adaptive algorithm. 
We observe that the MUAS method is the most efficient, whereas the SMUAS method is the least efficient. 
In fact, SMUAS requires approximately twice as much computing time as the next most inefficient method (BJK).  
Hence, although the SMUAS method requires the fewest adaptive loops, it is inefficient in terms of overall computing time.

\begin{table}[tbp]
\centering
\begin{tabular}{|c|c|c|}
\hline
\textbf{Limiter} & \textbf{Time (s)} & \textbf{Number of Adaptive Loops} \\
\hline\hline
BJK & 1947 & 22 \\
MC & 150 & 26 \\
BBK & 265 & 21 \\
MUAS & 136 & 23 \\
SMUAS &  3345 & 23 \\
\hline
\end{tabular}
\caption{Example~\ref{ex:corner_singularity}: Computing time (in seconds) and number of adaptive loops for the adaptive algorithm.}
\label{tab:corner_time}
\end{table}

Overall, in terms of adaptive grid quality, the SMUAS and MC limiters perform best, whereas in terms of nonlinear iteration efficiency, the MUAS and SMUAS methods yield the most favorable results. Regarding computing time, MUAS and the MC limiter take the least time. In both respects, the BJK limiter is the weakest performer among the tested methods.

\subsection{Solution with Different Regimes}\label{ex:volker_example}
This example is taken from \cite{BJK25}. The domain contains regions where convection dominates, 
regions where reaction dominates, and regions where diffusion is the primary mechanism.

Consider Eq.~\eqref{eq:cdr_eqn} on $\Omega = (0,1)^2$ with $\varepsilon = 10^{-6}$, $f = 0$, 
and a convective field given by
\[
\bb(x,y)
  = r\,\psi(r)
    \begin{bmatrix}
      y + 27/16 \\[0.3em]
      1 - x
    \end{bmatrix},
  \qquad
  r = \left( (x-1)^2 + \left( y + \tfrac{27}{16} \right)^2 \right)^{1/2},
\]
where
\[
\psi(r)
 = \begin{cases}
    \dfrac{1}{r}, & |r-r_0|\le r_d,\\[0.7em]
    \dfrac{1}{r_0+r_d}\exp\!\bigl(-1000\,(r-(r_0+r_d))^2\bigr), & r>r_0+r_d,\\[0.7em]
    \dfrac{1}{r_0-r_d}\exp\!\bigl(-1000\,(r-(r_0-r_d))^2\bigr), & r<r_0-r_d,
   \end{cases}
\]
with $r_0 = 17/8$ and $r_d = \sqrt{1217}/16 - r_0$.

The reaction coefficient is defined analogously:
\[
c =
  \begin{cases}
    1, & |r-r_0|\le r_d,\\[0.7em]
    \exp\!\bigl(-100\,(r-(r_0+r_d))^2\bigr), & r>r_0+r_d,\\[0.7em]
    \exp\!\bigl(-100\,(r-(r_0-r_d))^2\bigr), & r<r_0-r_d.
  \end{cases}
\]

\begin{figure}[tbp]
    \centering
    \includegraphics[width=0.4\linewidth]{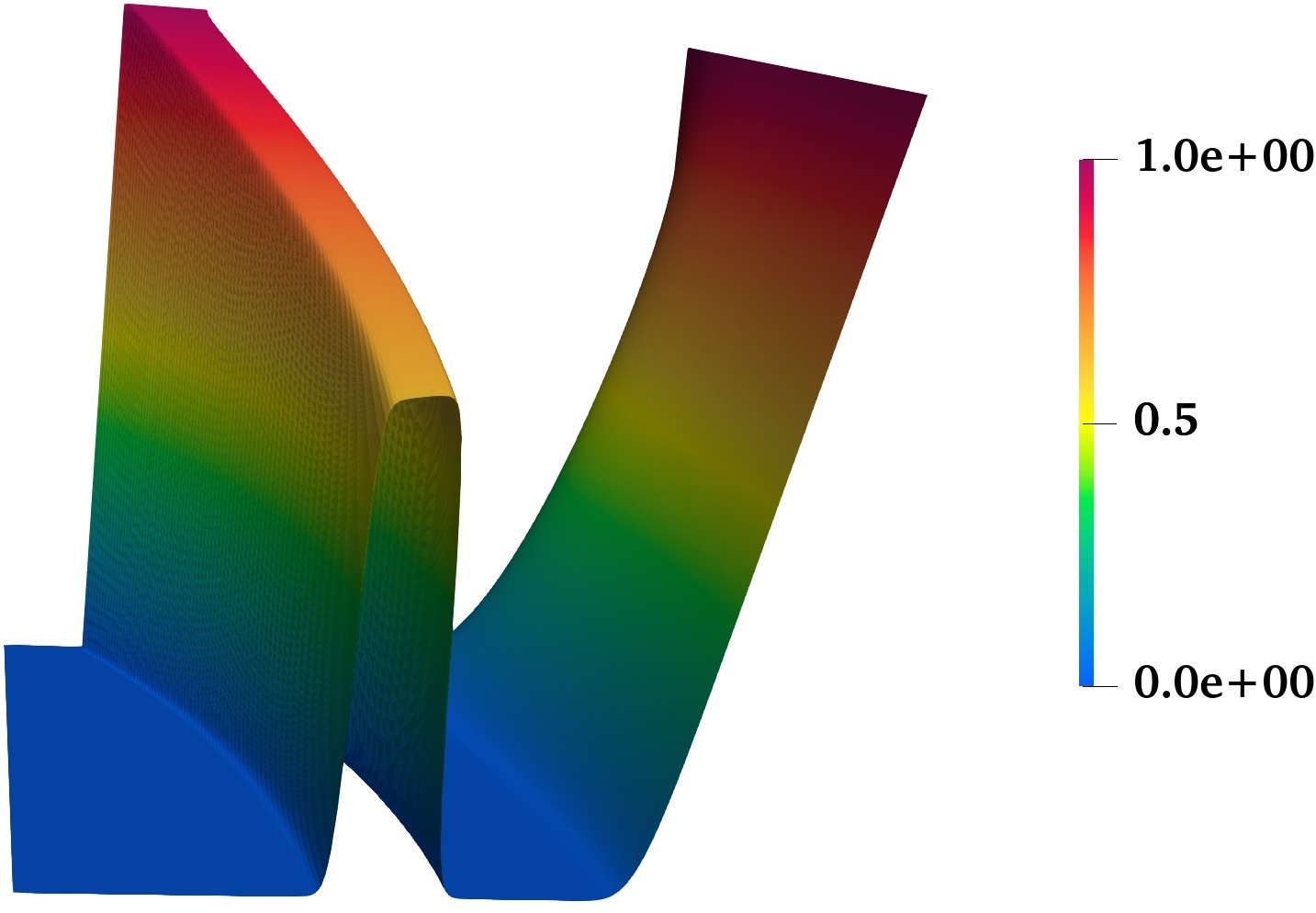}
    \caption{Example~\ref{ex:corner_singularity}: Numerical solution with algebraic stabilisation and MUAS limiter on Grid~1 with $66049$ degrees of freedom.}
    \label{fig:volker_example}
\end{figure}

Both the convection and reaction fields attain their largest values in the region
\[
\{(x,y)\in\Omega : |r-r_0|\le r_d\},
\]
and decay exponentially away from it, with the reaction field decaying more slowly 
(see Fig.~\ref{fig:volker_example}).

Homogeneous Neumann boundary conditions are imposed along the boundary $(1,y)$ for $y\in(0,1)$.
Dirichlet boundary conditions on the remaining boundary are prescribed as
\[
u_{\mD} = 
\begin{cases}
0, & y=0,\\[0.3em]
0, & x=0,\; y\in[0,0.125),\\[0.3em]
1, & x=0,\; y\in[0.125,0.25],\\[0.3em]
0, & x=0,\; y\in(0.25,1),\\[0.3em]
1, & y=1.
\end{cases}
\]

All simulations are performed on Grid~1.  
To assess the quality of the numerical solution, we compare solution values along several cutlines:
\[
x = 0.5, \qquad x = 1, \qquad x = y, \qquad x = 1-y.
\]

\begin{figure}[tbp]
\centering
\includegraphics[width=0.42\textwidth]{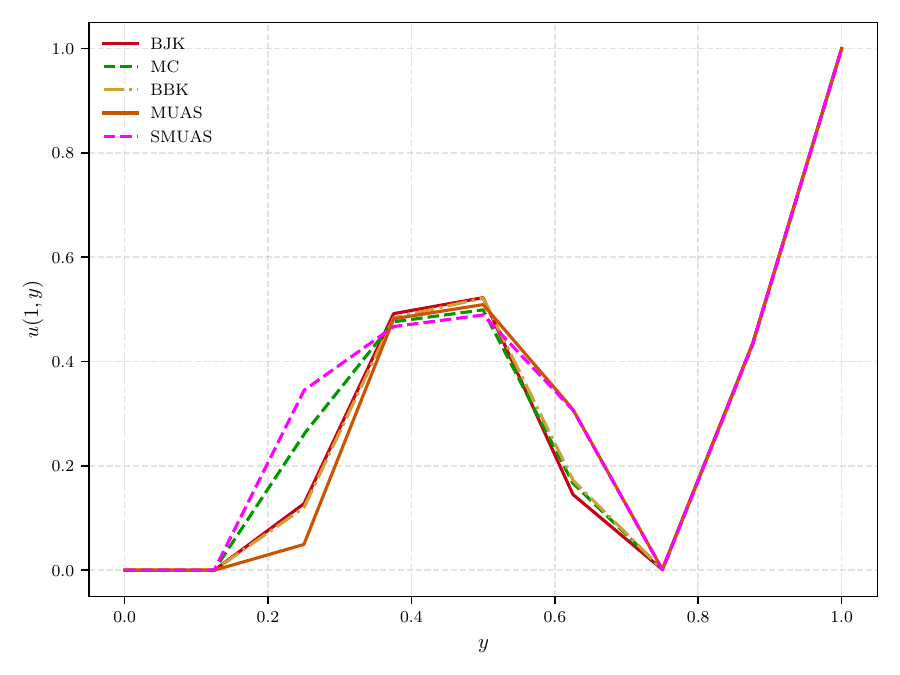}
\includegraphics[width=0.42\textwidth]{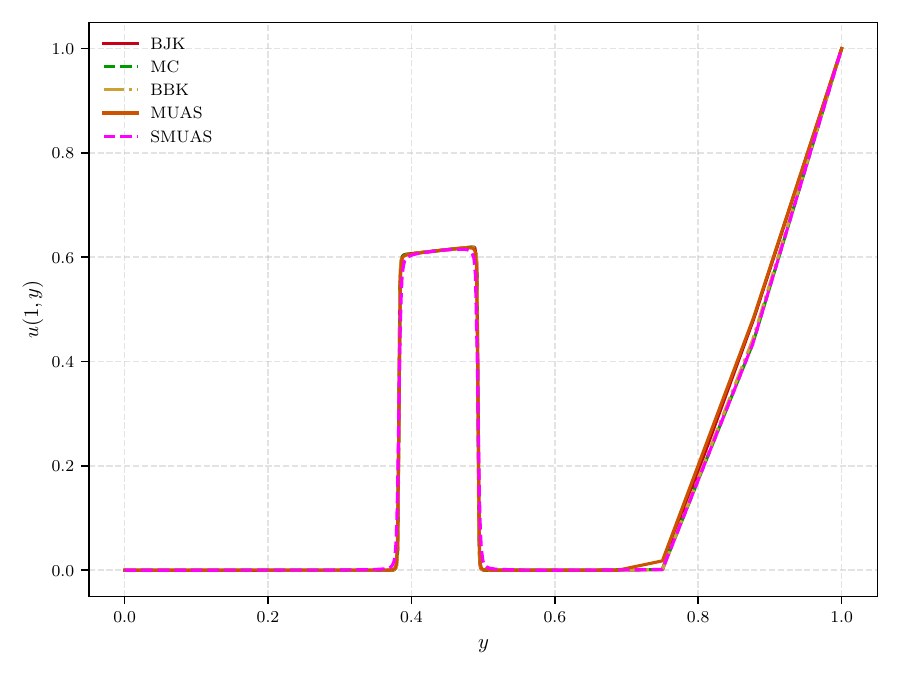}
\caption{Example~\ref{ex:volker_example}: Solution along the cutline $x=1$ with 
approximately $100$ \# dof (left) and $125{,}000$ \# dof (right).}
\label{fig:volker_example_cutline_1}
\end{figure}

In Fig.~\ref{fig:volker_example_cutline_1} we present the results along the cutline $x = 1$ 
for two degrees-of-freedom levels: approximately $100$ (left), and $125{,}000$ (right).  The solution obtained with $100$ \# dof is the solution on the first adaptively refined grid. Here we see that all the limiters capture the basic features of the example.  As the mesh becomes fine,  we observe that the solution is robust with respect to the choice of limiter, 
and only minimal (if any) differences are visible in the solution profiles.

\begin{figure}[tbp]
\centering
\includegraphics[width=0.42\textwidth]{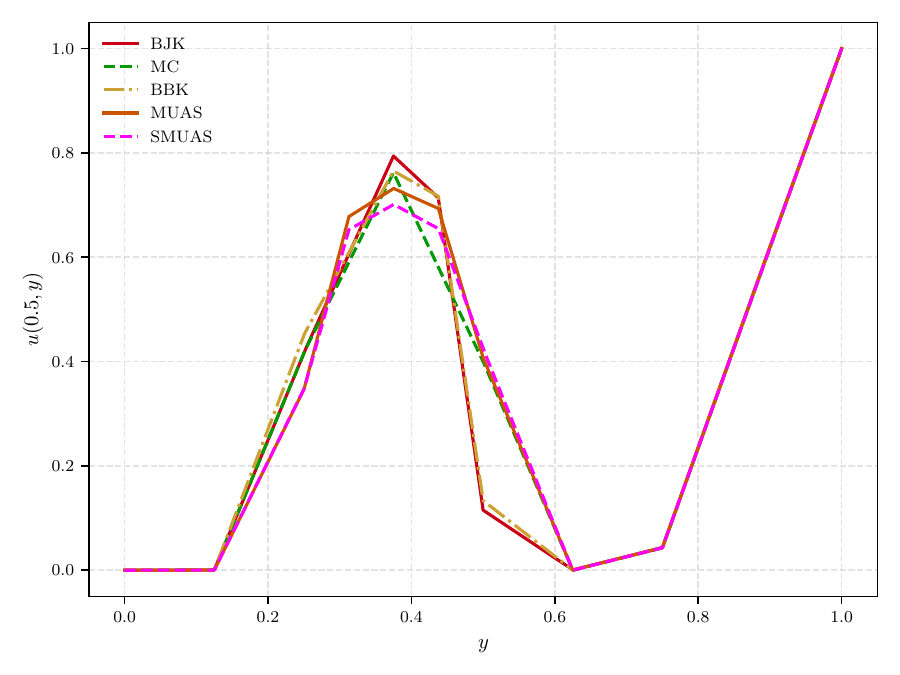}
\includegraphics[width=0.42\textwidth]{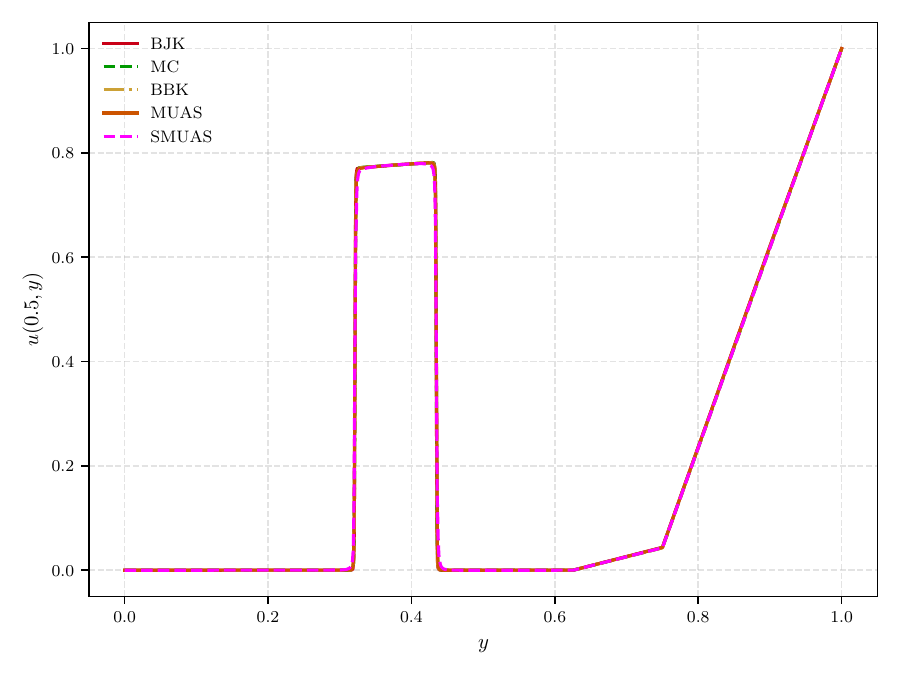}
\caption{Example~\ref{ex:volker_example}: Solution along the cutline $x=0.5$ with 
approximately $100$ \# dof (left) and $125{,}000$ \# dof (right).}
\label{fig:volker_example_cutline_05}
\end{figure}

\begin{figure}[tbp]
\centering
\includegraphics[width=0.42\textwidth]{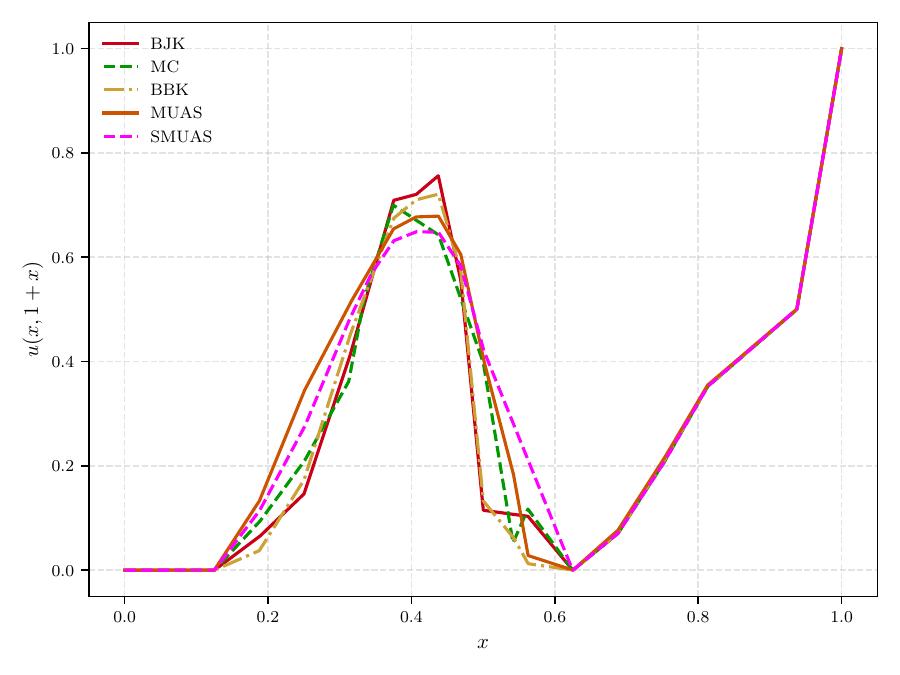}
\includegraphics[width=0.42\textwidth]{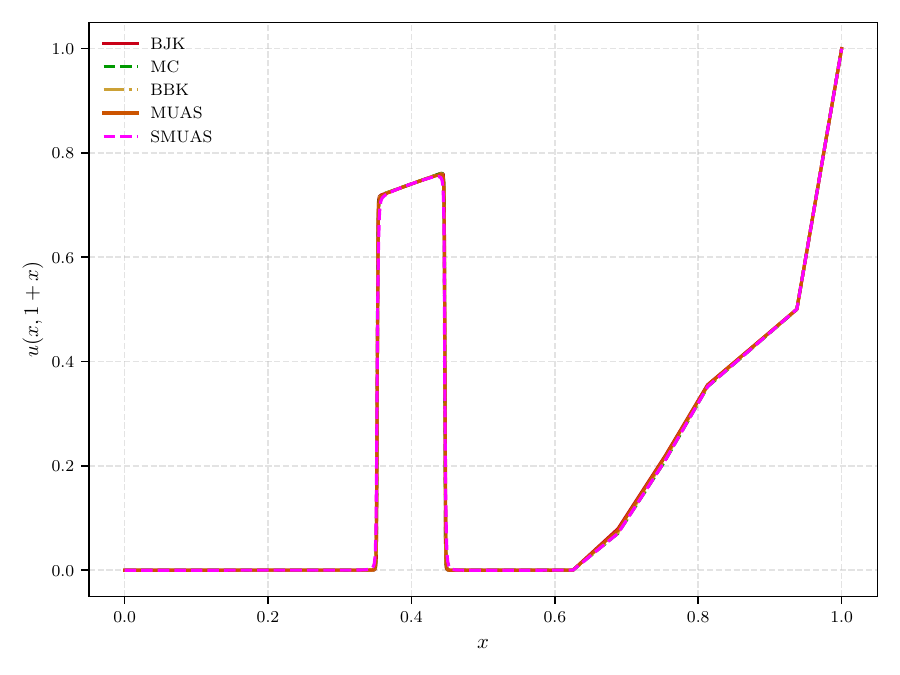}
\caption{Example~\ref{ex:volker_example}: Solution along the cutline $x=1-y$ with 
approximately $100$ \# dof (left) and $125{,}000$ \# dof (right).}
\label{fig:volker_example_cutline_x1_y}
\end{figure}

\begin{figure}[tbp]
\centering
\includegraphics[width=0.42\textwidth]{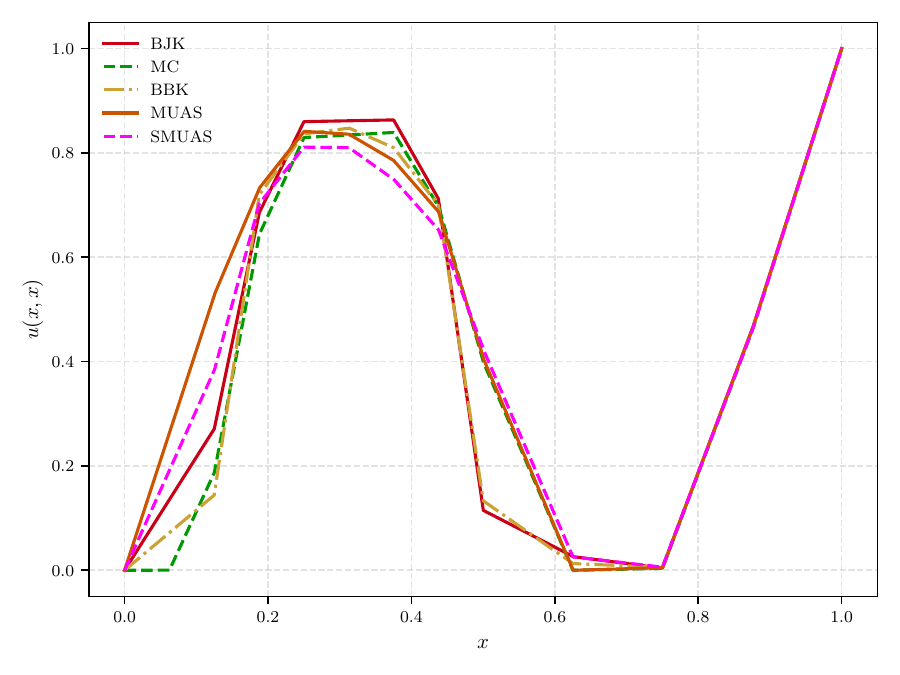}
\includegraphics[width=0.42\textwidth]{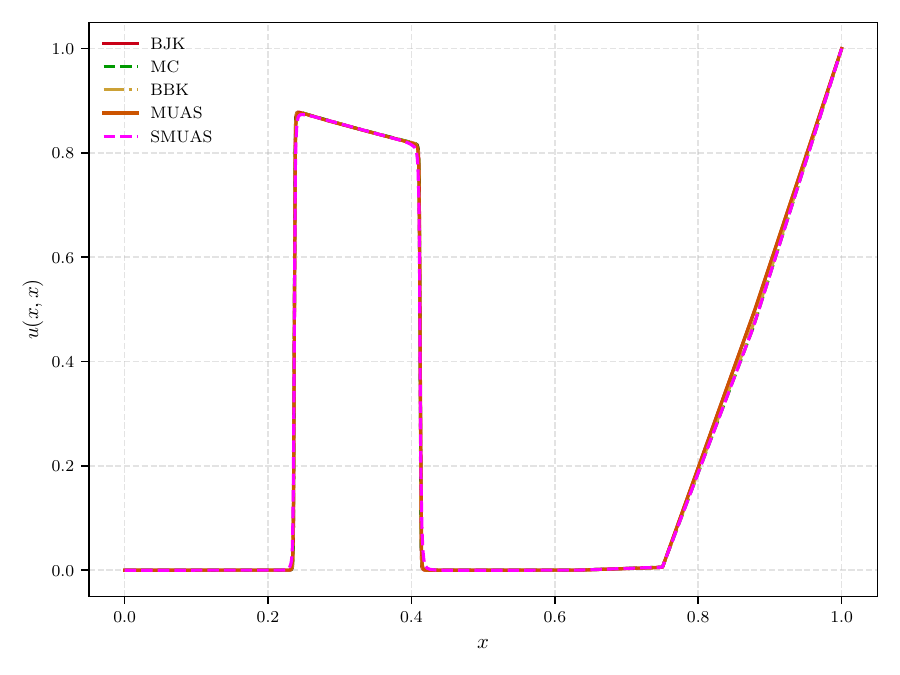}
\caption{Example~\ref{ex:volker_example}: Solution along the cutline $x=y$ with 
approximately $100$ \# dof (left) and $125{,}000$ \# dof (right).}
\label{fig:volker_example_cutline_x_y}
\end{figure}

Figs.~\ref{fig:volker_example_cutline_05}, \ref{fig:volker_example_cutline_x1_y}, and Fig.~\ref{fig:volker_example_cutline_x_y} shows the results along the cutline $x = 0.5$, $x=y$, and $x=1-y$, respectively.  
Here, similar observations can be made: the limiters produce nearly identical solution profiles, 
indicating that the method behaves robustly in all the configurations.

\begin{figure}[tbp]
    \centering{
    \includegraphics[width=0.3\linewidth]{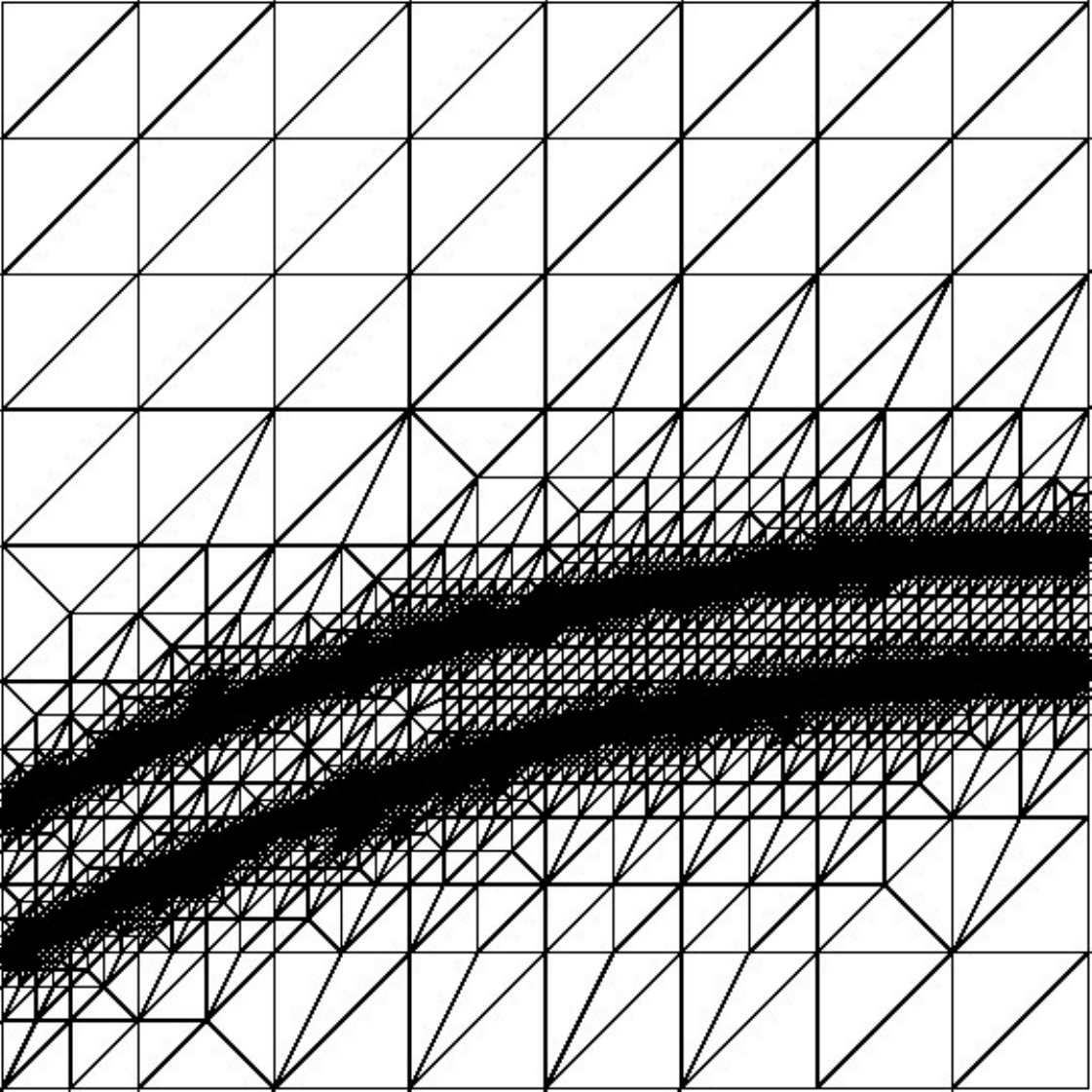}\hspace*{1em}
    \includegraphics[width=0.3\linewidth]{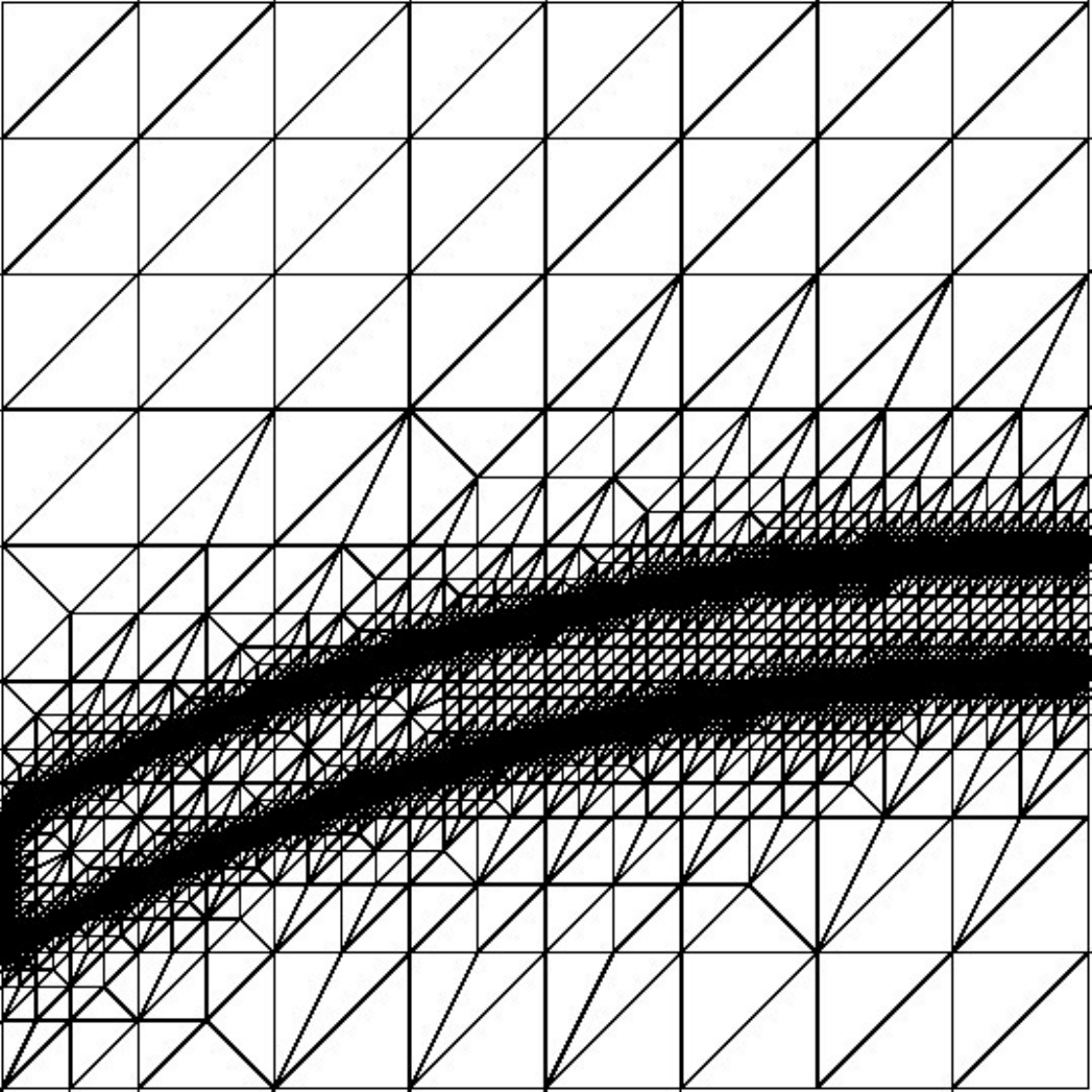}}
    \caption{Example~\ref{ex:volker_example}: Adaptive grids for the AFC methods with $\#\ \mathrm{dofs}\approx 2.5\times 10^5$.  BJK limiter (left) and MC limiter (right).}
    \label{fig:adaptive_afc_volker}
\end{figure}

\begin{figure}[tbp]
    \centering{
    \includegraphics[width=0.3\linewidth]{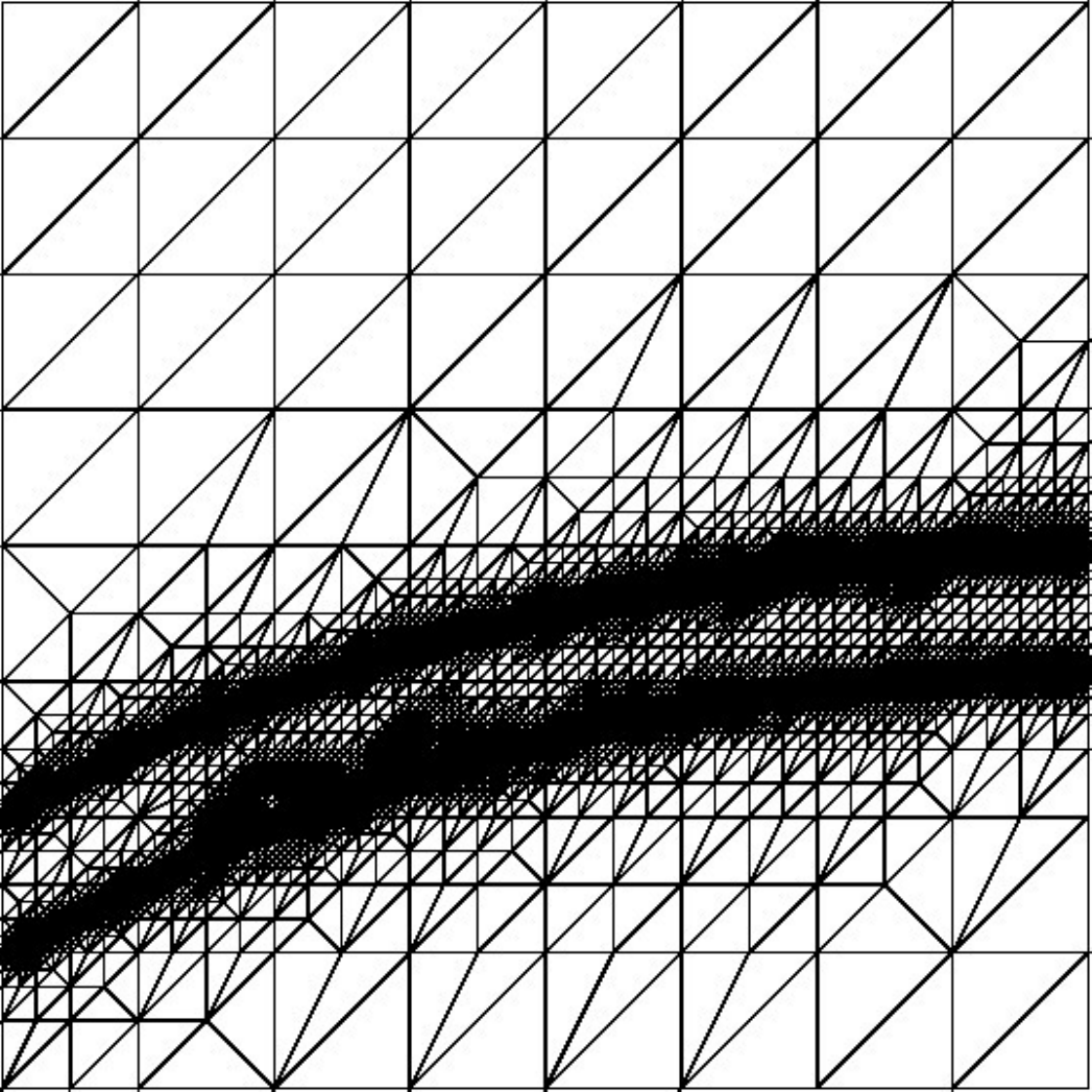}\hspace*{1em}
    \includegraphics[width=0.3\linewidth]{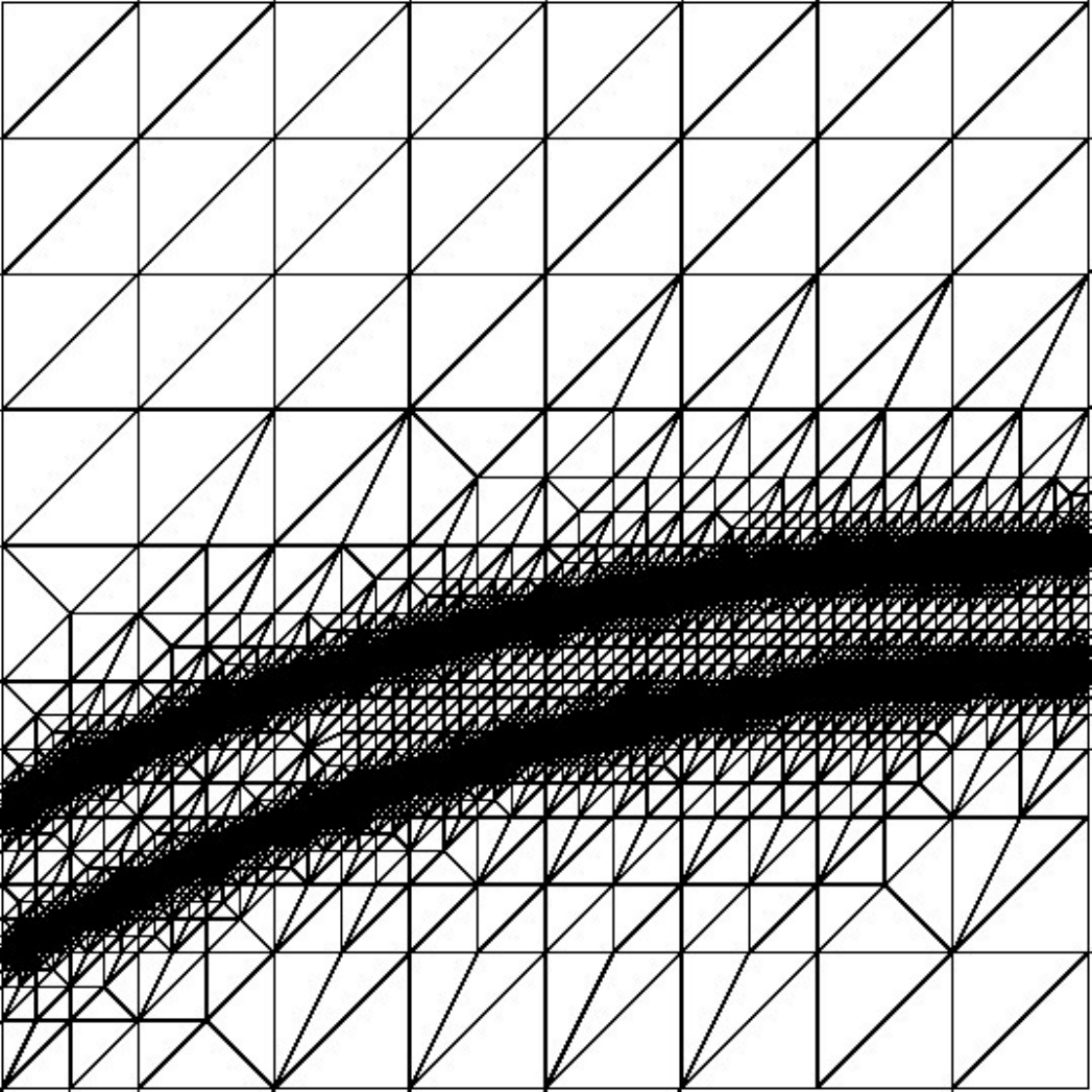}\hspace*{1em}
    \includegraphics[width=0.3\linewidth]{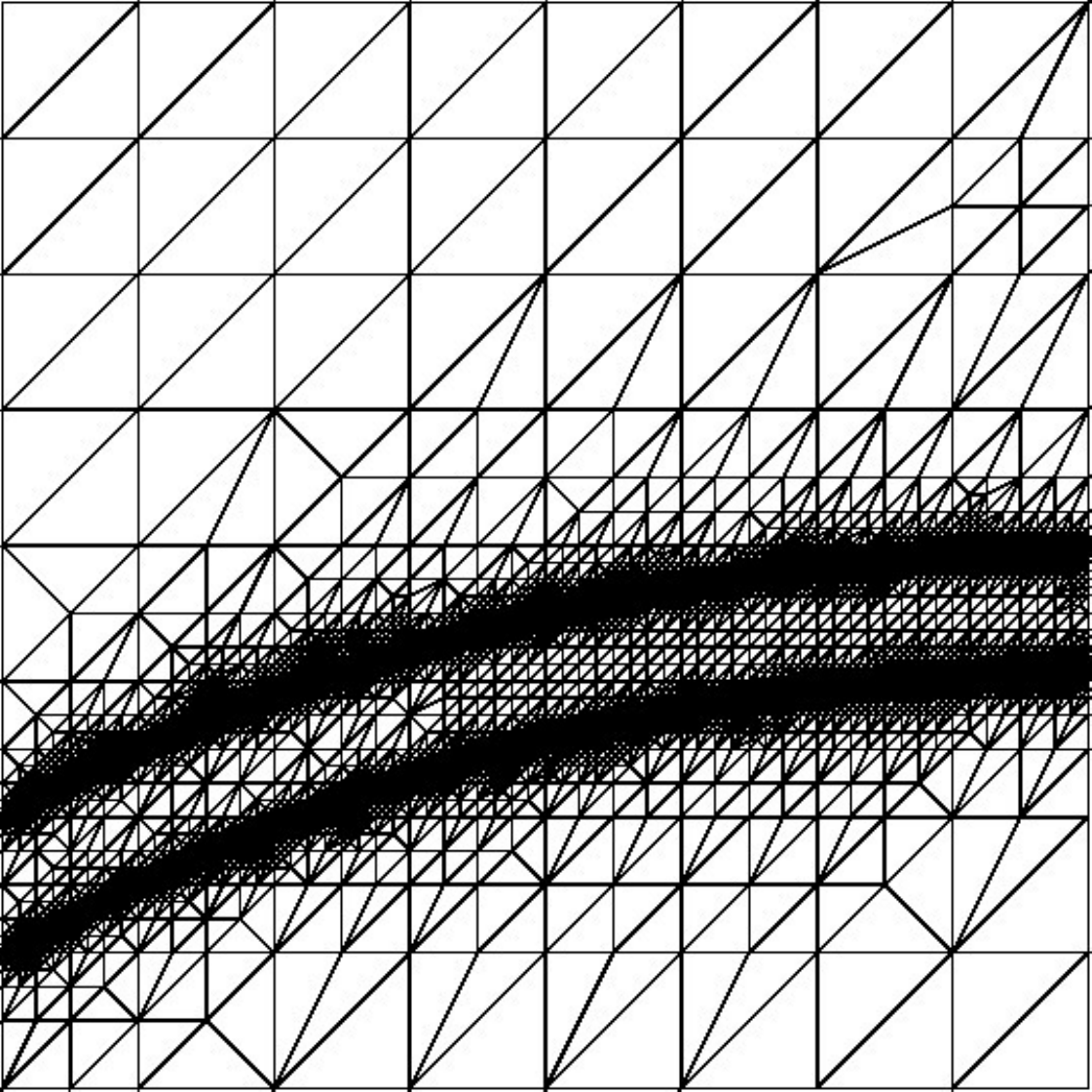}}
    \caption{Example~\ref{ex:volker_example}: Adaptive grids for the algebraically stabilized methods with $\#\ \mathrm{dofs}\approx 2\times 10^5$. 
    MUAS (left),  SMUAS (middle), BBK (right).}
    \label{fig:adaptive_as_volker}
\end{figure}

Lastly, in Figs.~\ref{fig:adaptive_afc_volker} and~\ref{fig:adaptive_as_volker}, we present the adaptive grids obtained using different limiters for approximately $2.5\times 10^{5}$ degrees of freedom. We observe that the layers captured by the MC limiter are sharper than those produced by the BJK limiter; however, the MC limiter tends to over-refine the mesh near the inlet boundary $(y=0)$. 

Among the algebraically stabilized methods, the SMUAS method appears to be the most promising, with the BBK method showing comparable performance. Overall, the quality of the adaptive grids is similar across the considered approaches, with the BJK and SMUAS methods providing the most favorable results.

\subsection{Hemker Problem}\label{ex:hemker}
The next example we consider is the standard Hemker problem, originally introduced in \cite{Hem96}. 
The computational domain is given by
$$
\Omega = (-3,9)\times(-3,3) \setminus \{(x,y)\in\mathbb{R}^2 : x^2+y^2\le 1\}.
$$
The convection field is prescribed as $\bb = (1,0)^{\top}$, while the reaction coefficient and the right-hand side
in Eq.~\eqref{eq:cdr_eqn} vanish, i.e.,
$$
c = f = 0.
$$
Dirichlet boundary conditions are imposed at $x=-3$, where $u_{\mD}=0$, and on the circular boundary, where $u_{\mD}=1$.
On all remaining boundaries, homogeneous Neumann boundary conditions are prescribed.

\begin{figure}[tbp]
\centerline{\includegraphics[width=0.5\textwidth]{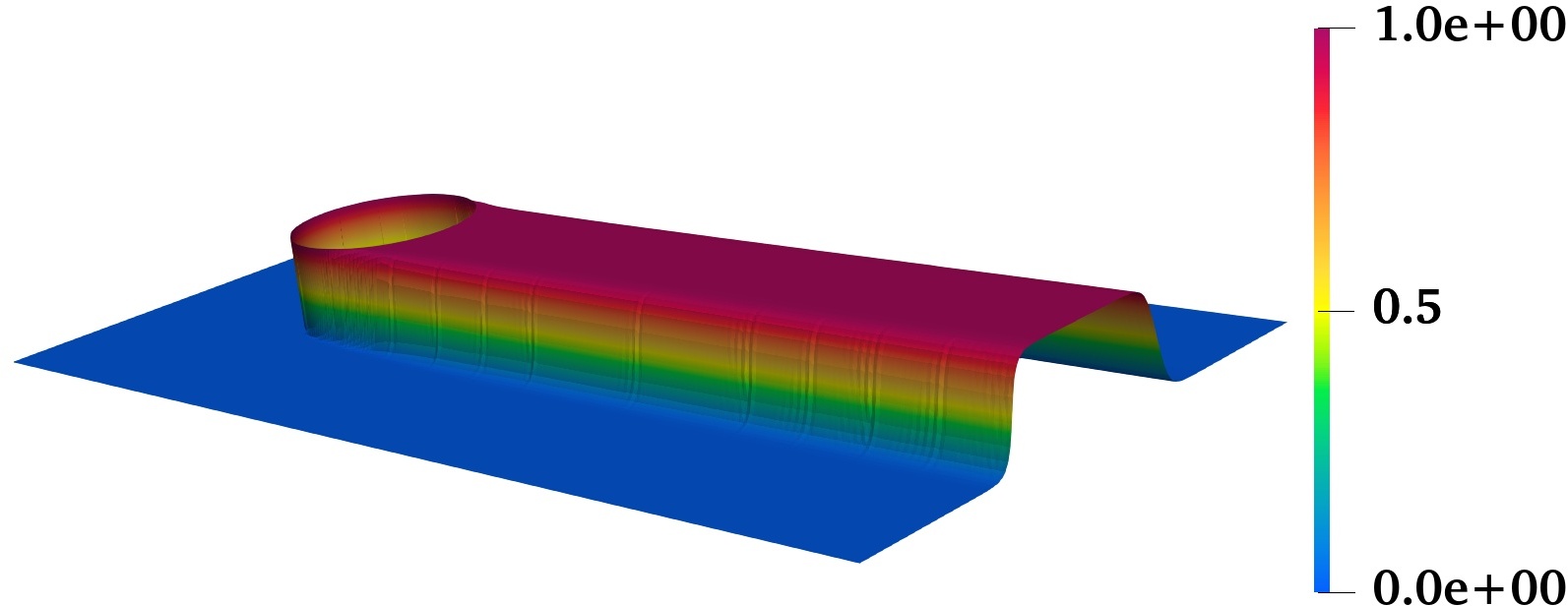}
\hspace*{1em}
\includegraphics[width=0.4\textwidth]{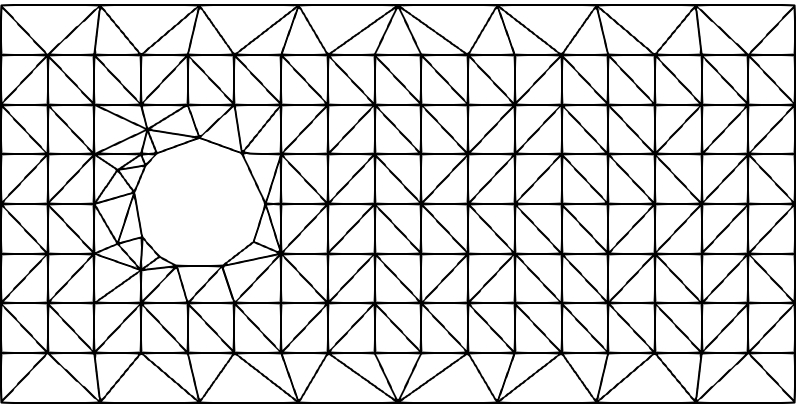}}
\caption{Example~\ref{ex:hemker}: Solution for $\varepsilon = 10^{-4}$ (left), computed with the MUAS limiter using $33{,}496$ degrees of freedom; initial grid (right), level~0.}
\label{fig:sol_hemker2d}
\end{figure}

This problem was studied comprehensively in \cite{ACJ11} for $\varepsilon = 10^{-4}$, where reference values for several
quantities of interest were reported. The same diffusion parameter is employed in our study.
Figure~\ref{fig:sol_hemker2d} (left) illustrates the solution, which takes values in the interval $[0,1]$,
while Figure~\ref{fig:sol_hemker2d} (right) shows the initial grid with $\#\mathrm{dof}=151$.

Although the exact solution to this problem is unknown, several quantities of interest are available.
In particular, the solution is known to be bounded between $0$ and $1$, and reference values exist for the
smearing of the internal layer along the cutline $x=4$.
To quantify these properties, we introduce the measures
$$
\mathrm{osc}(u_h) = \max(u_h) - \min(u_h) \approx 1,
$$
and
$$
\mathrm{smear}_{\mathrm{int}} = y_2 - y_1,
$$
(see Eq.~\eqref{eq:smear}). The cutline is discretised using $100{,}000$ equidistant intervals, and the discrete solution is evaluated at the corresponding nodes.

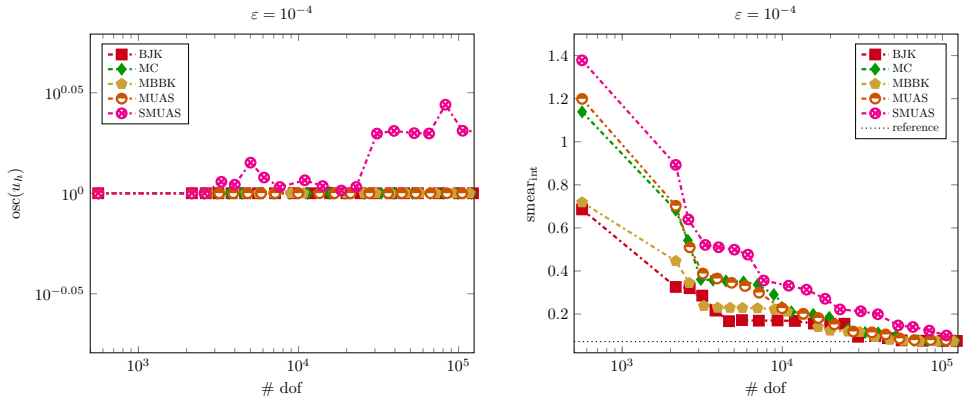
\begin{figure}[tbp]
\centering
\begin{tikzpicture}[scale=0.6]
\begin{loglogaxis}[
    xmin=500, xmax=125000,
    ymin=0.833, ymax=1.2,
    xlabel={$\#\ \mathrm{dof}$},
    ylabel={$\mathrm{osc}(u_h)$},
    title={$\varepsilon=10^{-4}$},
    legend pos=north west,
    legend cell align=left,
    legend style={nodes={scale=0.75, transform shape}}
]
\addplot[color=crimson,  mark=square*, line width = 0.5mm, dashdotted,,mark options = {scale= 1.5, solid}]
coordinates{( 561.0 , 1.0000000000000002 )( 2158.0 , 1.0000000000000002 )( 2635.0 , 1.0000000000100002 )( 3161.0 , 1.0000000000000002 )( 3823.0 , 1.0000000000000162 )( 4641.0 , 1.0000000000000002 )( 5591.0 , 1.0000000000000004 )( 7132.0 , 1.0000000000000002 )( 9332.0 , 1.0000000000000002 )( 11967.0 , 1.0000000000000018 )( 15753.0 , 1.0000000000000002 )( 19947.0 , 1.0000000000000007 )( 24431.0 , 1.0000000000000027 )( 29742.0 , 1.000000000000236 )( 38065.0 , 1.0000000000000322 )( 45682.0 , 1.0000000000000004 )( 55634.0 , 1.0 )( 75368.0 , 1.0000000001492337 )( 100684.0 , 1.0000000000001625 )( 122964.0 , 1.000000000000424 )( 156944.0 , 1.000000000185522 )( 202645.0 , 1.0000000167076557 )( 268857.0 , 1.0000000753941307 )};
\addlegendentry{BJK} 
\addplot[color=dark_green,  mark=diamond*, line width = 0.5mm, dashdotted,,mark options = {scale= 1.5, solid}]
coordinates{( 561.0 , 1.0000000000000002 )( 2158.0 , 1.0000000000000002 )( 2569.0 , 1.0000000000002198 )( 3097.0 , 1.0000000000006113 )( 3689.0 , 1.0000000000411016 )( 4442.0 , 1.000000000289738 )( 5721.0 , 1.0000000000067548 )( 7009.0 , 1.0000000000225153 )( 8839.0 , 1.000000000014486 )( 11392.0 , 1.0000000000000526 )( 15587.0 , 1.0000000003999998 )( 19641.0 , 1.0000000000564793 )( 25745.0 , 1.0000000000083837 )( 32866.0 , 1.0000000000000075 )( 39645.0 , 1.0000000000652791 )( 51477.0 , 1.0000000000088938 )( 66422.0 , 1.000000000548936 )( 80603.0 , 1.0000000021635762 )( 99109.0 , 1.0000000003814127 )( 120871.0 , 1.0000000001259721 )( 169344.0 , 1.0000000009080923 )( 207314.0 , 1.0000000003551754 )( 278667.0 , 1.0000000001055644 )};
\addlegendentry{MC} 
\addplot[color=gold,  mark=pentagon*, line width = 0.5mm, dashdotted,,mark options = {scale= 1.5, solid}]
coordinates{( 561.0 , 1.0000000000000002 )( 2158.0 , 1.00000000000486 )( 2623.0 , 1.000000000000409 )( 3233.0 , 1.000000001203706 )( 3941.0 , 1.0000000000062623 )( 4711.0 , 1.000000000000129 )( 5634.0 , 1.0000000005722878 )( 6994.0 , 1.0000000006953704 )( 8933.0 , 1.0000000012898742 )( 10732.0 , 1.0000000000017288 )( 13426.0 , 1.0000000002163163 )( 16713.0 , 1.0000000000000002 )( 19862.0 , 1.0000000528612487 )( 25665.0 , 1.000000133151704 )( 30418.0 , 1.000000124617466 )( 37855.0 , 1.0000000056223302 )( 46855.0 , 1.0000000489724588 )( 57203.0 , 1.0000000349378637 )( 72323.0 , 1.0000000000472964 )( 92637.0 , 1.0000000007679375 )( 117175.0 , 1.000000129976361 )( 140345.0 , 1.000001204000938 )( 169661.0 , 1.0000004213440827 )( 213336.0 , 1.0000001262839326 )( 265927.0 , 1.0000005291286014 )};
\addlegendentry{MBBK} 
\addplot[color=carrotorange,  mark=halfcircle*, line width = 0.5mm, dashdotted,,mark options = {scale= 1.5, solid}]
coordinates{( 561.0 , 1.0000000000000002 )( 2158.0 , 1.0000000000000002 )( 2637.0 , 1.0000000000000002 )( 3190.0 , 1.0000000000000002 )( 3915.0 , 1.0000000000000002 )( 4821.0 , 1.0000000000000002 )( 5845.0 , 1.0000000000000002 )( 7152.0 , 1.0000000000000002 )( 9959.0 , 1.0000000000000002 )( 13513.0 , 1.0000000000000002 )( 16786.0 , 1.0000000000000002 )( 21148.0 , 1.0000000000000002 )( 27714.0 , 1.0000000000000002 )( 36246.0 , 1.0000000000000002 )( 44177.0 , 1.0000000000000002 )( 53954.0 , 1.0000000000000022 )( 66779.0 , 1.0000000000000093 )( 84437.0 , 1.0 )( 105126.0 , 1.0000000000000007 )( 135151.0 , 1.0000000000000022 )( 171650.0 , 1.0000000000000384 )( 227865.0 , 1.000000000000827 )( 295950.0 , 1.0000000000003166 )};
\addlegendentry{MUAS} 
\addplot[color=magenta,  mark=otimes, line width = 0.5mm, dashdotted,,mark options = {scale= 1.5, solid}]
coordinates{( 561.0 , 1.0 )( 2158.0 , 1.0000000000000002 )( 2587.0 , 1.0000000000000002 )( 3299.0 , 1.0132146999592702 )( 4000.0 , 1.00983621534282 )( 5017.0 , 1.0356875303382 )( 6095.0 , 1.0181541001139 )( 7645.0 , 1.00716562063517 )( 10952.0 , 1.0149800225708001 )( 14159.0 , 1.00826851719085 )( 18529.0 , 1.00315880277046 )( 22929.0 , 1.0072368779766 )( 30905.0 , 1.0707769334257 )( 39491.0 , 1.0740437965 )( 52850.0 , 1.0712750142677 )( 65387.0 , 1.0708203841291999 )( 82785.0 , 1.106730234136 )( 106041.0 , 1.0741893228311 )( 144055.0 , 1.0730456756516 )( 185032.0 , 1.0731148279254 )( 236670.0 , 1.012892101049 )( 297532.0 , 1.0101170452358 )};
\addlegendentry{SMUAS} 
\end{loglogaxis}
\end{tikzpicture}
\hspace*{1em}
\begin{tikzpicture}[scale=0.6]
\begin{semilogxaxis}[
    xmin=500, xmax=125000,
    ymin=0.02, ymax=1.5,
    xlabel={$\#\ \mathrm{dof}$},
    ylabel={$\mathrm{smear}_{\mathrm{int}}$},
    title={$\varepsilon=10^{-4}$},
    legend pos=north east,
    legend cell align=left,
    legend style={nodes={scale=0.75, transform shape}}
]
\addplot[color=crimson,  mark=square*, line width = 0.5mm, dashdotted,,mark options = {scale= 1.5, solid}]
coordinates{( 561.0 , 0.685771 )( 2158.0 , 0.325858 )( 2635.0 , 0.32028 )( 3161.0 , 0.285568 )( 3823.0 , 0.216121 )( 4641.0 , 0.166813 )( 5591.0 , 0.171895 )( 7132.0 , 0.169152 )( 9332.0 , 0.170207 )( 11967.0 , 0.16794 )( 15753.0 , 0.155967 )( 19947.0 , 0.154606 )( 24431.0 , 0.154937 )( 29742.0 , 0.0945776 )( 38065.0 , 0.0986668 )( 45682.0 , 0.0883239 )( 55634.0 , 0.0777593 )( 75368.0 , 0.0762356 )( 100684.0 , 0.0753804 )( 122964.0 , 0.0750436 )( 156944.0 , 0.0747162 )( 202645.0 , 0.075 )( 268857.0 , 0.0747599 )};
\addlegendentry{BJK} 
\addplot[color=dark_green,  mark=diamond*, line width = 0.5mm, dashdotted,,mark options = {scale= 1.5, solid}]
coordinates{( 561.0 , 1.13885 )( 2158.0 , 0.682461 )( 2569.0 , 0.541969 )( 3097.0 , 0.360517 )( 3689.0 , 0.357901 )( 4442.0 , 0.353694 )( 5721.0 , 0.349108 )( 7009.0 , 0.336666 )( 8839.0 , 0.290334 )( 11392.0 , 0.2101 )( 15587.0 , 0.19828 )( 19641.0 , 0.184736 )( 25745.0 , 0.129669 )( 32866.0 , 0.114355 )( 39645.0 , 0.112422 )( 51477.0 , 0.1005 )( 66422.0 , 0.0767191 )( 80603.0 , 0.0764353 )( 99109.0 , 0.0750193 )( 120871.0 , 0.0742226 )( 169344.0 , 0.073858 )( 207314.0 , 0.07329 )( 278667.0 , 0.073015 )};
\addlegendentry{MC} 
\addplot[color=gold,  mark=pentagon*, line width = 0.5mm, dashdotted,,mark options = {scale= 1.5, solid}]
coordinates{( 561.0 , 0.718773 )( 2158.0 , 0.446392 )( 2623.0 , 0.344179 )( 3233.0 , 0.237562 )( 3941.0 , 0.229129 )( 4711.0 , 0.228974 )( 5634.0 , 0.228304 )( 6994.0 , 0.226889 )( 8933.0 , 0.222143 )( 10732.0 , 0.211347 )( 13426.0 , 0.197179 )( 16713.0 , 0.139937 )( 19862.0 , 0.123022 )( 25665.0 , 0.119034 )( 30418.0 , 0.116827 )( 37855.0 , 0.0933793 )( 46855.0 , 0.0802218 )( 57203.0 , 0.0763436 )( 72323.0 , 0.0745634 )( 92637.0 , 0.0734886 )( 117175.0 , 0.0733408 )( 140345.0 , 0.0731938 )( 169661.0 , 0.0729817 )( 213336.0 , 0.0727839 )( 265927.0 , 0.0726573 )};
\addlegendentry{MBBK} 
\addplot[color=carrotorange,  mark=halfcircle*, line width = 0.5mm, dashdotted,,mark options = {scale= 1.5, solid}]
coordinates{( 561.0 , 1.19893 )( 2158.0 , 0.702809 )( 2637.0 , 0.509842 )( 3190.0 , 0.388052 )( 3915.0 , 0.365082 )( 4821.0 , 0.345006 )( 5845.0 , 0.330084 )( 7152.0 , 0.299085 )( 9959.0 , 0.226746 )( 13513.0 , 0.201238 )( 16786.0 , 0.181363 )( 21148.0 , 0.151775 )( 27714.0 , 0.118896 )( 36246.0 , 0.11471 )( 44177.0 , 0.105415 )( 53954.0 , 0.0867222 )( 66779.0 , 0.0790042 )( 84437.0 , 0.0791571 )( 105126.0 , 0.0783651 )( 135151.0 , 0.0772726 )( 171650.0 , 0.0751012 )( 227865.0 , 0.0755705 )( 295950.0 , 0.0755113 )};
\addlegendentry{MUAS} 
\addplot[color=magenta,  mark=otimes, line width = 0.5mm, dashdotted,,mark options = {scale= 1.5, solid}]
coordinates{( 561.0 , 1.37798 )( 2158.0 , 0.892482 )( 2587.0 , 0.639355 )( 3299.0 , 0.520956 )( 4000.0 , 0.510257 )( 5017.0 , 0.498879 )( 6095.0 , 0.475772 )( 7645.0 , 0.35568 )( 10952.0 , 0.331568 )( 14159.0 , 0.313166 )( 18529.0 , 0.270935 )( 22929.0 , 0.221408 )( 30905.0 , 0.212912 )( 39491.0 , 0.198954 )( 52850.0 , 0.14697 )( 65387.0 , 0.139883 )( 82785.0 , 0.12324 )( 106041.0 , 0.100471 )( 144055.0 , 0.0971396 )( 185032.0 , 0.0909445 )( 236670.0 , 0.0802964 )( 297532.0 , 0.0722521 )};
\addlegendentry{SMUAS} 
\addplot[color=black, line width=0.25mm, dotted]
coordinates{(500,0.0723) (125000,0.0723)};
\addlegendentry{reference}
\end{semilogxaxis}
\end{tikzpicture}
\caption{Example~\ref{ex:hemker}: Spurious oscillations for different limiters (left), thickness of the internal layer along the cutline $x=4$ (right).}
\label{fig:hemker_osc}
\end{figure}

Figure~\ref{fig:hemker_osc} (left) presents the oscillation measure.
We observe that all methods, except the SMUAS limiter, satisfy the DMP. With mesh refinement, the grids become partially
non-Delaunay. As noted in \cite{BJK25}, the MC limiter requires a Delaunay
triangulation to guarantee the DMP; however, in the present setting, this condition
is not violated. This can be explained by the adaptive refinement, where only local
portions of the mesh are non-Delaunay, which still allows the DMP to hold in
practice. A closer inspection of the numerical solution shows that the violations
produced by the SMUAS method occur primarily in the vicinity of the circular
obstacle. The SMUAS limiter relies on a reconstructed state $u_{ij}$ based on a
recovered gradient of the numerical solution. In the Hemker problem, the solution
contains very sharp layers around the obstacle and along the convection direction,
and these layers are not aligned with the mesh. Consequently, the reconstructed
gradient may not accurately represent the local behaviour of the solution near the
layer and the limiter introduces insufficient artificial diffusion in some regions,
leading to mild over- and undershoots.

Figure~\ref{fig:hemker_osc} (right) shows the smearing of the internal layer.
A reference value of $0.0723$ is reported in \cite{ACJ11}. The BJK limiter provides
the sharpest resolution of the layer, followed closely by the BBK limiter. For
sufficiently fine meshes, all limiters approach values close to the reference
result, although the SMUAS method exhibits the largest amount of smearing.

\begin{figure}[tbp]
\centering
\includegraphics[width=0.32\linewidth]{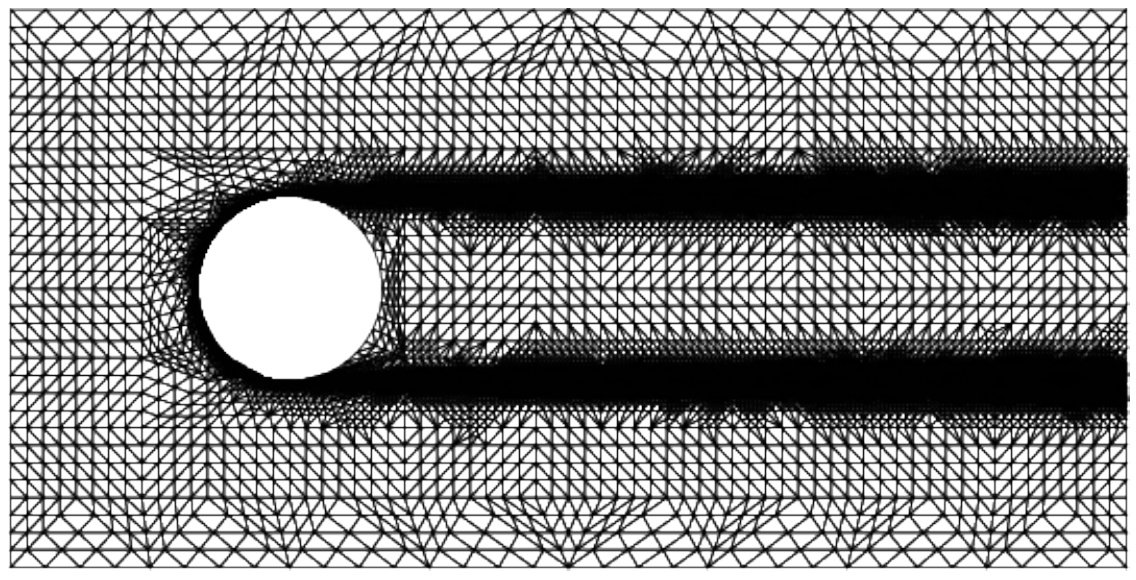}\hspace*{1em}
\includegraphics[width=0.32\linewidth]{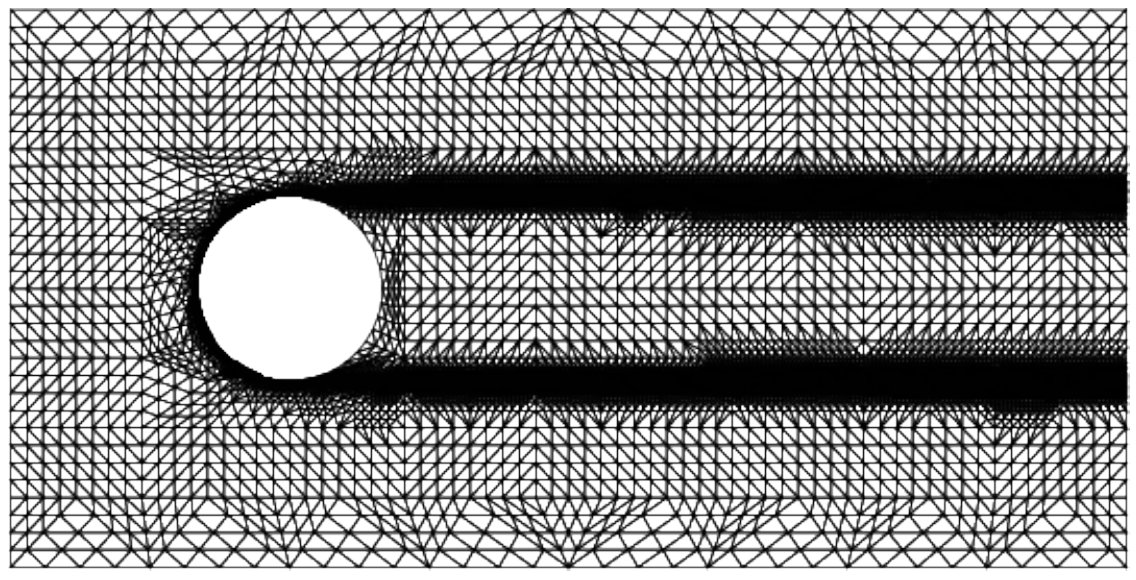}
\caption{Example~\ref{ex:hemker}: Adaptive grids for the AFC methods with $\#\mathrm{dofs}\approx 10^5$.
BJK limiter (left) and MC limiter (right).}
\label{fig:adaptive_afc_hemker}
\end{figure}

\begin{figure}[tbp]
\centering
\includegraphics[width=0.32\linewidth]{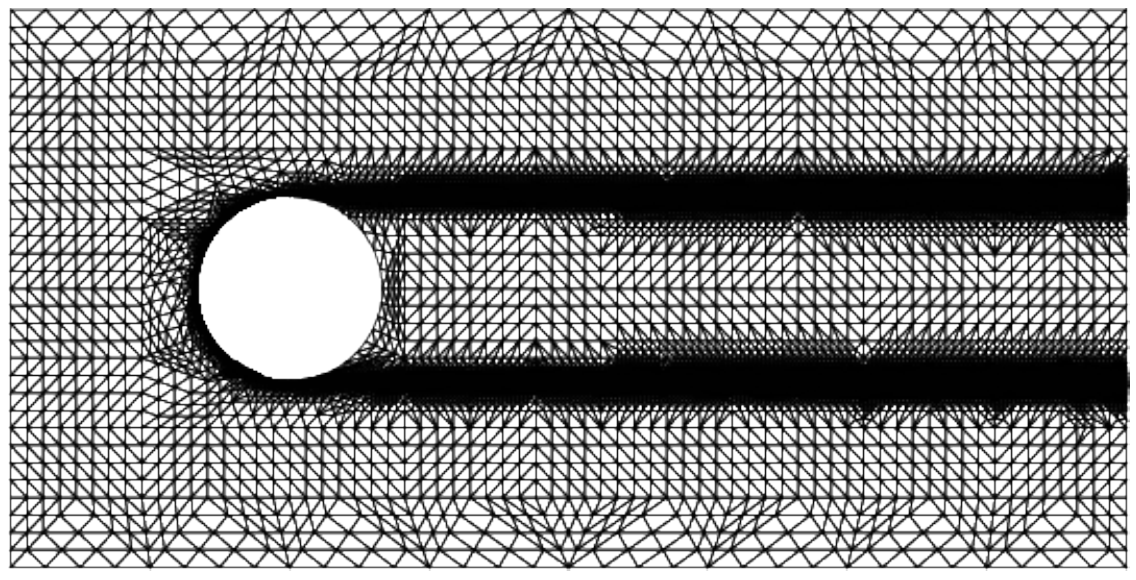}\hfill
\includegraphics[width=0.32\linewidth]{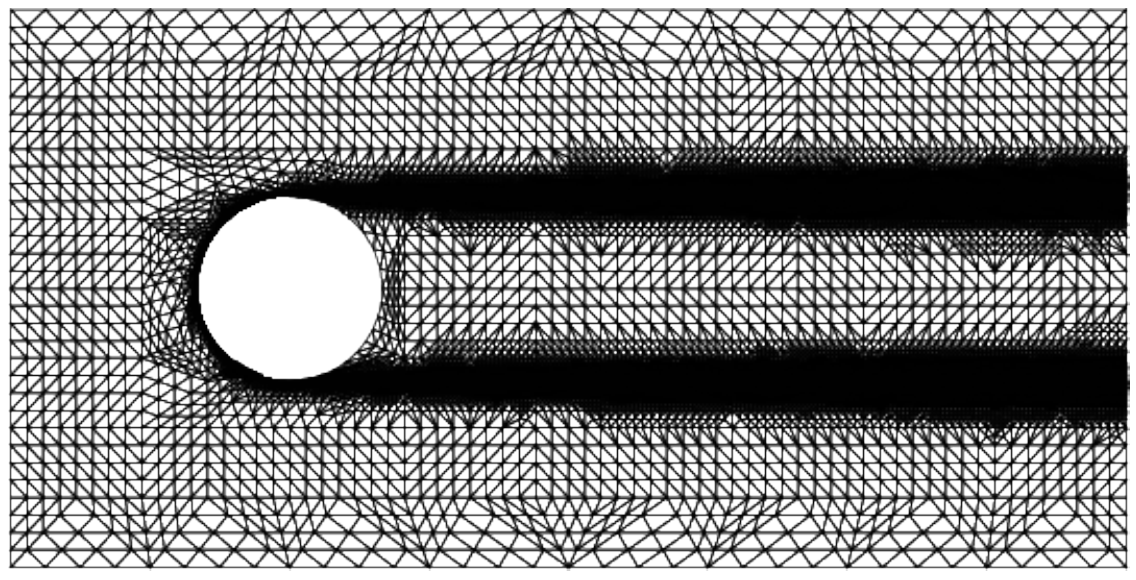}\hfill
\includegraphics[width=0.32\linewidth]{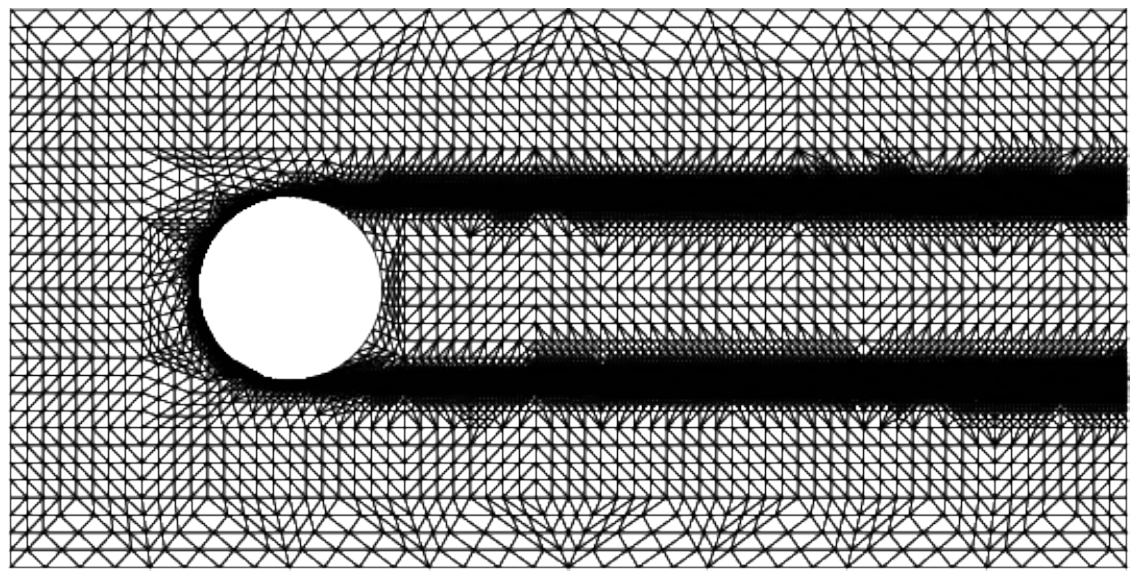}
\caption{Example~\ref{ex:hemker}: Adaptive grids for the algebraically stabilised methods with
$\#\mathrm{dofs}\approx 10^5$. MUAS (left), SMUAS (middle), BBK (right).}
\label{fig:adaptive_as_hemker}
\end{figure}

Finally, Figures~\ref{fig:adaptive_afc_hemker} and~\ref{fig:adaptive_as_hemker} present the adaptive grids obtained
for all methods using approximately $\#\mathrm{dofs}\approx 10^5$.
For the AFC methods, the MC limiter yields the sharpest internal layer,
whereas for the algebraically stabilised schemes, the MUAS and BBK methods provide the best resolution.

\subsection{Solution with Inner Layer, Non-Linear Convection and Linear Reaction}\label{ex:nonlinear_cdr}
\begin{figure}[tbp]
    \centering
    \includegraphics[width=0.4\linewidth]{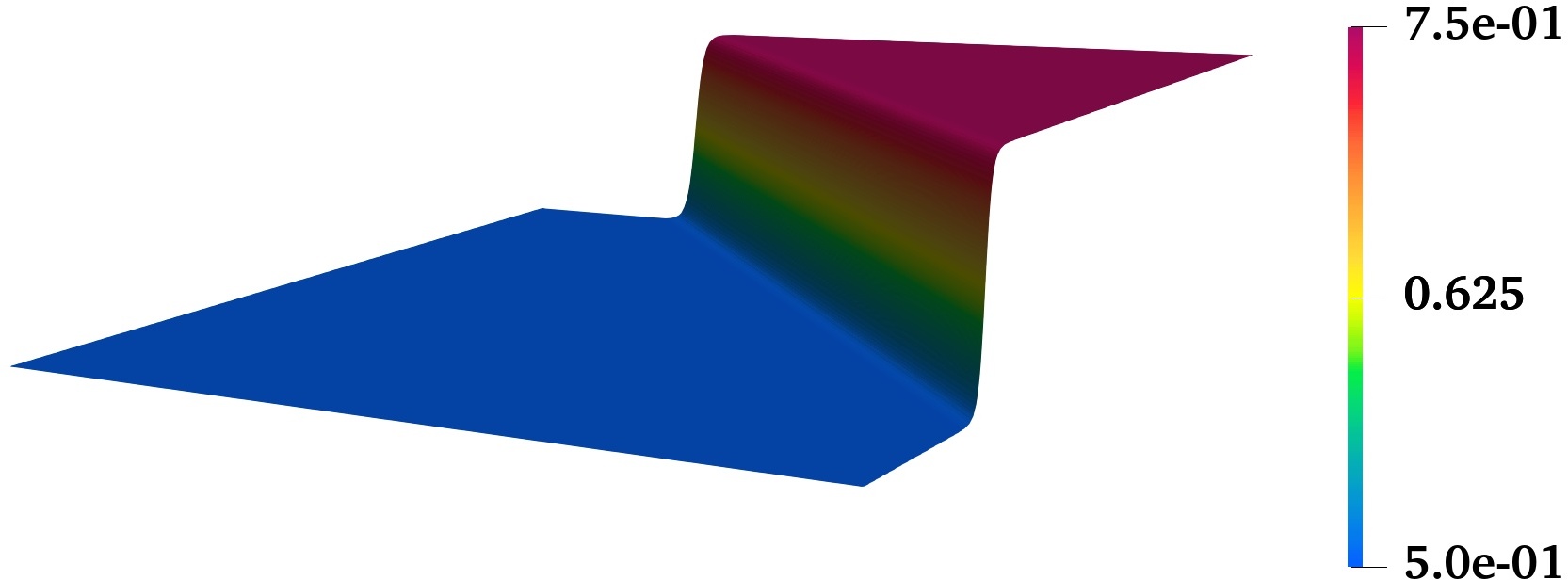}
    \caption{Example~\ref{ex:nonlinear_cdr}: Numerical solution obtained with the AFC scheme and BJK limiter on Grid~1 with $16{,}641$ degrees of freedom.}
    \label{fig:nonlinear_solution}
\end{figure}

Consider Eq.~\eqref{eq:cdr_eqn} with
$$
\varepsilon = 10^{-3}, \quad 
\bb = \left(u,\, u\right)^{\top}, \quad 
c = 0, \quad 
\Omega = (0,1)^2, \quad 
\Gamma_{\mD} = \Gamma.
$$
The right-hand side and boundary conditions are chosen such that
$$
u(x,y) = \frac{3}{4} - \frac{1}{4\left[1 + \exp\left(\frac{-4x + 4y - 1}{32\varepsilon}\right)\right]}
$$
is the exact solution (see Fig.~\ref{fig:nonlinear_solution}).  
This is a non-linear example, with the non-linearity arising from the convection term.  
Equation~\eqref{eq:cdr_cond} is still satisfied in this case.  
An interior layer develops along the line $-4x + 4y - 1 = 0$.  
Among the three grids, Grid~1 is aligned with this interior layer.  
This example serves as the steady-state counterpart of Example~2 from~\cite{ZLG25}.

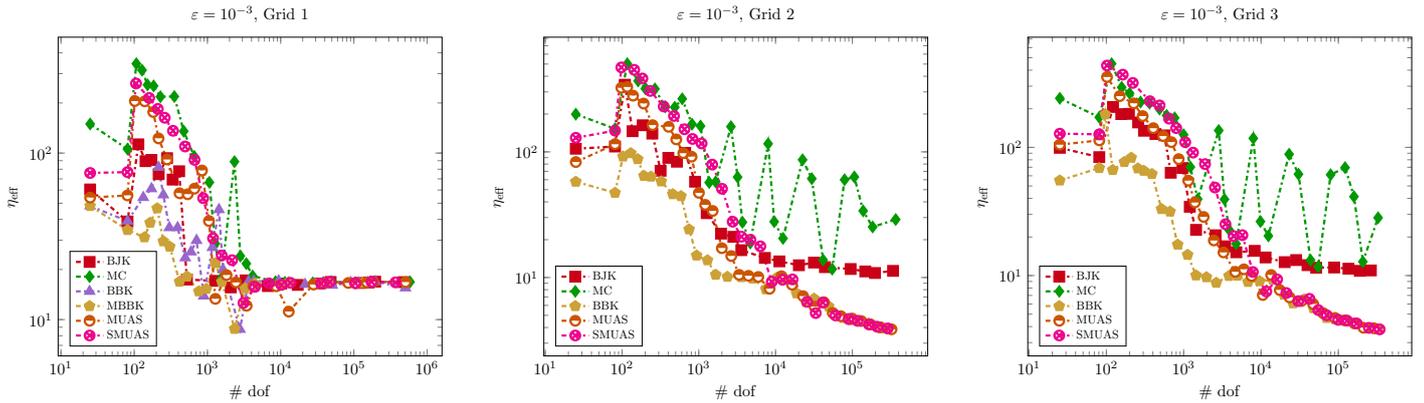
\begin{figure}[tbp]
\centerline{
\begin{tikzpicture}[scale=0.6]
\begin{loglogaxis}[
       legend pos=south west, xlabel = $\#\ \mathrm{dof}$, ylabel=$\eta_{\mathrm{eff}}$,
    legend cell align ={left}, title = {$\varepsilon=10^{-3}$, Grid~1},
    legend style={nodes={scale=0.75, transform shape}}]
]
\addplot[color=crimson,  mark=square*, line width = 0.5mm, dashdotted,,mark options = {scale= 1.5, solid}]
coordinates{( 25.0 , 60.4962 )( 81.0 , 38.8387 )( 114.0 , 112.899 )( 143.0 , 89.4261 )( 172.0 , 90.9943 )( 216.0 , 74.7065 )( 281.0 , 93.2348 )( 335.0 , 69.2778 )( 411.0 , 77.8123 )( 524.0 , 17.3985 )( 1280.0 , 17.0874 )( 2083.0 , 15.6075 )( 3336.0 , 17.2049 )( 6698.0 , 15.8992 )( 10179.0 , 16.463 )( 17096.0 , 16.1574 )( 34239.0 , 16.4757 )( 48118.0 , 16.8285 )( 123241.0 , 16.599 )( 157875.0 , 16.6945 )( 195258.0 , 16.8956 )( 477315.0 , 16.791 )};
\addlegendentry{BJK} 
\addplot[color=dark_green,  mark=diamond*, line width = 0.5mm, dashdotted,,mark options = {scale= 1.5, solid}]
coordinates{( 25.0 , 149.483 )( 81.0 , 105.71 )( 107.0 , 344.022 )( 128.0 , 314.809 )( 152.0 , 256.896 )( 184.0 , 253.106 )( 228.0 , 217.996 )( 353.0 , 218.671 )( 475.0 , 135.307 )( 660.0 , 95.5535 )( 1063.0 , 66.6816 )( 1277.0 , 28.6069 )( 1541.0 , 17.651 )( 2340.0 , 88.8403 )( 2815.0 , 24.0815 )( 3394.0 , 21.6775 )( 4140.0 , 18.2786 )( 8509.0 , 16.8854 )( 11410.0 , 16.7337 )( 27407.0 , 16.606 )( 38369.0 , 16.6497 )( 49157.0 , 16.888 )( 102811.0 , 16.5995 )( 155185.0 , 16.7051 )( 191034.0 , 16.8875 )( 584134.0 , 16.8352 )};
\addlegendentry{MC} 
\addplot[color=amethyst,  mark=triangle*, line width = 0.5mm, dashdotted,,mark options = {scale= 1.5, solid}]
coordinates{( 25.0 , 48.1422 )( 81.0 , 39.0392 )( 133.0 , 54.1285 )( 172.0 , 61.0896 )( 215.0 , 82.1809 )( 250.0 , 56.1262 )( 302.0 , 35.6644 )( 391.0 , 35.6762 )( 498.0 , 23.491 )( 596.0 , 25.4965 )( 710.0 , 29.8831 )( 891.0 , 13.8598 )( 1170.0 , 27.1064 )( 1436.0 , 45.7514 )( 1631.0 , 19.9879 )( 2785.0 , 8.73895 )( 3631.0 , 16.9546 )( 4295.0 , 16.9836 )( 9381.0 , 17.0952 )( 12823.0 , 16.6449 )( 21575.0 , 16.2793 )( 29622.0 , 16.3605 )( 37646.0 , 16.5478 )( 50775.0 , 16.0705 )( 90671.0 , 16.6775 )( 117003.0 , 16.5978 )( 144823.0 , 16.6516 )( 197654.0 , 16.8968 )( 506061.0 , 15.5009 )};
\addlegendentry{BBK} 
\addplot[color=gold,  mark=pentagon*, line width = 0.5mm, dashdotted,,mark options = {scale= 1.5, solid}]
coordinates{( 25.0 , 48.1351 )( 81.0 , 34.5281 )( 138.0 , 31.2099 )( 165.0 , 38.1805 )( 206.0 , 46.5777 )( 250.0 , 29.532 )( 306.0 , 27.4239 )( 416.0 , 16.9034 )( 526.0 , 18.1361 )( 763.0 , 14.7824 )( 966.0 , 15.2932 )( 1280.0 , 21.8037 )( 1573.0 , 16.8307 )( 2349.0 , 8.79583 )( 3263.0 , 15.4001 )( 4433.0 , 16.401 )( 9604.0 , 15.9875 )( 12715.0 , 16.697 )( 29736.0 , 16.3039 )( 46034.0 , 16.7556 )( 133437.0 , 16.5986 )( 180029.0 , 16.853 )( 505732.0 , 16.7825 )};
\addlegendentry{MBBK} 
\addplot[color=carrotorange,  mark=halfcircle*, line width = 0.5mm, dashdotted,,mark options = {scale= 1.5, solid}]
coordinates{( 25.0 , 54.017 )( 81.0 , 55.6594 )( 102.0 , 205.27 )( 139.0 , 204.497 )( 182.0 , 176.469 )( 215.0 , 122.476 )( 281.0 , 91.5753 )( 416.0 , 57.3806 )( 537.0 , 56.5757 )( 665.0 , 61.1037 )( 843.0 , 78.609 )( 1049.0 , 39.091 )( 1283.0 , 13.3177 )( 1822.0 , 18.4974 )( 2445.0 , 16.7187 )( 3429.0 , 12.1126 )( 4236.0 , 15.8904 )( 8685.0 , 15.8522 )( 12894.0 , 11.1657 )( 28075.0 , 16.2397 )( 38921.0 , 16.5332 )( 83973.0 , 16.5969 )( 107152.0 , 16.5829 )( 150964.0 , 16.7226 )( 181723.0 , 16.8298 )( 506284.0 , 16.7813 )};
\addlegendentry{MUAS} 
\addplot[color=magenta,  mark=otimes, line width = 0.5mm, dashdotted,,mark options = {scale= 1.5, solid}]
coordinates{( 25.0 , 75.9669 )( 81.0 , 76.8929 )( 106.0 , 261.786 )( 160.0 , 214.353 )( 210.0 , 184.165 )( 261.0 , 162.958 )( 339.0 , 136.102 )( 495.0 , 109.556 )( 667.0 , 91.0847 )( 875.0 , 53.6081 )( 1193.0 , 30.7948 )( 1567.0 , 24.2929 )( 2169.0 , 22.7434 )( 3102.0 , 12.5918 )( 4463.0 , 15.721 )( 6910.0 , 16.5021 )( 10230.0 , 16.1581 )( 13583.0 , 16.624 )( 35613.0 , 16.4939 )( 46890.0 , 16.8276 )( 103883.0 , 16.5696 )( 174329.0 , 16.8469 )( 369795.0 , 16.7288 )};
\addlegendentry{SMUAS} 
\end{loglogaxis}
\end{tikzpicture}
\hspace*{1em}
\begin{tikzpicture}[scale=0.6]
\begin{loglogaxis}[
       legend pos=south west, xlabel = $\#\ \mathrm{dof}$, ylabel=$\eta_{\mathrm{eff}}$,
    legend cell align ={left}, title = {$\varepsilon=10^{-3}$, Grid~2},
    legend style={nodes={scale=0.75, transform shape}}]
]
\addplot[color=crimson,  mark=square*, line width = 0.5mm, dashdotted,,mark options = {scale= 1.5, solid}]
coordinates{( 25.0 , 105.894 )( 81.0 , 110.276 )( 109.0 , 342.527 )( 137.0 , 146.05 )( 186.0 , 162.872 )( 248.0 , 139.627 )( 318.0 , 70.7055 )( 401.0 , 89.542 )( 519.0 , 83.4591 )( 664.0 , 98.2549 )( 895.0 , 57.8683 )( 1268.0 , 32.5135 )( 1950.0 , 22.3451 )( 2826.0 , 21.0776 )( 3584.0 , 16.504 )( 7264.0 , 14.322 )( 11228.0 , 13.4395 )( 20157.0 , 12.5196 )( 31661.0 , 13.1491 )( 44007.0 , 12.092 )( 95047.0 , 11.6868 )( 147461.0 , 11.1763 )( 197431.0 , 10.8591 )( 339594.0 , 11.2913 )};
\addlegendentry{BJK} 
\addplot[color=dark_green,  mark=diamond*, line width = 0.5mm, dashdotted,,mark options = {scale= 1.5, solid}]
coordinates{( 25.0 , 199.483 )( 81.0 , 152.777 )( 117.0 , 502.902 )( 162.0 , 367.001 )( 201.0 , 318.479 )( 268.0 , 317.111 )( 366.0 , 228.016 )( 486.0 , 226.477 )( 608.0 , 264.375 )( 816.0 , 166.068 )( 1053.0 , 160.514 )( 1366.0 , 56.9993 )( 1659.0 , 57.7266 )( 2624.0 , 158.739 )( 3206.0 , 62.9469 )( 3706.0 , 27.5386 )( 4685.0 , 18.8196 )( 8011.0 , 115.944 )( 9605.0 , 27.8266 )( 12517.0 , 20.472 )( 22489.0 , 86.4796 )( 29651.0 , 61.1917 )( 41626.0 , 13.7937 )( 54688.0 , 11.6711 )( 79302.0 , 60.0181 )( 106273.0 , 63.1512 )( 138433.0 , 33.9206 )( 184010.0 , 25.2483 )( 367913.0 , 28.9447 )};
\addlegendentry{MC} 
%\addplot[color=amethyst,  mark=triangle*, line width = 0.5mm, dashdotted,,mark options = {scale= 1.5, solid}]
%coordinates{( 25.0 , 58.2986 )( 81.0 , 50.042 )( 103.0 , 127.477 )( 123.0 , 120.868 )( 145.0 , 156.827 )( 183.0 , 147.25 )( 269.0 , 111.26 )( 354.0 , 98.7477 )( 479.0 , 74.7907 )( 577.0 , 70.5214 )( 722.0 , 55.3231 )( 901.0 , 54.5963 )( 1119.0 , 38.6345 )( 1336.0 , 23.8761 )( 1697.0 , 12.0682 )( 2171.0 , 11.1338 )( 2850.0 , 10.0014 )( 3950.0 , 10.5096 )( 5269.0 , 9.93373 )( 8010.0 , 8.58995 )( 10466.0 , 9.9227 )( 14204.0 , 9.74641 )( 19395.0 , 7.35027 )( 25490.0 , 6.39967 )( 35795.0 , 6.78683 )( 46381.0 , 6.09846 )( 63784.0 , 5.22956 )( 83793.0 , 4.66134 )( 116537.0 , 4.46011 )( 146638.0 , 4.35452 )( 194479.0 , 4.15983 )( 288709.0 , 3.90459 )};
%\addlegendentry{BBK} 
\addplot[color=gold,  mark=pentagon*, line width = 0.5mm, dashdotted,,mark options = {scale= 1.5, solid}]
coordinates{( 25.0 , 57.8454 )( 81.0 , 47.2407 )( 103.0 , 91.8043 )( 130.0 , 97.6828 )( 160.0 , 87.562 )( 193.0 , 64.4271 )( 238.0 , 63.5645 )( 324.0 , 58.2298 )( 458.0 , 46.0496 )( 586.0 , 44.4751 )( 748.0 , 24.1603 )( 936.0 , 14.9661 )( 1297.0 , 13.715 )( 1661.0 , 10.5154 )( 2358.0 , 10.1297 )( 3585.0 , 10.0736 )( 5029.0 , 9.81597 )( 7488.0 , 8.0528 )( 10403.0 , 9.97052 )( 13348.0 , 9.73936 )( 18766.0 , 7.48167 )( 25017.0 , 6.38016 )( 31753.0 , 6.79031 )( 49117.0 , 5.90209 )( 76136.0 , 4.77197 )( 97199.0 , 4.57381 )( 121299.0 , 4.52459 )( 162292.0 , 4.30256 )( 237909.0 , 3.93419 )( 300613.0 , 3.90307 )};
\addlegendentry{BBK} 
\addplot[color=carrotorange,  mark=halfcircle*, line width = 0.5mm, dashdotted,,mark options = {scale= 1.5, solid}]
coordinates{( 25.0 , 82.9689 )( 81.0 , 115.106 )( 98.0 , 323.508 )( 117.0 , 328.472 )( 139.0 , 280.083 )( 191.0 , 242.332 )( 250.0 , 163.137 )( 405.0 , 158.163 )( 509.0 , 125.668 )( 628.0 , 97.6489 )( 805.0 , 91.5346 )( 1004.0 , 47.1646 )( 1212.0 , 37.8765 )( 1517.0 , 33.8699 )( 1996.0 , 17.113 )( 2651.0 , 14.8162 )( 3354.0 , 10.4981 )( 4437.0 , 10.3096 )( 5960.0 , 10.0997 )( 8374.0 , 8.13565 )( 11967.0 , 10.2402 )( 16410.0 , 8.75961 )( 22873.0 , 7.10086 )( 32208.0 , 5.86585 )( 42671.0 , 6.35305 )( 55614.0 , 5.22608 )( 71021.0 , 4.92212 )( 99115.0 , 4.7129 )( 130994.0 , 4.50347 )( 172908.0 , 4.23461 )( 232524.0 , 4.01848 )( 324309.0 , 3.86752 )};
\addlegendentry{MUAS} 
\addplot[color=magenta,  mark=otimes, line width = 0.5mm, dashdotted,,mark options = {scale= 1.5, solid}]
coordinates{( 25.0 , 129.392 )( 81.0 , 147.626 )( 98.0 , 470.115 )( 144.0 , 449.775 )( 184.0 , 383.862 )( 233.0 , 305.428 )( 349.0 , 229.766 )( 475.0 , 193.067 )( 651.0 , 151.362 )( 828.0 , 126.683 )( 1088.0 , 116.705 )( 1494.0 , 79.2471 )( 2007.0 , 50.8482 )( 2764.0 , 27.8018 )( 3600.0 , 21.1246 )( 4715.0 , 19.9604 )( 6337.0 , 17.7142 )( 8935.0 , 9.26097 )( 12345.0 , 9.59328 )( 16648.0 , 9.66285 )( 25755.0 , 6.40322 )( 33395.0 , 5.21123 )( 41261.0 , 6.32131 )( 60405.0 , 4.99048 )( 89383.0 , 4.67334 )( 117360.0 , 4.51414 )( 164595.0 , 4.22844 )( 210479.0 , 4.05745 )( 288803.0 , 3.9353 )};
\addlegendentry{SMUAS} 
\end{loglogaxis}
\end{tikzpicture}
\hspace*{1em}
\begin{tikzpicture}[scale=0.6]
\begin{loglogaxis}[
       legend pos=south west, xlabel = $\#\ \mathrm{dof}$, ylabel=$\eta_{\mathrm{eff}}$,
    legend cell align ={left}, title = {$\varepsilon=10^{-3}$, Grid~3},
    legend style={nodes={scale=0.75, transform shape}}]
]
\addplot[color=crimson,  mark=square*, line width = 0.5mm, dashdotted,,mark options = {scale= 1.5, solid}]
coordinates{( 25.0 , 99.1236 )( 81.0 , 83.9134 )( 122.0 , 206.492 )( 156.0 , 181.665 )( 196.0 , 182.396 )( 258.0 , 156.468 )( 307.0 , 134.674 )( 432.0 , 127.054 )( 564.0 , 124.04 )( 681.0 , 63.3684 )( 965.0 , 69.0206 )( 1201.0 , 34.3943 )( 1447.0 , 22.7259 )( 2586.0 , 20.6329 )( 3436.0 , 16.9175 )( 4795.0 , 15.1908 )( 8274.0 , 15.5929 )( 11449.0 , 13.8442 )( 20244.0 , 12.7293 )( 28484.0 , 13.2128 )( 40008.0 , 12.0585 )( 51737.0 , 11.5086 )( 88425.0 , 11.562 )( 135342.0 , 11.3147 )( 189597.0 , 10.8848 )( 263573.0 , 10.9653 )};
\addlegendentry{BJK} 
\addplot[color=dark_green,  mark=diamond*, line width = 0.5mm, dashdotted,,mark options = {scale= 1.5, solid}]
coordinates{( 25.0 , 241.308 )( 81.0 , 171.64 )( 117.0 , 448.914 )( 158.0 , 293.657 )( 200.0 , 263.988 )( 277.0 , 225.652 )( 362.0 , 221.309 )( 487.0 , 199.135 )( 603.0 , 177.327 )( 765.0 , 169.84 )( 999.0 , 125.387 )( 1220.0 , 70.1865 )( 1518.0 , 39.4925 )( 2838.0 , 135.345 )( 3371.0 , 39.3292 )( 3877.0 , 22.6887 )( 4754.0 , 17.7021 )( 7951.0 , 117.619 )( 9828.0 , 26.3904 )( 12456.0 , 20.3881 )( 23093.0 , 88.3724 )( 30548.0 , 61.7525 )( 43723.0 , 13.1854 )( 54867.0 , 11.7167 )( 79573.0 , 61.1891 )( 122392.0 , 69.4899 )( 159471.0 , 41.2681 )( 207854.0 , 12.8168 )( 326961.0 , 28.2465 )};
\addlegendentry{MC} 
%\addplot[color=amethyst,  mark=triangle*, line width = 0.5mm, dashdotted,,mark options = {scale= 1.5, solid}]
%coordinates{( 25.0 , 55.2015 )( 81.0 , 86.5959 )( 106.0 , 267.887 )( 132.0 , 180.569 )( 164.0 , 143.081 )( 195.0 , 122.845 )( 237.0 , 105.592 )( 327.0 , 100.194 )( 423.0 , 113.279 )( 526.0 , 81.3179 )( 678.0 , 77.9043 )( 846.0 , 49.1845 )( 1081.0 , 27.2331 )( 1318.0 , 14.202 )( 1762.0 , 12.9438 )( 2402.0 , 17.3747 )( 3280.0 , 9.23095 )( 4641.0 , 11.211 )( 6150.0 , 10.4854 )( 8791.0 , 7.91832 )( 12095.0 , 10.2702 )( 15435.0 , 8.56018 )( 25481.0 , 6.35227 )( 36086.0 , 6.25986 )( 51763.0 , 5.32959 )( 71977.0 , 4.63772 )( 100287.0 , 4.51827 )( 137797.0 , 4.3371 )( 183672.0 , 4.17094 )( 238317.0 , 4.01748 )( 324958.0 , 3.85662 )};
%\addlegendentry{BBK} 
\addplot[color=gold,  mark=pentagon*, line width = 0.5mm, dashdotted,,mark options = {scale= 1.5, solid}]
coordinates{( 25.0 , 55.1705 )( 81.0 , 69.013 )( 97.0 , 180.836 )( 121.0 , 66.627 )( 171.0 , 76.9161 )( 210.0 , 83.0373 )( 256.0 , 68.5502 )( 305.0 , 65.664 )( 389.0 , 62.0234 )( 517.0 , 33.0452 )( 667.0 , 31.5619 )( 837.0 , 17.4401 )( 1125.0 , 14.5167 )( 1445.0 , 10.0269 )( 1970.0 , 9.70817 )( 2637.0 , 8.78856 )( 3851.0 , 9.9343 )( 4943.0 , 9.97784 )( 6794.0 , 8.92344 )( 9617.0 , 9.04406 )( 14271.0 , 9.28831 )( 19813.0 , 7.08117 )( 25975.0 , 6.0748 )( 36188.0 , 6.52743 )( 47667.0 , 5.56554 )( 70391.0 , 4.67047 )( 91582.0 , 4.56477 )( 128514.0 , 4.44619 )( 167238.0 , 4.20417 )( 236889.0 , 3.90066 )( 325222.0 , 3.81809 )};
\addlegendentry{BBK} 
\addplot[color=carrotorange,  mark=halfcircle*, line width = 0.5mm, dashdotted,,mark options = {scale= 1.5, solid}]
coordinates{( 25.0 , 104.928 )( 81.0 , 113.384 )( 102.0 , 353.561 )( 149.0 , 250.922 )( 226.0 , 220.171 )( 295.0 , 175.282 )( 405.0 , 141.096 )( 542.0 , 122.537 )( 693.0 , 109.795 )( 863.0 , 80.9278 )( 1147.0 , 54.9655 )( 1397.0 , 37.176 )( 1811.0 , 28.6575 )( 2456.0 , 18.7829 )( 3270.0 , 15.152 )( 4653.0 , 10.6898 )( 6004.0 , 11.1851 )( 7549.0 , 9.91739 )( 10548.0 , 7.05216 )( 13539.0 , 10.0953 )( 17488.0 , 7.70369 )( 23394.0 , 6.75701 )( 33425.0 , 6.33283 )( 47393.0 , 5.99468 )( 63970.0 , 5.08493 )( 90458.0 , 4.67081 )( 120666.0 , 4.47575 )( 161591.0 , 4.23624 )( 211845.0 , 3.91718 )( 288261.0 , 3.88125 )};
\addlegendentry{MUAS} 
\addplot[color=magenta,  mark=otimes, line width = 0.5mm, dashdotted,,mark options = {scale= 1.5, solid}]
coordinates{( 25.0 , 127.801 )( 81.0 , 126.379 )( 102.0 , 434.516 )( 163.0 , 369.015 )( 220.0 , 318.053 )( 366.0 , 229.023 )( 488.0 , 212.393 )( 647.0 , 167.565 )( 801.0 , 141.294 )( 1048.0 , 109.779 )( 1322.0 , 90.9723 )( 1889.0 , 74.0366 )( 2539.0 , 48.6255 )( 3472.0 , 25.1622 )( 4432.0 , 20.322 )( 5743.0 , 20.659 )( 7839.0 , 10.6624 )( 11717.0 , 7.51744 )( 16150.0 , 9.33807 )( 21255.0 , 7.30965 )( 30321.0 , 6.31927 )( 42331.0 , 6.58593 )( 55162.0 , 5.37844 )( 70952.0 , 4.89903 )( 99795.0 , 4.51446 )( 129697.0 , 4.45348 )( 173068.0 , 4.2343 )( 249514.0 , 3.90839 )( 342796.0 , 3.80729 )};
\addlegendentry{SMUAS} 
\end{loglogaxis}
\end{tikzpicture}}
\caption{Example~\ref{ex:nonlinear_cdr}: Effectivity index for Grids~1–3 (left to right). }
\label{fig:effec_index_nonlinear_cdr}
\end{figure}

Fig.~\ref{fig:effec_index_nonlinear_cdr} presents the effectivity index.
For Grid~1, the values obtained for all methods on sufficiently fine meshes are
close and approximately equal to $16$. The decay is smooth for most limiters,
whereas the MC limiter exhibits several noticeable jumps. On Grids~2 and~3, this
behaviour becomes more pronounced. For these grids, the effectivity index for all
methods except the BJK limiter approaches values around $4$, while the BJK limiter remains closer to $10$.

This behaviour can be attributed to the nonlinear nature of the problem. In this
example, the convection field depends on the solution itself, $\bb=(u,u)^{\top}$, and hence
the effective transport direction changes during the nonlinear iteration. On
Grid~1, the interior layer is approximately aligned with the mesh, and therefore
all methods achieve similar estimator behaviour. On Grids~2 and~3, however, the
layer is not aligned with the mesh, which makes the stabilization more sensitive
to the limiter. The BJK limiter enforces stronger upwinding and therefore reacts
more robustly to the moving transport direction, whereas the other limiters reduce
artificial diffusion and consequently produce smaller estimator values. The jumps
observed for the MC limiter occur when the limiter switches between active and
inactive states during the adaptive refinement, which causes abrupt changes in the stabilization contribution to the residual.

\begin{figure}[tbp]
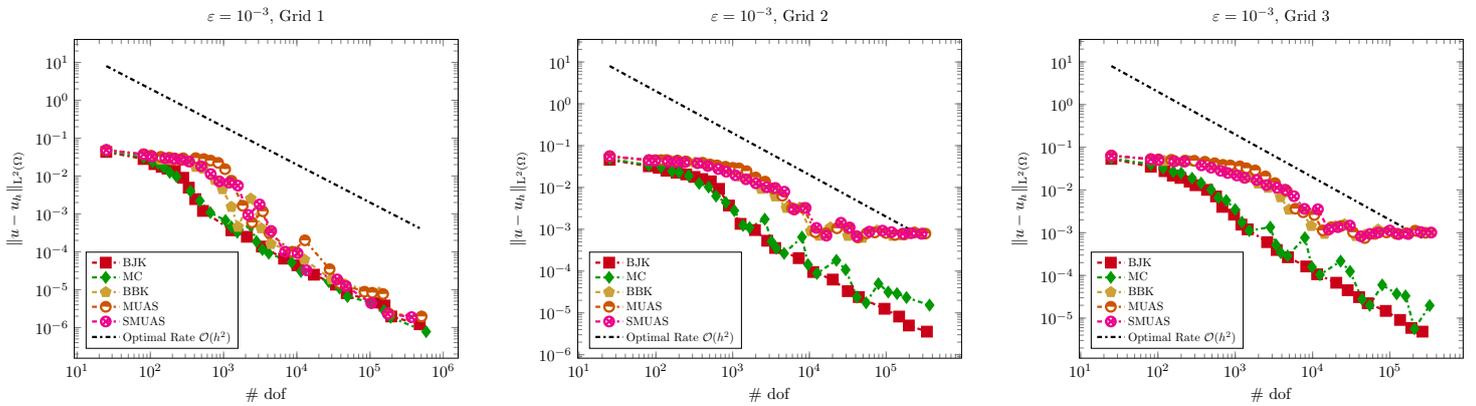

\centerline{
}
\caption{Example~\ref{ex:nonlinear_cdr}: Error in the $\mL^2(\Omega)$ norm for Grids~1–3 (left to right).}
\label{fig:l2_error_nonlinear_cdr}
\end{figure}

\begin{figure}[tbp]
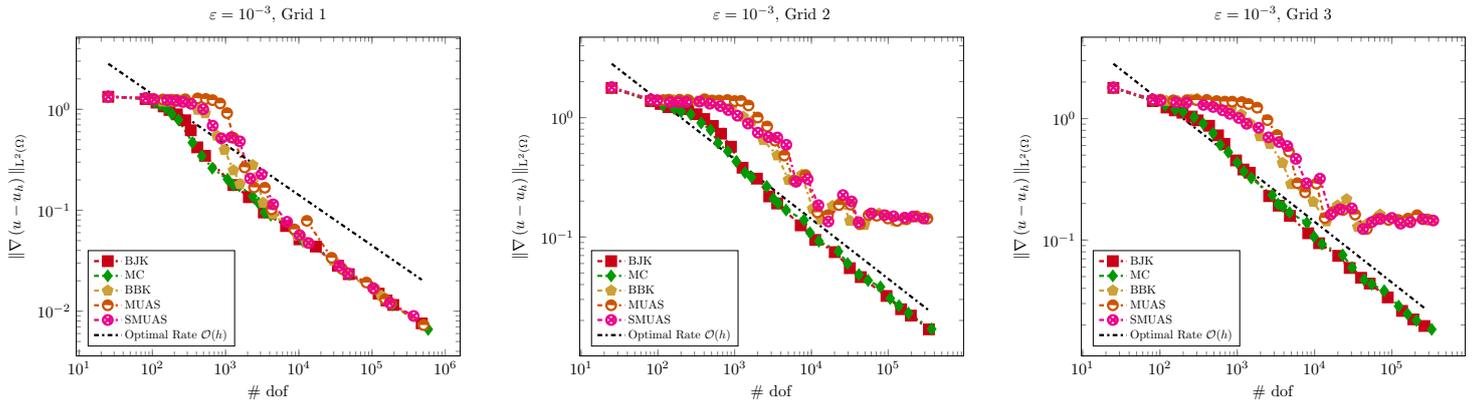

\centerline{
%
}
\caption{Example~\ref{ex:nonlinear_cdr}: Error of the gradient in the $\mL^2(\Omega)$ norm for Grids~1–3 (left to right).}
\label{fig:h1_error_nonlinear_cdr}
\end{figure}

Figs.~\ref{fig:l2_error_nonlinear_cdr} and~\ref{fig:h1_error_nonlinear_cdr} show
the $\mL^2$ error of the solution and the $\mL^2$ error of its gradient,
respectively. For Grid~1, the error decays optimally for all methods. In contrast,
on Grids~2 and~3, all methods except the BJK and MC limiters lose their optimal
convergence rate. This again reflects the influence of mesh alignment. Since the
location and orientation of the interior layer depend on the solution, the problem
resembles a transport-dominated equation with a moving layer. The stronger
upwinding behaviour of the BJK limiter allows it to capture this layer more
accurately on non-aligned meshes and therefore yields the most accurate results in this example.

We also note that, in contrast to the Hemker problem, all limiters satisfy the
discrete maximum principle for this test case. The numerical solution remains
bounded within the physically relevant range $[0.5,\,0.75]$ for all meshes and all
methods.

\begin{figure}[tbp]
\centerline{
\begin{tikzpicture}[scale=0.6]
\begin{semilogxaxis}[
    legend pos=north east, xlabel = $\#\ \mathrm{dof}$, ylabel = iterations+rejections,
    legend cell align ={left}, title = {$\varepsilon=10^{-3}$, Grid~1},
     legend style={nodes={scale=0.75, transform shape}}]
\addplot[color=crimson,  mark=square*, line width = 0.5mm, dashdotted,,mark options = {scale= 1.5, solid}]
coordinates{( 25.0 , 139.0 )( 81.0 , 125.0 )( 114.0 , 94.0 )( 143.0 , 83.0 )( 172.0 , 73.0 )( 216.0 , 73.0 )( 281.0 , 53.0 )( 335.0 , 54.0 )( 411.0 , 42.0 )( 524.0 , 42.0 )( 1280.0 , 22.0 )( 2083.0 , 15.0 )( 3336.0 , 11.0 )( 6698.0 , 7.0 )( 10179.0 , 8.0 )( 17096.0 , 2.0 )( 34239.0 , 2.0 )( 48118.0 , 2.0 )( 123241.0 , 2.0 )( 157875.0 , 2.0 )( 195258.0 , 2.0 )( 477315.0 , 2.0 )};
\addlegendentry{BJK} 
\addplot[color=dark_green,  mark=diamond*, line width = 0.5mm, dashdotted,,mark options = {scale= 1.5, solid}]
coordinates{( 25.0 , 43.0 )( 81.0 , 42.0 )( 107.0 , 43.0 )( 128.0 , 48.0 )( 152.0 , 42.0 )( 184.0 , 34.0 )( 228.0 , 40.0 )( 353.0 , 34.0 )( 475.0 , 29.0 )( 660.0 , 30.0 )( 1063.0 , 23.0 )( 1277.0 , 18.0 )( 1541.0 , 13.0 )( 2340.0 , 18.0 )( 2815.0 , 15.0 )( 3394.0 , 10.0 )( 4140.0 , 6.0 )( 8509.0 , 3.0 )( 11410.0 , 2.0 )( 27407.0 , 2.0 )( 38369.0 , 2.0 )( 49157.0 , 2.0 )( 102811.0 , 2.0 )( 155185.0 , 2.0 )( 191034.0 , 2.0 )( 584134.0 , 2.0 )};
\addlegendentry{MC} 
%\addplot[color=amethyst,  mark=triangle*, line width = 0.5mm, dashdotted,,mark options = {scale= 1.5, solid}]
%coordinates{( 25.0 , 5.0 )( 81.0 , 8.0 )( 133.0 , 7.0 )( 172.0 , 10.0 )( 215.0 , 9.0 )( 250.0 , 6.0 )( 302.0 , 5.0 )( 391.0 , 6.0 )( 498.0 , 5.0 )( 596.0 , 5.0 )( 710.0 , 5.0 )( 891.0 , 4.0 )( 1170.0 , 4.0 )( 1436.0 , 3.0 )( 1631.0 , 3.0 )( 2785.0 , 4.0 )( 3631.0 , 2.0 )( 4295.0 , 2.0 )( 9381.0 , 2.0 )( 12823.0 , 2.0 )( 21575.0 , 2.0 )( 29622.0 , 2.0 )( 37646.0 , 2.0 )( 50775.0 , 2.0 )( 90671.0 , 2.0 )( 117003.0 , 2.0 )( 144823.0 , 2.0 )( 197654.0 , 2.0 )( 506061.0 , 2.0 )};
%\addlegendentry{BBK} 
\addplot[color=gold,  mark=pentagon*, line width = 0.5mm, dashdotted,,mark options = {scale= 1.5, solid}]
coordinates{( 25.0 , 5.0 )( 81.0 , 5.0 )( 138.0 , 5.0 )( 165.0 , 5.0 )( 206.0 , 5.0 )( 250.0 , 5.0 )( 306.0 , 5.0 )( 416.0 , 5.0 )( 526.0 , 5.0 )( 763.0 , 4.0 )( 966.0 , 4.0 )( 1280.0 , 3.0 )( 1573.0 , 2.0 )( 2349.0 , 4.0 )( 3263.0 , 2.0 )( 4433.0 , 2.0 )( 9604.0 , 2.0 )( 12715.0 , 2.0 )( 29736.0 , 2.0 )( 46034.0 , 2.0 )( 133437.0 , 2.0 )( 180029.0 , 2.0 )( 505732.0 , 2.0 )};
\addlegendentry{BBK} 
\addplot[color=carrotorange,  mark=halfcircle*, line width = 0.5mm, dashdotted,,mark options = {scale= 1.5, solid}]
coordinates{( 25.0 , 5.0 )( 81.0 , 5.0 )( 102.0 , 5.0 )( 139.0 , 5.0 )( 182.0 , 5.0 )( 215.0 , 4.0 )( 281.0 , 4.0 )( 416.0 , 4.0 )( 537.0 , 4.0 )( 665.0 , 4.0 )( 843.0 , 4.0 )( 1049.0 , 4.0 )( 1283.0 , 4.0 )( 1822.0 , 3.0 )( 2445.0 , 3.0 )( 3429.0 , 3.0 )( 4236.0 , 2.0 )( 8685.0 , 2.0 )( 12894.0 , 2.0 )( 28075.0 , 2.0 )( 38921.0 , 2.0 )( 83973.0 , 2.0 )( 107152.0 , 2.0 )( 150964.0 , 2.0 )( 181723.0 , 2.0 )( 506284.0 , 2.0 )};
\addlegendentry{MUAS} 
\addplot[color=magenta,  mark=otimes, line width = 0.5mm, dashdotted,,mark options = {scale= 1.5, solid}]
coordinates{( 25.0 , 6.0 )( 81.0 , 5.0 )( 106.0 , 5.0 )( 160.0 , 5.0 )( 210.0 , 5.0 )( 261.0 , 4.0 )( 339.0 , 4.0 )( 495.0 , 4.0 )( 667.0 , 4.0 )( 875.0 , 4.0 )( 1193.0 , 4.0 )( 1567.0 , 4.0 )( 2169.0 , 3.0 )( 3102.0 , 4.0 )( 4463.0 , 2.0 )( 6910.0 , 2.0 )( 10230.0 , 2.0 )( 13583.0 , 2.0 )( 35613.0 , 2.0 )( 46890.0 , 2.0 )( 103883.0 , 2.0 )( 174329.0 , 2.0 )( 369795.0 , 2.0 )};
\addlegendentry{SMUAS}  
\end{semilogxaxis}
\end{tikzpicture}
\hspace*{1em}
\begin{tikzpicture}[scale=0.6]
\begin{semilogxaxis}[
    legend pos=north east, xlabel = $\#\ \mathrm{dof}$, ylabel = iterations+rejections,
    legend cell align ={left}, title = {$\varepsilon=10^{-3}$, Grid~2},
    ymax= 280,
     legend style={nodes={scale=0.75, transform shape}}]
\addplot[color=crimson,  mark=square*, line width = 0.5mm, dashdotted,,mark options = {scale= 1.5, solid}]
coordinates{( 25.0 , 188.0 )( 81.0 , 152.0 )( 109.0 , 142.0 )( 137.0 , 109.0 )( 186.0 , 117.0 )( 248.0 , 116.0 )( 318.0 , 89.0 )( 401.0 , 89.0 )( 519.0 , 70.0 )( 664.0 , 221.0 )( 895.0 , 60.0 )( 1268.0 , 50.0 )( 1950.0 , 75.0 )( 2826.0 , 40.0 )( 3584.0 , 24.0 )( 7264.0 , 18.0 )( 11228.0 , 17.0 )( 20157.0 , 14.0 )( 31661.0 , 4.0 )( 44007.0 , 2.0 )( 95047.0 , 2.0 )( 147461.0 , 2.0 )( 197431.0 , 2.0 )( 339594.0 , 2.0 )};
\addlegendentry{BJK} 
\addplot[color=dark_green,  mark=diamond*, line width = 0.5mm, dashdotted,,mark options = {scale= 1.5, solid}]
coordinates{( 25.0 , 37.0 )( 81.0 , 40.0 )( 117.0 , 44.0 )( 162.0 , 43.0 )( 201.0 , 58.0 )( 268.0 , 55.0 )( 366.0 , 43.0 )( 486.0 , 43.0 )( 608.0 , 42.0 )( 816.0 , 39.0 )( 1053.0 , 30.0 )( 1366.0 , 27.0 )( 1659.0 , 30.0 )( 2624.0 , 26.0 )( 3206.0 , 20.0 )( 3706.0 , 20.0 )( 4685.0 , 14.0 )( 8011.0 , 24.0 )( 9605.0 , 13.0 )( 12517.0 , 8.0 )( 22489.0 , 11.0 )( 29651.0 , 7.0 )( 41626.0 , 2.0 )( 54688.0 , 2.0 )( 79302.0 , 6.0 )( 106273.0 , 5.0 )( 138433.0 , 4.0 )( 184010.0 , 2.0 )( 367913.0 , 2.0 )};
\addlegendentry{MC} 
%\addplot[color=amethyst,  mark=triangle*, line width = 0.5mm, dashdotted,,mark options = {scale= 1.5, solid}]
%coordinates{( 25.0 , 5.0 )( 81.0 , 8.0 )( 103.0 , 17.0 )( 123.0 , 7.0 )( 145.0 , 7.0 )( 183.0 , 10.0 )( 269.0 , 7.0 )( 354.0 , 6.0 )( 479.0 , 6.0 )( 577.0 , 6.0 )( 722.0 , 7.0 )( 901.0 , 8.0 )( 1119.0 , 6.0 )( 1336.0 , 5.0 )( 1697.0 , 4.0 )( 2171.0 , 4.0 )( 2850.0 , 4.0 )( 3950.0 , 4.0 )( 5269.0 , 3.0 )( 8010.0 , 4.0 )( 10466.0 , 3.0 )( 14204.0 , 2.0 )( 19395.0 , 3.0 )( 25490.0 , 3.0 )( 35795.0 , 2.0 )( 46381.0 , 2.0 )( 63784.0 , 2.0 )( 83793.0 , 2.0 )( 116537.0 , 2.0 )( 146638.0 , 2.0 )( 194479.0 , 2.0 )( 288709.0 , 2.0 )};
%\addlegendentry{BBK} 
\addplot[color=gold,  mark=pentagon*, line width = 0.5mm, dashdotted,,mark options = {scale= 1.5, solid}]
coordinates{( 25.0 , 5.0 )( 81.0 , 5.0 )( 103.0 , 5.0 )( 130.0 , 5.0 )( 160.0 , 5.0 )( 193.0 , 5.0 )( 238.0 , 5.0 )( 324.0 , 5.0 )( 458.0 , 5.0 )( 586.0 , 6.0 )( 748.0 , 4.0 )( 936.0 , 4.0 )( 1297.0 , 4.0 )( 1661.0 , 4.0 )( 2358.0 , 4.0 )( 3585.0 , 4.0 )( 5029.0 , 3.0 )( 7488.0 , 4.0 )( 10403.0 , 3.0 )( 13348.0 , 2.0 )( 18766.0 , 3.0 )( 25017.0 , 3.0 )( 31753.0 , 2.0 )( 49117.0 , 2.0 )( 76136.0 , 2.0 )( 97199.0 , 2.0 )( 121299.0 , 2.0 )( 162292.0 , 2.0 )( 237909.0 , 2.0 )( 300613.0 , 2.0 )};
\addlegendentry{BBK} 
\addplot[color=carrotorange,  mark=halfcircle*, line width = 0.5mm, dashdotted,,mark options = {scale= 1.5, solid}]
coordinates{( 25.0 , 5.0 )( 81.0 , 5.0 )( 98.0 , 5.0 )( 117.0 , 5.0 )( 139.0 , 5.0 )( 191.0 , 5.0 )( 250.0 , 5.0 )( 405.0 , 5.0 )( 509.0 , 4.0 )( 628.0 , 4.0 )( 805.0 , 4.0 )( 1004.0 , 4.0 )( 1212.0 , 4.0 )( 1517.0 , 4.0 )( 1996.0 , 4.0 )( 2651.0 , 4.0 )( 3354.0 , 4.0 )( 4437.0 , 4.0 )( 5960.0 , 3.0 )( 8374.0 , 3.0 )( 11967.0 , 3.0 )( 16410.0 , 2.0 )( 22873.0 , 2.0 )( 32208.0 , 2.0 )( 42671.0 , 2.0 )( 55614.0 , 2.0 )( 71021.0 , 2.0 )( 99115.0 , 2.0 )( 130994.0 , 2.0 )( 172908.0 , 2.0 )( 232524.0 , 2.0 )( 324309.0 , 2.0 )};
\addlegendentry{MUAS} 
\addplot[color=magenta,  mark=otimes, line width = 0.5mm, dashdotted,,mark options = {scale= 1.5, solid}]
coordinates{( 25.0 , 5.0 )( 81.0 , 5.0 )( 98.0 , 5.0 )( 144.0 , 5.0 )( 184.0 , 5.0 )( 233.0 , 5.0 )( 349.0 , 4.0 )( 475.0 , 4.0 )( 651.0 , 4.0 )( 828.0 , 4.0 )( 1088.0 , 4.0 )( 1494.0 , 4.0 )( 2007.0 , 4.0 )( 2764.0 , 4.0 )( 3600.0 , 4.0 )( 4715.0 , 4.0 )( 6337.0 , 3.0 )( 8935.0 , 3.0 )( 12345.0 , 3.0 )( 16648.0 , 2.0 )( 25755.0 , 3.0 )( 33395.0 , 3.0 )( 41261.0 , 2.0 )( 60405.0 , 2.0 )( 89383.0 , 2.0 )( 117360.0 , 2.0 )( 164595.0 , 2.0 )( 210479.0 , 2.0 )( 288803.0 , 2.0 )};
\addlegendentry{SMUAS} 
\end{semilogxaxis}
\end{tikzpicture}
\hspace*{1em}
\begin{tikzpicture}[scale=0.6]
\begin{semilogxaxis}[
    legend pos=north east, xlabel = $\#\ \mathrm{dof}$, ylabel = iterations+rejections,
    legend cell align ={left}, title = {$\varepsilon=10^{-3}$, Grid~3},
    ymax= 280,
     legend style={nodes={scale=0.75, transform shape}}]
\addplot[color=crimson,  mark=square*, line width = 0.5mm, dashdotted,,mark options = {scale= 1.5, solid}]
coordinates{( 25.0 , 225.0 )( 81.0 , 181.0 )( 122.0 , 124.0 )( 156.0 , 196.0 )( 196.0 , 136.0 )( 258.0 , 115.0 )( 307.0 , 114.0 )( 432.0 , 125.0 )( 564.0 , 81.0 )( 681.0 , 70.0 )( 965.0 , 202.0 )( 1201.0 , 56.0 )( 1447.0 , 67.0 )( 2586.0 , 71.0 )( 3436.0 , 30.0 )( 4795.0 , 20.0 )( 8274.0 , 19.0 )( 11449.0 , 21.0 )( 20244.0 , 8.0 )( 28484.0 , 6.0 )( 40008.0 , 2.0 )( 51737.0 , 2.0 )( 88425.0 , 2.0 )( 135342.0 , 2.0 )( 189597.0 , 2.0 )( 263573.0 , 2.0 )};
\addlegendentry{BJK} 
\addplot[color=dark_green,  mark=diamond*, line width = 0.5mm, dashdotted,,mark options = {scale= 1.5, solid}]
coordinates{( 25.0 , 39.0 )( 81.0 , 71.0 )( 117.0 , 48.0 )( 158.0 , 49.0 )( 200.0 , 50.0 )( 277.0 , 40.0 )( 362.0 , 43.0 )( 487.0 , 31.0 )( 603.0 , 38.0 )( 765.0 , 33.0 )( 999.0 , 31.0 )( 1220.0 , 31.0 )( 1518.0 , 26.0 )( 2838.0 , 28.0 )( 3371.0 , 22.0 )( 3877.0 , 15.0 )( 4754.0 , 13.0 )( 7951.0 , 18.0 )( 9828.0 , 10.0 )( 12456.0 , 8.0 )( 23093.0 , 12.0 )( 30548.0 , 7.0 )( 43723.0 , 2.0 )( 54867.0 , 2.0 )( 79573.0 , 6.0 )( 122392.0 , 5.0 )( 159471.0 , 4.0 )( 207854.0 , 2.0 )( 326961.0 , 2.0 )};
\addlegendentry{MC} 
%\addplot[color=amethyst,  mark=triangle*, line width = 0.5mm, dashdotted,,mark options = {scale= 1.5, solid}]
%coordinates{( 25.0 , 5.0 )( 81.0 , 38.0 )( 106.0 , 29.0 )( 132.0 , 11.0 )( 164.0 , 8.0 )( 195.0 , 10.0 )( 237.0 , 10.0 )( 327.0 , 9.0 )( 423.0 , 8.0 )( 526.0 , 8.0 )( 678.0 , 7.0 )( 846.0 , 7.0 )( 1081.0 , 5.0 )( 1318.0 , 5.0 )( 1762.0 , 4.0 )( 2402.0 , 4.0 )( 3280.0 , 4.0 )( 4641.0 , 3.0 )( 6150.0 , 3.0 )( 8791.0 , 4.0 )( 12095.0 , 2.0 )( 15435.0 , 3.0 )( 25481.0 , 3.0 )( 36086.0 , 2.0 )( 51763.0 , 2.0 )( 71977.0 , 2.0 )( 100287.0 , 2.0 )( 137797.0 , 2.0 )( 183672.0 , 2.0 )( 238317.0 , 2.0 )( 324958.0 , 2.0 )};
%\addlegendentry{BBK} 
\addplot[color=gold,  mark=pentagon*, line width = 0.5mm, dashdotted,,mark options = {scale= 1.5, solid}]
coordinates{( 25.0 , 5.0 )( 81.0 , 7.0 )( 97.0 , 7.0 )( 121.0 , 6.0 )( 171.0 , 6.0 )( 210.0 , 5.0 )( 256.0 , 6.0 )( 305.0 , 7.0 )( 389.0 , 5.0 )( 517.0 , 4.0 )( 667.0 , 4.0 )( 837.0 , 4.0 )( 1125.0 , 4.0 )( 1445.0 , 4.0 )( 1970.0 , 4.0 )( 2637.0 , 4.0 )( 3851.0 , 4.0 )( 4943.0 , 3.0 )( 6794.0 , 3.0 )( 9617.0 , 3.0 )( 14271.0 , 2.0 )( 19813.0 , 3.0 )( 25975.0 , 3.0 )( 36188.0 , 2.0 )( 47667.0 , 2.0 )( 70391.0 , 2.0 )( 91582.0 , 2.0 )( 128514.0 , 2.0 )( 167238.0 , 2.0 )( 236889.0 , 2.0 )( 325222.0 , 2.0 )};
\addlegendentry{BBK} 
\addplot[color=carrotorange,  mark=halfcircle*, line width = 0.5mm, dashdotted,,mark options = {scale= 1.5, solid}]
coordinates{( 25.0 , 5.0 )( 81.0 , 5.0 )( 102.0 , 5.0 )( 149.0 , 5.0 )( 226.0 , 4.0 )( 295.0 , 5.0 )( 405.0 , 5.0 )( 542.0 , 5.0 )( 693.0 , 4.0 )( 863.0 , 4.0 )( 1147.0 , 4.0 )( 1397.0 , 4.0 )( 1811.0 , 4.0 )( 2456.0 , 4.0 )( 3270.0 , 4.0 )( 4653.0 , 4.0 )( 6004.0 , 3.0 )( 7549.0 , 3.0 )( 10548.0 , 3.0 )( 13539.0 , 2.0 )( 17488.0 , 3.0 )( 23394.0 , 3.0 )( 33425.0 , 2.0 )( 47393.0 , 2.0 )( 63970.0 , 2.0 )( 90458.0 , 2.0 )( 120666.0 , 2.0 )( 161591.0 , 2.0 )( 211845.0 , 2.0 )( 288261.0 , 2.0 )};
\addlegendentry{MUAS} 
\addplot[color=magenta,  mark=otimes, line width = 0.5mm, dashdotted,,mark options = {scale= 1.5, solid}]
coordinates{( 25.0 , 5.0 )( 81.0 , 5.0 )( 102.0 , 5.0 )( 163.0 , 5.0 )( 220.0 , 5.0 )( 366.0 , 5.0 )( 488.0 , 4.0 )( 647.0 , 4.0 )( 801.0 , 4.0 )( 1048.0 , 4.0 )( 1322.0 , 4.0 )( 1889.0 , 4.0 )( 2539.0 , 4.0 )( 3472.0 , 4.0 )( 4432.0 , 4.0 )( 5743.0 , 4.0 )( 7839.0 , 3.0 )( 11717.0 , 4.0 )( 16150.0 , 2.0 )( 21255.0 , 3.0 )( 30321.0 , 3.0 )( 42331.0 , 2.0 )( 55162.0 , 2.0 )( 70952.0 , 2.0 )( 99795.0 , 2.0 )( 129697.0 , 2.0 )( 173068.0 , 2.0 )( 249514.0 , 2.0 )( 342796.0 , 2.0 )};
\addlegendentry{SMUAS} 
\end{semilogxaxis}
\end{tikzpicture}}
\caption{Example~\ref{ex:nonlinear_cdr}: Number of iterations and rejections for Grids~1–3 (left to right).}
\label{fig:iteration_nonlinear_cdr}
\end{figure}
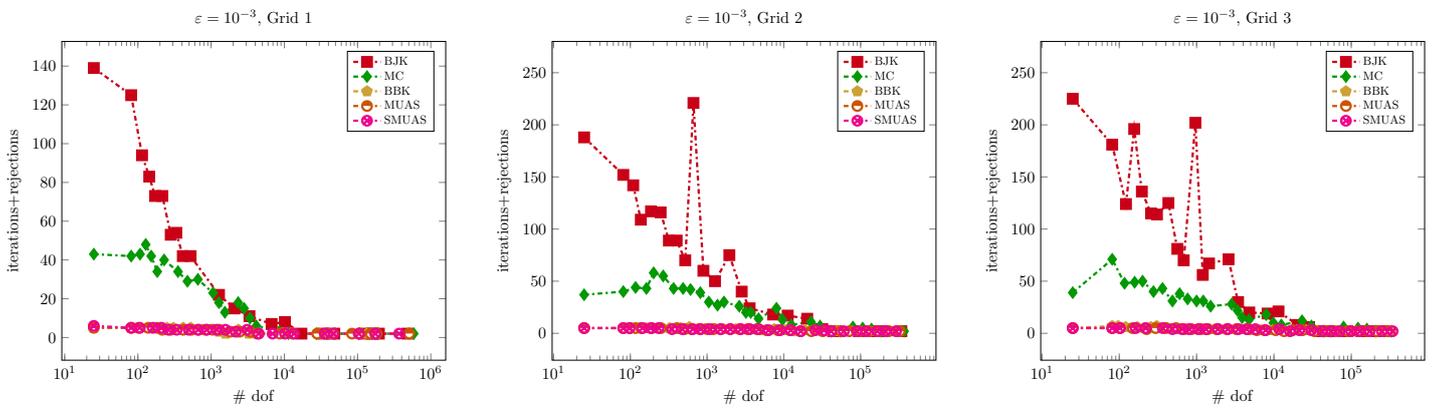

Figure~\ref{fig:iteration_nonlinear_cdr} shows the number of iterations and
rejections obtained using the \emph{fixed-point matrix} method. The BBK, MUAS, and
SMUAS limiters require the smallest number of nonlinear iterations across all
grids, whereas the BJK limiter performs the worst.  In contrast, the algebraic
stabilization methods depend on smoother indicators and therefore lead to a more stable nonlinear iteration. Consequently, algebraically stabilized schemes are
more efficient from the solver perspective for nonlinear problems.

\begin{figure}[tbp]
    \centering{
    \includegraphics[width=0.3\linewidth]{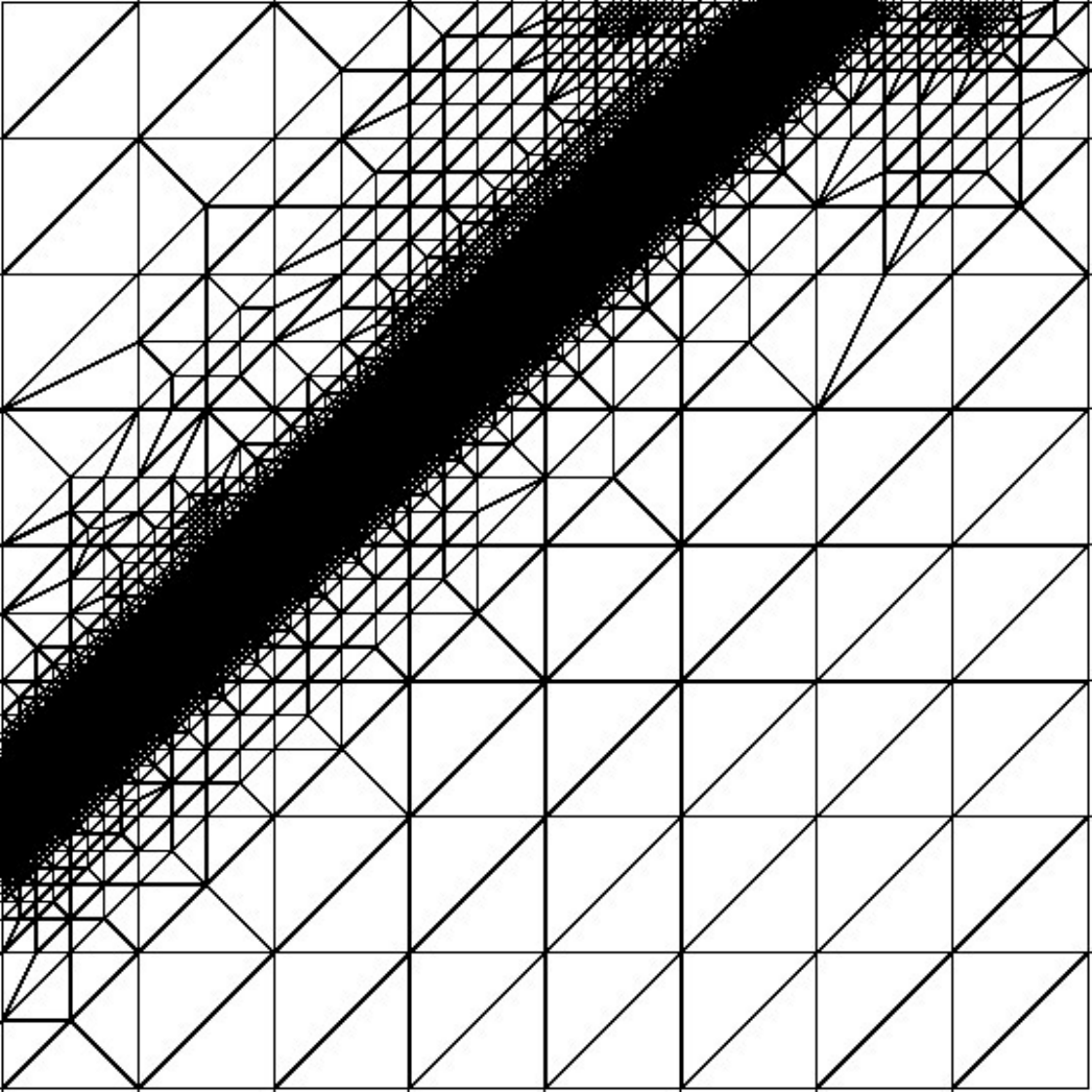}\hspace*{1em}
    \includegraphics[width=0.3\linewidth]{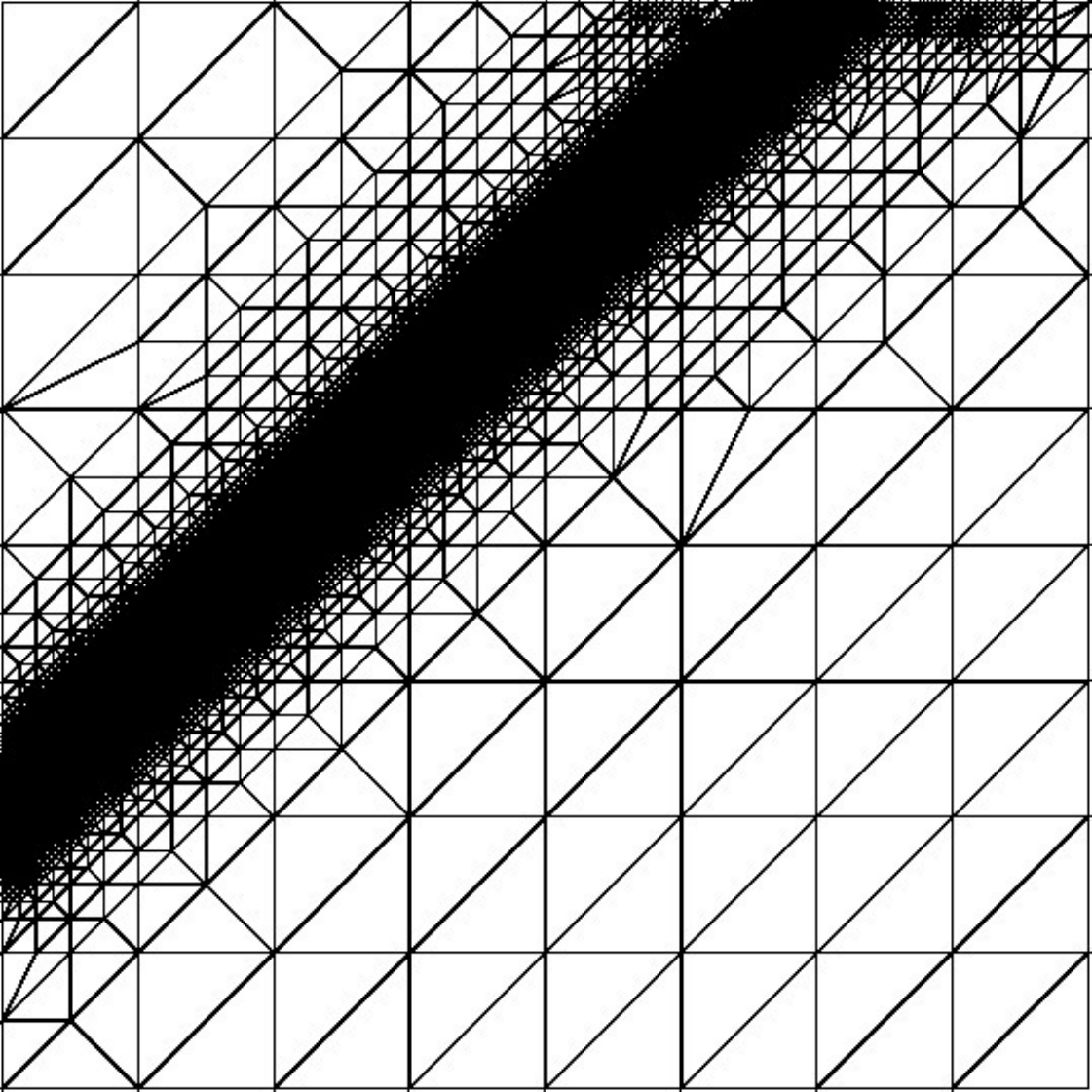}}
    \caption{Example~\ref{ex:nonlinear_cdr} Adaptive grids for the AFC methods with $\#\ \mathrm{dofs}\approx 2.5\times 10^5$. BJK limiter (left) and MC limiter (right).}
    \label{fig:nl_adaptive_afc}
\end{figure}

\begin{figure}[tbp]
    \centering{
    \includegraphics[width=0.3\linewidth]{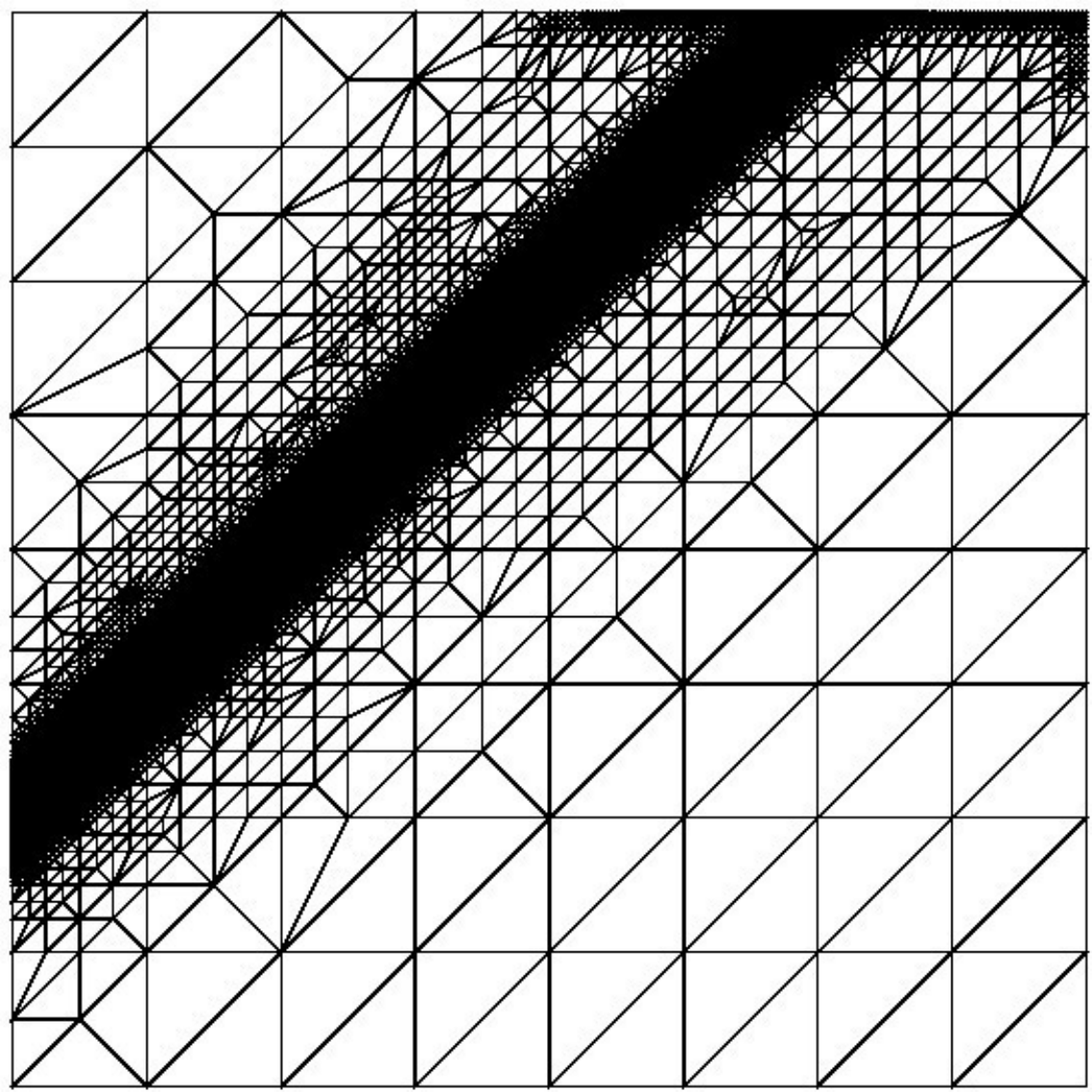}\hspace*{1em}
    \includegraphics[width=0.3\linewidth]{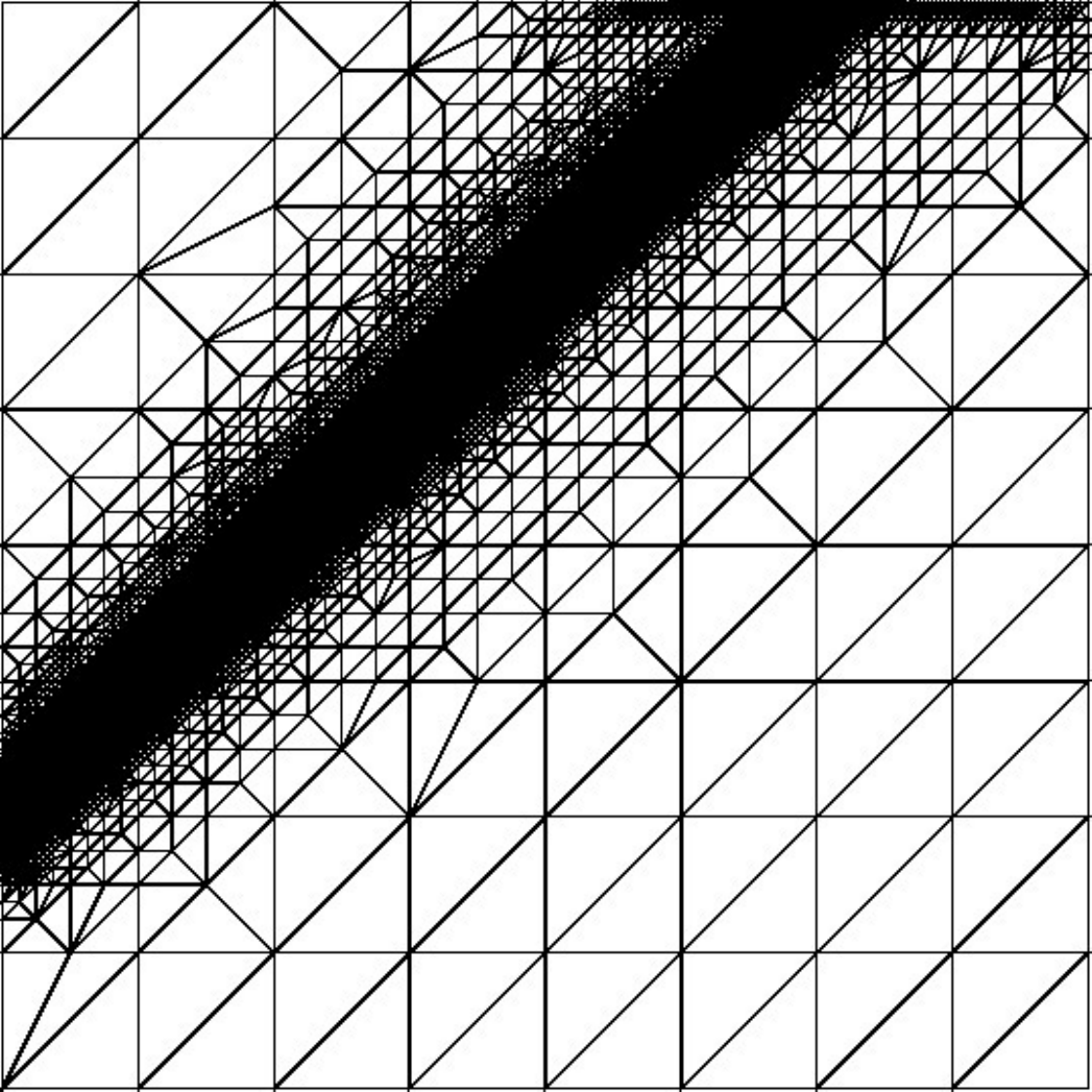}\hspace*{1em}
    \includegraphics[width=0.3\linewidth]{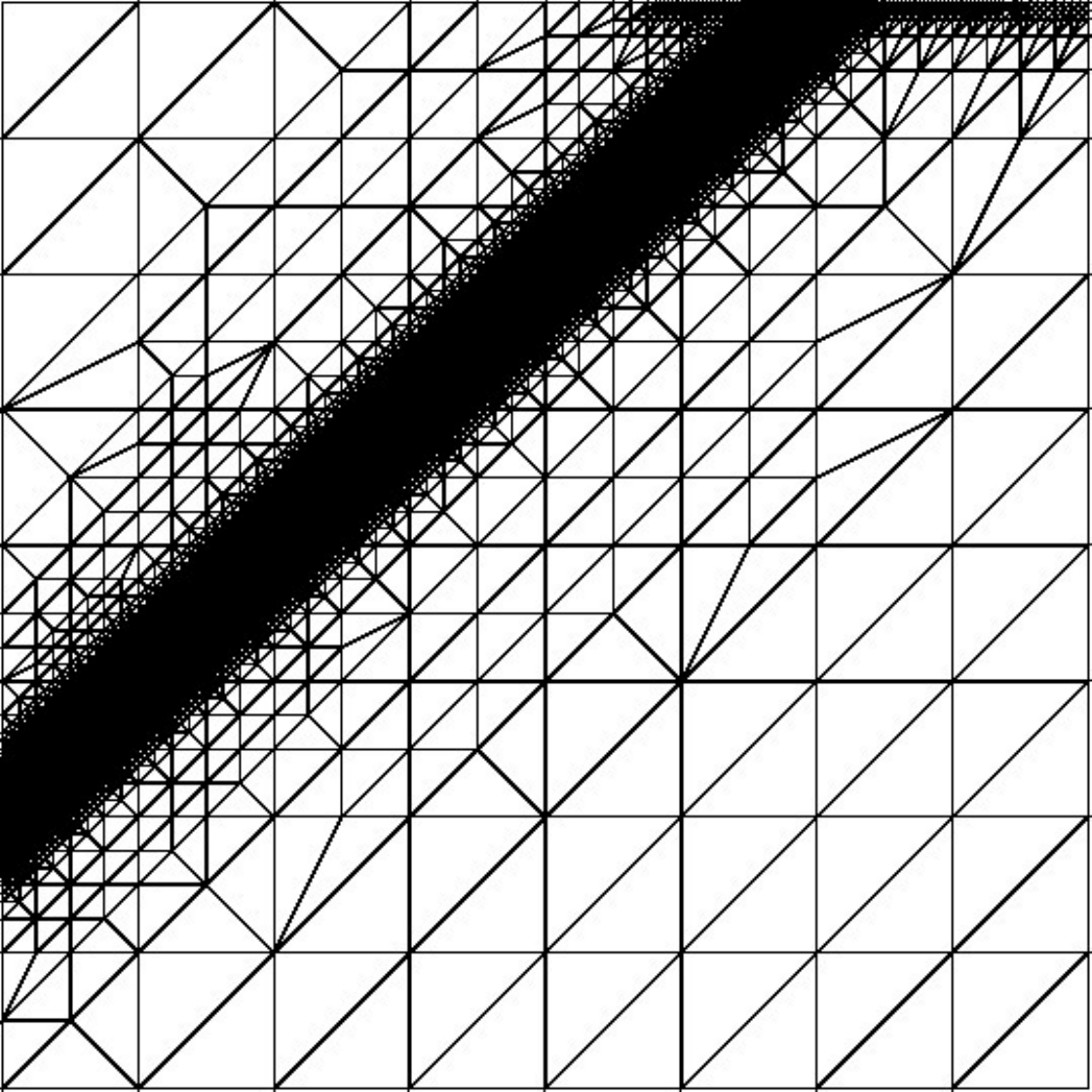}}
    \caption{Example~\ref{ex:nonlinear_cdr} Adaptive grids for the algebraic stabilized methods with $\#\ \mathrm{dofs}\approx 2.5\times 10^5$. MUAS (left) , SMUAS (middle),  BBK (right).}
    \label{fig:nl_adaptive_as}
\end{figure}

Figs.~\ref{fig:nl_adaptive_afc} and~\ref{fig:nl_adaptive_as} show the adaptively
refined grids for the AFC and algebraically stabilized methods, respectively.
Among the AFC approaches, the BJK limiter resolves the interior layer most sharply.
Among the algebraically stabilized methods, the BBK limiter appears most effective,
although all methods exhibit some degree of over-refinement near the top-right
corner where the layer terminates.

Finally, the computing times are reported in Table~\ref{tab:nl_time}. The BBK
limiter is the most efficient overall, while the MC and MUAS limiters also show
good performance. In contrast, the SMUAS method is significantly more expensive,
requiring approximately twenty-six times more computational time than the next
slowest method. This high cost arises from the use of reconstructed states
involving fictitious points and gradient evaluations, which must be performed
repeatedly as the mesh is refined. Overall, the BBK method provides the best
balance between accuracy and computational efficiency for this nonlinear example.

\begin{table}[tbp]
\centering
\begin{tabular}{|c|c|c|}
\hline
\textbf{Limiter} & \textbf{Time (s)} & \textbf{Number of Adaptive Loops}\\
\hline \hline
BJK & 190 & 22 \\
MC & 205 & 26\\
BBK & 160 & 23\\
MUAS & 193 & 26\\
SMUAS & 5025 & 24\\
\hline
\end{tabular}
\caption{Example~\ref{ex:nonlinear_cdr}: Computing time for the adaptive algorithm (in seconds) and the number of adaptive loops for Grid~1.}
\label{tab:nl_time}
\end{table}

\section{Summary}\label{sec:summary}
In this section, we summarize the results presented throughout the paper.
We investigated the numerical behaviour of several limiting strategies on
adaptively refined meshes. In particular, we considered AFC schemes with
the BJK and MC limiters and algebraically stabilized schemes, including
MUAS, SMUAS, and the BBK method. The numerical assessment is organized
into three categories.

\begin{enumerate}

\item \textbf{Accuracy:}
We evaluate accuracy using both the effectivity index and the error in the
numerical solution.

From Example~\ref{ex:boundary_layer}, all methods achieve optimal
convergence rates. The BJK, BBK, and SMUAS methods produce the smallest
effectivity indices, while differences in the actual error remain minor.

In the nonlinear Example~\ref{ex:nonlinear_cdr}, the behaviour changes.
The transport direction depends on the solution, and the interior layer is
not aligned with the mesh on Grids~2 and~3. In this case, the BJK limiter
retains the most accurate solution, while the remaining methods experience
a loss of optimal convergence. The MC limiter also performs robustly,
whereas the algebraically stabilized methods yield slightly larger errors.
Overall, the BJK limiter provides the most reliable accuracy across all
tests.

\item \textbf{Approximation:}
Here, we consider both the sharpness of internal layers and the quality of
adaptive mesh refinement.

In Example~\ref{ex:corner_singularity}, the SMUAS method and the MC limiter
produce the most effective adaptive meshes. In
Example~\ref{ex:volker_example}, all methods behave similarly, although the
SMUAS and BJK limiters resolve the layer slightly more sharply.
For the Hemker problem (Example~\ref{ex:hemker}), the MUAS and BBK methods
lead to well-distributed adaptive meshes, while the SMUAS method shows
increased smearing. The BJK limiter, on the other hand, produces a less
favorable adaptive grid, but still captures the layer sharply in the
numerical solution.

Taken together, these results indicate that no single method is uniformly
superior in approximation behaviour. The BJK limiter tends to produce the
sharpest layers and therefore highly accurate solutions, whereas the MC,
MUAS, and BBK methods interact more favorably with the residual-based
estimator and therefore generate more effective adaptive meshes. Thus, the
differences between the methods primarily reflect the interaction between
stabilization and mesh refinement rather than a clear superiority of one
method in approximation quality.

\item \textbf{Efficiency:}
Efficiency is measured using both the number of nonlinear iterations and
the overall computational time.

With respect to nonlinear iterations, the BJK limiter requires the largest
number of iterations, whereas the MUAS, SMUAS, and BBK methods converge
more rapidly. However, computational time tells a different story. The SMUAS
method is significantly more expensive due to the gradient reconstruction
used in the limiter, while the BBK method is the fastest. The MUAS method
also performs efficiently and provides a good balance between iteration
count and runtime.
\end{enumerate}

Table~\ref{tab:summary} summarizes these observations. The symbols $++$,
$+$, and $\circ$ indicate strong, moderate, and comparatively weaker
performance, respectively. Overall, the results suggest that the BJK limiter
is the most accurate method, the BBK and MUAS methods are the most efficient,
and the MC limiter provides consistently good adaptive approximation.

\begin{table}[tbp]
\centering
\begin{tabular}{|c|c|c|c|}
\hline
\textbf{Limiter} & \textbf{Accuracy} & \textbf{Approximation} & \textbf{Efficiency} \\
\hline\hline
\textbf{BJK} & $++$ & $+$ & $\circ$ \\
\textbf{MC} & $+$ & $+$ & $+$ \\
\textbf{BBK} & $+$ & $+$ & $++$ \\
\textbf{MUAS} & $+$ & $+$ & $++$ \\
\textbf{SMUAS} & $\circ$ & $+$ & $\circ$ \\
\hline
\end{tabular}
\caption{Summary of the numerical performance of the considered limiters.}
\label{tab:summary}
\end{table}

The numerical study also reveals several open questions. In particular, the
behaviour of the SMUAS limiter in the Hemker problem indicates that discrete
maximum principle properties may depend not only on the mesh but also on the
alignment of sharp layers with respect to the grid. Furthermore, the nonlinear
example with solution–dependent convection demonstrates that limiter performance
changes when the transport direction evolves during the nonlinear iteration.
A theoretical understanding of algebraic stabilization schemes in such
quasilinear settings, as well as the interaction between nonlinear stabilization
and residual-based a posteriori estimators, remains an interesting topic for
future research.

\noindent{\bf Acknowledgement}. 
The authors would like to thank Prof. Dr. Volker John for discussions regarding the implementation of the SMUAS method.  N.A. acknowledges support from the Gulf University for Science and Technology through grant
 ISG Case No. 96. A.J. acknowledges support from the Indian Institute of Technology  Gandhinagar through grant IP/52016.
\bibliographystyle{plain}
\bibliography{Posteriori_MUAS}
\end{document}